\documentclass[a4paper,12pt]{amsart}
\usepackage{import}
\usepackage[utf8]{inputenc}
\usepackage{datetime} 
\usepackage[margin=1in]{geometry}
\usepackage{amsmath, amssymb, amsthm, extarrows, float,mathrsfs, mathtools, amsfonts}
\usepackage{color}
\usepackage{tikz, tikz-cd}
\usepackage{quiver}
\usepackage{graphicx}
\usepackage{tensor}
\usepackage[colorinlistoftodos]{todonotes}
\usepackage{latexsym}
\usepackage{mdwlist}
\usepackage[backend=biber,
            isbn=false,
            doi=false,
            eprint=false,
            url=false,
            maxbibnames=9,
            giveninits=true,
            style=alphabetic,
            citestyle=alphabetic]{biblatex}
\AtEveryBibitem{\clearlist{language}}
\bibliography{References}
\renewbibmacro{in:}{}

\usepackage{color}
\definecolor{e-mail}{rgb}{0,.40,.80}
\definecolor{burgundy}{RGB}{128,0,32}
\definecolor{citation}{rgb}{0,.40,.80}

\usepackage[colorlinks=true,
            linkcolor=burgundy,
            citecolor=citation,
            urlcolor=e-mail]{hyperref}
\usepackage{preamble}
\title{Exponential map in DT theory }
\author{\v{S}arūnas Kaubrys}
\address{KAVLI IPMU (WPI), UTIAS, THE UNIVERSITY OF TOKYO,
KASHIWA, CHIBA 277-8583, JAPAN}
\email{sarunas.kaubrys@ipmu.jp}
\begin{document}

\begin{abstract}
    This paper studies the Cohomological Donaldson-Thomas theory of loop stacks of $0$-shifted symplectic stacks. In particular, we compare $(-1)$-shifted tangent stacks of these moduli problems, which we view as additive, to loop stacks, which we view as multiplicative, via an exponential map that preserves induced $(-1)$-shifted symplectic structures. As an application, we prove for certain moduli of objects of $2$-Calabi-Yau categories a loop dimensional reduction theorem for the loop stacks of these moduli spaces. Finally, we prove a loop version of nonabelian Hodge theory for stacks in the $\GL_n$ case.
\end{abstract}
\maketitle

\setcounter{tocdepth}{1}
\tableofcontents
\newpage

\section{Introduction}
Cohomological Donaldson-Thomas theory is a categorification of sheaf counting invariants on Calabi-Yau $3$-folds, originally defined by Thomas \cite{thomas_dt_og}. In \cite{ben2015darboux} these invariants were generalised to oriented $(-1)$-shifted symplectic stacks in the sense of \cite{PTVV}. In particular, the authors define a canonical perverse sheaf on $X$ called the DT sheaf $\varphi_{X}$. The study of the cohomology of this perverse sheaf is known as cohomological DT theory. In this paper, we will compare additive and multiplicative cohomological DT invariants by using a type of exponential map. As an application to nonabelian Hodge theory, we will prove a multiplicative version of nonabelian Hodge theory for stacks as in \cite{davison2022bps}.
\subsection{Dimensional reduction}
In certain cases, cohomological DT invariants recover known cohomology theories. Consider a quasismooth derived artin stack $X$. By analogy to classical symplectic geometry in \cite{calaque_cot} it was proven that the shifted cotangent stack $\To^{*}[-1]X$ is canonically $(-1)$-shifted symplectic. Furthermore, it can be shown that this stack is oriented and therefore, we can define $\varphi_{\To^{*}[-1]X}$ .  In \cite{dim_red_kinjo}, it was proven that there is an equivalence
\begin{equation} \label{dim_red_intro}
    \pi_{*}\varphi_{\To^{*}[-1]X} = \omega_{X}[- \mathrm{vdim}X],
\end{equation}
with $\pi \colon \To^{*}[-1]X \to X$ the projection and $\omega_{X}$ the dualizing sheaf of $X$. In particular, by taking cohomology of both sides we get 
\begin{equation} \label{intro_dim_redeq}
    \mathrm{H}^{*}(X,\varphi_{\To^{*}[-1]X}) = \mathrm{H}^{\mathrm{BM}}_{ \mathrm{vdim}-*}(X).
\end{equation}
Therefore, we can compute the Cohomological DT invariants of $\To^{*}[-1]Y$ in terms of the Borel-Moore homology of the base. This process is known as dimensional reduction. This refers to the fact that a natural source of $(-1)$-shifted symplectic stacks comes from coherent sheaves on $3$-Calabi-Yau varieties or more generally $3$-Calabi-Yau categories. Denote by $\mathrm{Coh}X$ the derived stack of coherent sheaves on a variety. If we consider $\To^*[-1] \mathrm{Coh}S$ for a surface $S$, we get the isomorphism  
\begin{equation}
    \To^{*}[-1]\mathrm{Coh}S  \cong \mathrm{Coh}_{c}(\mathrm{Tot}_{S}(K_{S})).
\end{equation}
Dimensional reduction relates the DT invariants of the $3$-fold $\mathrm{Tot}_{S}(K_{S})$ to the Borel-Moore homology of the stack of coherent sheaves on the base surface $S$.
Therefore, dimensional reduction relates $3$-dimensional moduli problems to $2$-dimensional moduli problems.
\subsection{Additive versus multiplicative}
Consider now the case of a $K3$ surface $S$, then in this case 
\begin{equation}\label{intro_k3surface}
    \To^{*}[-1] \mathrm{Coh}S \cong \mathrm{Coh}_{c}(S \times \mathbb{A}^{1}).
\end{equation}
The appearance of the additive group $\mathbb{A}^{1}$ suggests that we view this moduli space as an  ``additive'' moduli space. In fact, it can be shown that for any algebraic variety $X$, if we take the $(-1)$-shifted \emph{tangent}, we get
\begin{equation} 
    \To[-1] \mathrm{Coh}X \cong \mathrm{Coh}_{c}(X \times \mathbb{A}^{1}).
\end{equation}
 Furthermore, the stack $\To[-1]Y$ can be equivalently viewed as a mapping stack $\Map(\B \widehat{\mathbb{G}}_{a},Y)$ and a point in $\To[-1]X$ as given by pairs $(y, a \in \mathfrak{g}_{y})$ where $\mathfrak{g}_{y} = \mathrm{Lie}(G_{y})$.
Now consider $Y$ equipped with a $0$-shifted symplectic structure. Then the isomorphism
\begin{equation}
    \mathbb{T}_{Y} \cong \mathbb{L}_{Y}
\end{equation}
induces an isomorphism of stacks
\begin{equation}
    \To^{*}[-1]Y \cong \To[-1]Y.
\end{equation}
This equivalence recovers equation \eqref{intro_k3surface} mentioned before.
Therefore, in this paper we will rather study the stack $\To[-1]Y$.
We can also consider the ``multiplicative'' version, using the loop stack $\mathcal{L}X = \Map(S^{1},X)$. This is a derived version of the inertia stack used in classical algebraic geometry. Again, a point in $\mathcal{L}Y$ is a pair $(y, A \in G_{y})$ where $G_{y}$ is the stabilizer of the point $y$. Similarly, if $X$ is an algebraic variety then 
\begin{equation}
    \mathcal{L} \mathrm{Coh}X = \mathrm{Coh}_{c}(X \times \mathbb{G}_{m}).
\end{equation}
 Using the AKSZ construction in \cite{PTVV}, for any $0$-shifted symplectic $Y$, we can construct a canonical $(-1)$-shifted symplectic structure on $\mathcal{L}Y$. By work of \cite{naef2023torsion} both $\To[-1]Y$ and $\mathcal{L}Y$ are oriented so we can define the DT sheaf and study cohomological DT invariants of these moduli spaces. In this paper we will consider the following questions
\begin{enumerate}
    \item \textbf{Question:} For $X$ a $0$-shifted symplectic stack can we compare the cohomological DT invariants of $\To[-1]X$ and $\mathcal{L}X$?
    \item \textbf{Question:} \text{Can we formulate a similar dimensional reduction theorem for} $\mathcal{L}X$ ?
\end{enumerate}
\subsection{Exponential map} \label{exp_map_general_intro}
We will compare the cohomological DT invariants of $\To[-1]X$ and $\mathcal{L}X$ using certain exponential maps. For derived Artin stacks one can define a formal exponential map using the work of \cite{ben2012loop}
\begin{equation}
    \exp \colon \widehat{\To}^{0_{X}}[-1]X \to \widehat{\mathcal{L}}^{\mathrm{const}}X
\end{equation}
where $\mathrm{const} \colon X \to \mathcal{L}X$ is the map of constant loops and $0_{X} \colon X \to \To[-1]X$ is the zero section. We prove that if we start with an $n$-shifted symplectic stack $X$, then this map is compatible with the induced $(n-1)$-shifted symplectic structures.
\begin{thm}[= Theorem \ref{exp_general_thm}] \label{exp_general_intro}
    Let $X$ be a derived artin stack equipped with an $n$-shifted symplectic structure. The exponential map preserves $(n-1)$-shifted symplectic structures on both sides.
\end{thm}
The proof is an application of the AKSZ formalism applied to the correspondence
\begin{equation} 
\begin{tikzcd}
	{\widehat{\To}^{0_{X}}[-1]X} && {\widehat{\mathcal{L}}^{\mathrm{const}}X} \\
	& {\widehat{\mathcal{L}}^{u,\mathrm{const}}X} \\
	& {\mathcal{L}^{u}X} \\
	{\To[-1]X} && {\mathcal{L}X}
	\arrow["\exp"', from=1-1, to=1-3]
	\arrow["\exp"', from=1-1, to=2-2]
	\arrow[from=1-1, to=4-1]
	\arrow[shift right=2, from=1-3, to=2-2]
	\arrow[from=1-3, to=4-3]
	\arrow["\cong"{description}, shift right=2, from=2-2, to=1-3]
	\arrow[from=2-2, to=3-2]
	\arrow["{q_{a}}", from=3-2, to=4-1]
	\arrow["{q_{m}}"', from=3-2, to=4-3]
\end{tikzcd}
\end{equation}
Indeed, we view all three spaces 
\begin{align*}
    \To[-1]X &\cong \Map(\B \widehat{\mathbb{G}}_{a},X) \\
    \mathcal{L}X &\cong \Map(S^{1},X) \\
    \mathcal{L}^{u}X &\cong \Map(\B \mathbb{G}_{a},X)
\end{align*}
as mapping stacks from $3$-different types of ``circles''. Where $S^{1}$ is viewed as the multiplicative circle and $\B \widehat{\mathbb{G}}_{a}$ is viewed as the additive circle.
See \cite{scherotzke2025fourier} for more on this perspective.
The compatibility with symplectic structures follows from the AKSZ formalism applied to the cospan
\begin{equation}
    \B \widehat{\mathbb{G}}_{a} \to \B \mathbb{G}_{a} \xleftarrow{} S^{1}.
\end{equation}
\subsection{DT sheaves and exponential map}
The above exponential map is formal and only gives us information on certain parts of our stacks. For applications in DT theory, in special cases, we define a global complex analytic exponential map. Central to our approach to study the exponential map is the interplay between the algebraic, complex analytic and formal settings, both derived and classical. We will now give the moduli problems we will consider. These are all examples of substacks of moduli of objects in $2$-Calabi-Yau categories. 
\begin{ex}[Moduli problems we consider]\label{intro_list}
We construct an exponential map for the following moduli problems $\mathfrak{M}$
\begin{enumerate}
    \item (twisted) Local systems on a surface 
    \item semistable Higgs sheaves/bundles on a curve 
    \item (deformed) preprojective algebra representations of a quiver
    \item semistable coherent sheaves on a smooth projective symplectic surface.
\end{enumerate}
\end{ex}
Furthermore, restricting to a certain \'etale locus $\To^{\et}[-1]\mathfrak{M}$ the map $\exp$ is \'etale and surjective. This culminates in the following theorem
\begin{thm}[=Theorem \ref{main_dcrit_thm}]
    For examples of $\mathfrak{M}$ $(1)$ to $(3)$ above,  the map
    \begin{equation}
        \exp \colon \To^{\et}[-1]\mathfrak{M} \to \mathcal{L}\mathfrak{M}
    \end{equation}
    \begin{enumerate}
        \item is \'etale and surjective
        \item preserves d-critical structures
        \item preserves the canonical orientations.
    \end{enumerate}
Therefore, we obtain the following relation for DT sheaves
\begin{equation}
    \exp^{*} \varphi_{\mathcal{L}\mathfrak{M}} \cong \varphi_{\To^{\et}[-1]\mathfrak{M}}.
\end{equation}
In particular, we obtain a canonical map
\begin{equation}
    \mathrm{H}^{*}(\mathcal{L}\mathfrak{M}, \varphi_{\mathcal{L}\mathfrak{M}}) \xrightarrow{\exp} \mathrm{H}^{*}(\To^{\et}[-1]\mathfrak{M}, \varphi_{\To[-1]\mathfrak{M}}).
\end{equation}
\end{thm}
For coherent sheaves on symplectic surfaces all the results above also hold, except preservation of orientations. This means that for coherent sheaves we can relate the DT sheaf of the loop stack to a twisted version of the DT sheaf of the shifted tangent. See Remark \ref{coherent_sheaves_orient} for some discussion. \par 
The preservation of d-critical structures and orientations is proved by using the compatibility of shifted symplectic structures of the formal exponential map of Theorem \ref{exp_general_intro} and compatibility with direct sum or more generally graded points.
\subsection{Loop dimensional reduction}
Consider now $\mathfrak{M}$ as in \ref{intro_list} with good moduli space $M$ and denote by $\pi_{a} \colon \To[-1]\mathfrak{M} \to M_{a}$ and $\pi_{m} \colon \mathcal{L}\mathfrak{M} \to M_{m}$ the respective good moduli space maps. In particular, we can also define an exponential on good moduli spaces
\begin{equation}
\begin{tikzcd}
	{M^{\et}_{a}} & {M_{a}} & {M_{m}} \\
	& {M \times \mathbb{G}_{a}} & {M \times \mathbb{G}_{m}}
	\arrow[from=1-1, to=1-2]
	\arrow["\exp", from=1-2, to=1-3]
	\arrow["{\iota_{a}}"', from=2-2, to=1-2]
	\arrow["{\id \times \exp}", from=2-2, to=2-3]
	\arrow["{\iota_m}"', from=2-3, to=1-3]
\end{tikzcd}
\end{equation}
where $M^{\et}_{a}$ is a complex analytic open inside $M_{a}$ and $\iota_a$ is the inclusion of an object with a central endomorphism and $\iota_{m}$ with a central automorphism.
Then via the work of \cite{bu2025cohomology} we can define the BPS sheaves 
\begin{align*}
   \BPS_{a} = \pH^{1}(\pi_{a*} \varphi_{\To[-1]\mathfrak{M}}) \\
   \BPS_{m} = \pH^{1}(\pi_{m*}\varphi_{\mathcal{L}\mathfrak{M}}).
\end{align*}
by taking $1$st perverse cohomology of pushforwards to the good moduli space.
 The results of \cite{bu2025cohomology} give cohomological integrality  decompositions, under which the full DT cohomology of $\To[-1]\mathfrak{M}$ and $\mathcal{L}\mathfrak{M}$ is generated by cohomology of the BPS sheaves. \footnote{For the exact form of the decompositions we will use, see Section \ref{loop_naht_proof}.} 
Furthermore, we have the following support lemma \cite[Theorem 7.2.15]{bu2025cohomology} for $\BPS_{a}$
\begin{equation}
    \BPS_{a} \cong \iota_{a*}(\BPS_{M} \boxtimes \mathrm{IC}(\mathbb{G}_{a}))
\end{equation}
Therefore, instead of trying to find a direct analogue of equation \eqref{intro_dim_redeq} we instead prove a dimensional reduction or support lemma for BPS sheaves of $\mathcal{L}\mathfrak{M}$. 
\begin{thm}[Loop dimensional reduction = Theorem \ref{loop_dim_red}] \label{loop_dim_red_intro}
Let $\mathfrak{M}$ be as in Example \ref{intro_list} (1-3). \footnote{Again for coherent sheaves we have a caveat and we obtain a twisted version of this dimensional reduction. See Remark \ref{coh_warning}.}Then we have
\begin{enumerate}
    \item $\exp^{*} \BPS_{m} \cong \BPS_{a}|_{M^{\et}_{a}}$
    \item an isomorphism of perverse sheaves
    \begin{equation}
    \BPS_{M_m} \cong \iota_{m*}(\BPS_{M} \boxtimes \IC(\mathbb{G}_{m})).
\end{equation}
\end{enumerate}
\end{thm}
\subsection{ Loop nonabelian Hodge theory}
We give an application of the loop dimensional reduction. In \cite{davison2022bps} the authors prove a nonabelian Hodge theory theorem for Borel-Moore homology on stacks.
\begin{thm}
Let $g \geq 0, r \geq 1, d \in  \mathbb{Z}$. There is an equivalence of graded vector spaces
\begin{equation}
    \mathrm{H}^{\mathrm{BM}}_{*}( \mathrm{Higgs}^{\mathrm{ss}}_{g,r,d}) \cong \mathrm{H}^{\mathrm{BM}}_{*}( \mathrm{Loc}^{\omega_d}_{\GL_r}(\Sigma_g)
\end{equation}
of the stack of rank $r$, degree $d$ semistable Higgs bundles $\mathrm{Higgs}^{\mathrm{ss}}_{g,r,d}$ and rank $r$, $\omega_{d} = \exp(\frac{2 \pi i d}{r})$ twisted local systems $ \mathrm{Loc}^{\omega_d}_{\GL_r}(\Sigma_{g})$.
\end{thm}
In this paper we prove a loop version of this theorem using the cohomological integrality isomorphism and Theorem \ref{loop_dim_red}.
\begin{thm}[Loop nonabelian Hodge theory = Theorem \ref{loop_naht}] \label{loop_naht_intro}
Let $g \geq 1$, $r \geq 1$, $d \in \mathbb{Z}$. Denote by $\BPSo^{\mathrm{Dol},m}_{r,d}$ the BPS cohomology of $\mathcal{L}\mathrm{Higgs}^{\mathrm{ss}}_{g,r,d}$ and by $\BPSo^{\mathrm{B},m}_{r,d}$  the BPS cohomology of $\mathcal{L}\Loc^{\omega_{d}}_{\GL_r}(\Sigma_g)$.
We have an equivalence of BPS cohomology
\begin{equation}
     \BPSo^{\mathrm{Dol},m}_{r,d} \cong \BPSo^{\mathrm{B},m}_{r,d}.
\end{equation}
Furthermore, we have an isomorphism of graded vector spaces
\begin{equation}
    \mathrm{H}^{*}(\mathcal{L}\mathrm{Higgs}^{\mathrm{ss}}_{g,r,d}, \varphi_{\mathrm{Dol}}) \cong \mathrm{H}^{*}(\mathcal{L}\Loc^{\omega_{d}}_{\GL_r}(\Sigma_g), \varphi_{\mathrm{B}}).
\end{equation}
\end{thm}
One motivation for studying this theorem comes from the relation between Langlands duality of $3$-manifolds and topological mirror symmetry for Higgs bundles as in \cite[Section 10.3.17]{bu2025cohomology}. In particular for a semisimple group $G$ we have the following expectation of equivalences
\begin{equation}
\begin{tikzcd}
	{\mathrm{H}^{*}(\mathcal{L}\mathrm{Higgs}^{\mathrm{ss}}_{g,G}, \varphi_{\mathrm{Dol},G})} & {\mathrm{H}^{*}(\mathcal{L}\mathrm{Loc}_{G}(\Sigma_{g}), \varphi_{\mathrm{B},G})} \\
	{\mathrm{H}^{*}(\mathcal{L}\mathrm{Higgs}^{\mathrm{ss}}_{g,G^{\vee}}, \varphi_{\mathrm{Dol},G^{\vee}})} & {\mathrm{H}^{*}(\mathcal{L}\mathrm{Loc}_{G^{\vee}}(\Sigma_{g}), \varphi_{\mathrm{B},G^{\vee}})}
	\arrow["{\text{LNAHT}}", tail reversed, from=1-1, to=1-2]
	\arrow["{\text{TMS}}"', tail reversed, from=1-1, to=2-1]
	\arrow["{\text{3d Langlands}}", tail reversed, from=1-2, to=2-2]
	\arrow["{\text{LNAHT}}"', tail reversed, from=2-1, to=2-2]
\end{tikzcd}
\end{equation}
So a loop nonabelian Hodge type theorem gives a way to connect topological mirror symmetry with Langlands for DT invariants of $3$-manifolds. In particular, one can show naive $2d$ version of Langlands duality between Borel-Moore homology of $\Loc_{G}( \Sigma_g)$ or $\mathrm{Higgs}^{ss}_{g,G}$ fails.
\subsection{BPS cohomology of $\Sigma_g \times S^{1}$ in terms of twisted character varieties}
We can use the above loop dimensional reduction to compute the BPS cohomology of the stack $\Loc_{\GL_r}(\Sigma_{g} \times S^{1}) = \mathcal{L} \Loc_{\GL_r}(\Sigma_{g})$.  In particular, we can use nonabelian Hodge theory for stacks in \cite{davison2022bps} and $\chi$-independence for Higgs bundles as in \cite{kinjokoseki2024cohomological} to compute
\begin{prop}[= Proposition \ref{sigma_g_bps}]
Let $g \geq 2$
    \begin{equation}
    \BPSo^{\B,m}_{g,r,0} =  \HHf(\mathrm{M}^{\mathrm{B}}_{g,r,1}, \mathbb{Q}_{\mathrm{vir}}) \otimes \mathrm{H}^{*}(\mathbb{G}_{m}). 
\end{equation}
\end{prop}
The cohomology of twisted character varieties has been computed in \cite{hausel_mixed}. The BPS cohomology for $g=1$ has been computed in \cite{kaubrys2024cohomological}. In future work, the author hopes to use this computation to obtain computations of Skein module dimensions to obtain further evidence of \cite[Conjecture C]{gunningham2023deformation}.
\subsection{Future work}
In this paper we have mostly restricted to concrete applications for specific moduli problems over $\GL_n$. This is ultimately due to the fact that the exponential map for other groups is not surjective, which means we cannot immediately produce a cover $\To[-1]X \to \mathcal{L}X$. Furthermore, as far as the author is aware there is no theory of good moduli spaces for complex analytic stacks, which makes some arguments more difficult. The author hopes to give a unified approach to the exponential map in derived complex analytic geometry in the future with more general applications to DT theory. In particular, the results of this paper should be strengthened to handle general moduli of objects of $2$-Calabi-Yau categories with good moduli theory.
\subsection{Relation to other work}
While this paper was being finished a general multiplicative dimensional reduction theorem for BPS sheaves on loop stacks was proven in \cite{kinjo_mult}. In particular, this implies Theorem \ref{loop_dim_red_intro}. However, our method is different since we work with the complex analytic exponential map and directly compare the DT sheaves on $\To[-1]X$ and $\mathcal{L}X$ using it in the cases we consider.
\subsection*{Acknowledgements}
This project was started during my PhD thesis and I want to thank my advisers Ben Davison and Pavel Safronov for several discussions. I would also like to thank Sebastian Schlegel Mejia, Pierre Descombes, Elsa Maneval, Eric Chen and Tanguy Vernet for useful discussions. Special thanks are due to Mauro Porta for help with derived analytic geometry. The author was
supported by the Carnegie Trust for the Universities of Scotland for part of the duration of this
research. This work was also supported by World Premier International Research Center Initiative (WPI), MEXT, Japan.
\section{Derived algebraic geometry and shifted symplectic structures}
In this paper we work with derived prestacks $\operatorname{dPreStk}$ over $\mathbb{C}$. 
Namely functors $\operatorname{cdga}^{\leq 0,\op} \to \operatorname{Spc}$. Here $\operatorname{Spc}$ is the $\infty$-category of spaces. There is a truncation functor $\tt_{0}\colon \operatorname{dPreStk} \to \operatorname{PreStk}$ to classical (higher) prestacks, which has a fully faithful right adjoint $\iota\colon \operatorname{PreStk} \to \operatorname{dPreStk}$. We can further truncate to $1$-prestacks $\operatorname{PreStk}^{\leq 1}$, where $\operatorname{PreStk}^{\leq 1}$ are functors $\operatorname{CAlg}^{\op}_{\mathbb{C}} \to \operatorname{Gpd}$. We will use the $\infty$-category of quasicoherent sheaves $\QCoh(X) = \lim_{\Spec R \to X} (\Mod R)$ and the subcategory of perfect complexes $\Perf (X)$. There is an internal Hom in $\operatorname{dPreStk}$ denoted by $\Map(X, Y) \in \operatorname{dPreStk}$. $\Map(X, Y)$ is defined by sending $R$ to the mapping space $\Hom(X \times \Spec R, Y) \in \operatorname{Spc}$.
We denote the $\infty$-category of derived stacks for the \'etale topology by $\operatorname{dStk}$. A derived stack $X$ is an Artin stack if it is a geometric stack and locally of finite presentation. In particular, $X$ admits a perfect cotangent complex.
\subsection{Classical complex analytic geometry}
We will also need to work with complex analytic stacks since we will work with exponential maps. We define complex analytic stacks $\operatorname{Stk}^{\an}$ as in \cite[Definition 3.1.1]{sun_analytic}. In particular, $X$ is a complex analytic stack if it is a stack over the site of complex analytic spaces with analytic topology, there is a smooth surjective map $U \to X$ from a complex analytic space $U$ and the diagonal of $X$ satisfies a representability and finiteness condition. There is an analytification functor from finite type Artin stacks $(-)_{\an}\colon \operatorname{ArtStk}^{\leq 1} \to \operatorname{Stk}^{\an}$, see \cite[Section 3.2.2]{sun_analytic}. 
\begin{ex}
In this paper the main example of stacks we will use is the following. Let $G$ be an algebraic group acting on a finite type scheme $X$. Then $([X/G])_{\an} = [X_{\an} / G_{\an}] $. Here $[X_{\an} / G_{\an}]$ is the quotient of the groupoid $G_{\an} \times X_{\an} \rightrightarrows X_{\an}$.
\end{ex}
\begin{notation}
If it is clear from context we will abuse notation and denote a stack $X$ and its analytification $X_{\an}$ by the same symbol.
\end{notation}
\subsection{Derived complex analytic geometry}
We follow \cite{holstein2018analytification} for background on derived analytic geometry. In particular, there exists a derived analytification functor
\begin{equation}
    \an \colon \mathrm{dSt} \xrightarrow{j^{s}} \mathrm{dSt}^{\mathrm{afp}} \to \mathrm{dAnSt}
\end{equation}
where $\mathrm{afp}$ is for almost of finite presentation and $j^{s}$ is the restriction induced by the inclusion 
\begin{equation}
\mathrm{dAff}^{\mathrm{aft}} \xrightarrow{} \mathrm{dAff}.
\end{equation}
We will only use derived analytic geometry in subsection \ref{derived_complex exponential}.
\subsection{Differential forms and shifted symplectic structures}
In this section, we give a brief recollection of $n$-shifted closed forms and structures. For the original material, see \cite{PTVV} and \cite{Park2025} for a great introduction.
\subsection{de Rham algebras}
Recall from \cite{PTVV} the construction of the de Rham algebra $\mathrm
DR(A)$ of a cdga $A $, which is a graded mixed algebra.
We can now define a functor
\begin{align*}
    \DR(-)\colon \operatorname{cdga}^{\leq 0 } & \to \Mod^{\mathrm{gr}} \mathbb{C}[\epsilon] \\
    R & \mapsto \DR(R)
\end{align*}
It can be shown that this functor satisfies \'etale descent. For a general derived stack $X$, the graded mixed cdga $\DR(X)$ is then defined by right Kan extension along the map $\operatorname{cdga}^{\leq 0 } \to \operatorname{dStk}^{\op}$.  Concretely we have 
\begin{equation}
    \DR(X) = \lim_{\Spec R \to X} \DR(R)
\end{equation}
 Let $X$ be a derived stack that admits a cotangent complex, then there is a canonical map
\begin{equation}
    \Gamma(X, \Sym(\mathbb{L}_{X}[-1])) \to \DR(X).
\end{equation}
\begin{thm}\cite[Theorem 2.6]{cal_saf} \label{forms_artin}
Let $X$ be a derived  prestack that admits a cotangent complex. Then the above map is an equivalence of graded cdgas. 
\end{thm}
For Artin stacks locally of finite type, this theorem was already proven in \cite[Proposition 1.14]{PTVV}. Apart from Section \ref{s0_formal_sect}, we will work with Artin stacks locally of finite type. 
\subsection{Closed forms}
We can define as in \cite[Section 1.2]{PTVV} the spaces of $p$-forms
\begin{defn}
    We define
    \begin{enumerate}
        \item $n$-shifted $p$-forms 
        \begin{align*}
            \mathcal{A}^{p}(-,n) & \colon \cdgacon \to \Spc  \\
            R & \mapsto \Map(\mathbb{C}[-n],\DR(R))
        \end{align*}
        \item $n$-shifted closed $p$-forms
        \begin{align*}
            \mathcal{A}^{p,cl}(-,n) & \colon \cdgacon \to \Spc  \\
            R & \mapsto \Map(\mathbb{C}[p-n](p),\DR(R))
        \end{align*}
    \end{enumerate}
\end{defn}
These functors satisfy \'etale descent hence one can define the same spaces for a derived stack $X$.
We now define exact forms
\begin{defn}[Closed exact forms]
    We define 
    \begin{equation}
        \mathcal{A}^{p, ex}(-,n) = \fib(\mathcal{A}^{p,cl}(-,n) \to \mathcal{A}^{0,cl}(-,n+p))
    \end{equation}
\end{defn}
In particular, we can view the de Rham cohomology $\mathcal{A}^{p}(-,p+n)$ as the obstruction to the existence of a lift of a general $n$-shifted $p$-form to an exact form.
We are now particularly interested in $(-1)$-shifted closed forms
\begin{prop}\cite[Proposition 3.2]{kinjo2024cohomological} \cite[Proposition 2.2.11]{hennion2024gluing} \label{exact_neg1forms}
    We have a fiber sequence given by rotating the defining sequence of $\mathcal{A}^{2,ex}(-,-1)$
    \begin{equation}
        \mathcal{A}^{0,cl}(-,0) \to \mathcal{A}^{2,ex}(-,-1) \to \mathcal{A}^{2,cl}(-,-1)
    \end{equation}
    Furthermore this sequence splits giving that
    \begin{equation} \label{minus_one_split}
        \mathcal{A}^{2,ex}(-,-1) \cong \mathcal{A}^{0,cl}(-,0) \oplus \mathcal{A}^{2,cl}(-,-1)
    \end{equation}
\end{prop}
Concretely, this means that specifying an exact closed $2$-form on $X$ is the same as specifying a pair
\begin{align*}
    & (f, \alpha) \\
    & df=0 \\
    & d_{\mathrm{dR}}f +d\alpha = 0
\end{align*}
and the map giving the equivalence \eqref{minus_one_split} is given by
\begin{equation}
    \mathcal{A}^{2,ex}(X,-1) \to  \mathcal{A}^{0,cl}(X,0) \oplus \mathcal{A}^{2,cl}(X,-1)
\end{equation}
$(f, \alpha)$ is mapped to $d_{\mathrm{dR}} \alpha \in \mathcal{A}^{2,cl}(X,-1)$ under the projection.
\subsection{Shifted symplectic structures}
In classical algebraic geometry, a symplectic structure is a nondegenerate closed $2$-form $\omega$ on a scheme or manifold $X$. We can express the nondegeneracy condition by saying that the form gives an isomorphism $\To X \cong \To^{*}X$ between the tangent and cotangent bundles. This version of symplectic structure can be readily generalised as follows.
\begin{defn}[Shifted symplectic structure]\label{shifted_symp_def}
An $n$-shifted symplectic structure on a derived Artin stack $X$ is a closed $2$-form $\omega \in \mathcal{A}^{2,\cl}(X,n) $ along with the \textbf{non-degeneracy} condition that the induced map 
$$\mathbb{T}_{X} \to \mathbb{L}_{X}[n]$$
is a quasi-isomorphism. \par
Let $f\colon X \to Y$ be a map of derived stacks with $n$-shifted closed forms $\omega_{X}$ and $\omega_{Y}$. Then we say that the map preserves closed forms if $f^{*} \omega_{Y} \sim \omega_{X}$.
\end{defn}
\subsection{Volume forms}
We recall the definition of volume form in \cite{naef2023torsion}.
\begin{defn}[Volume form] \label{volume_forms_ns}
    Let $X$ be a derived stack with a perfect cotangent complex. Define the dimension of $X$ to be $\chi(\mathbb{L}_{X})$.
    \begin{enumerate}
        \item A volume form is an isomorphism 
        \begin{equation}
            \mathrm{vol}_{X} \colon \mathcal{O}_{X} \to \det \mathbb{L}_{X}.
        \end{equation}
        \item Let $f \colon X \to Y$ be a formally \'etale map of derived stacks equipped with volume forms $\mathrm{vol}_{X}$ and $\mathrm{vol}_{Y}$. Then there is an induced pullback volume form $f^{*} \mathrm{vol}_{Y}$ on $X$, which differs from $\mathrm{vol}_{X}$ by the invertible function $(df)^{\vee}$ induced by the isomorphism
        \begin{equation}
            \det f^{*} \mathbb{L_{X}} \xrightarrow{\det (df)^{\vee}}  \det \mathbb{L}_{Y}.
        \end{equation}
        \end{enumerate}
\end{defn}
\subsection{AKSZ construction}\label{aksz_sect}
Now let us describe how one can induce $n$-forms, closed $n$- forms and symplectic structures from $Y$ to $\Map(X, Y)$ via transgression. This is called in \cite{PTVV} the AKSZ construction. We recall several notions of orientation, which will be important to us in Chapter \ref{exp_map_chapter}.
\begin{defn}[Orientations and fundamental classes] \label{orientations_csh}
We summarise various definitions of orientations in \cite{calaque2022aksz}.
    \begin{enumerate}
    \item A $d$-preorientation on a map $f \colon X \to S$ is a morphism
    \begin{equation}
        [X] \colon f_{*} \mathcal{O}_{X} \to \mathcal{O}_{S}[-d]
    \end{equation}
    \item A stack $X$ over $S$ equipped with a $d$-preorientation is \textit{weakly} $d$-\textit{oriented} if for dualizable objects $\mathcal{E}$ the morphism
    \begin{equation}
        f_{*} \mathcal{E}^{\vee} \to (f_{*} \mathcal{E})^{\vee}[-d],
    \end{equation}
which is the adjoint of the morphism
\begin{equation}
    f_{*} \mathcal{E} \otimes f_{*} \mathcal{E}^{\vee} \xrightarrow{} f_{*}(\mathcal{E} \otimes \mathcal{E}^{\vee}) \xrightarrow{\mathrm{ev}} f_{*} \mathcal{O}_{X} \xrightarrow{[X]} \mathcal{O}_{S}[-d].
    \end{equation}
    \item \cite[Proposition 2.2.5]{calaque2022aksz} A stack $X$ over $S$ equipped with a $d$-preorientation is $d$-\textit{oriented} if for any morphism $\sigma \colon \mathrm{Spec} A \to S$, and a pullback diagram
    \begin{equation}
\begin{tikzcd}
	{X^{'}} & X \\
	{S^{'}} & S
	\arrow[from=1-1, to=1-2]
	\arrow[from=1-1, to=2-1]
	\arrow[from=1-2, to=2-2]
	\arrow["\sigma", from=2-1, to=2-2]
\end{tikzcd}
    \end{equation}
    we get a weak $d$-orientation on $X^{'}$ from the morphism
\begin{equation}
    f^{'}_{*} \mathcal{O}_{X^{'}} \cong \sigma^{*}f_{*} \mathcal{O}_{X} \xrightarrow{\sigma^{*} [X]} \sigma^{*} \mathcal{O}_{S}[-d] \cong  \mathcal{O}_{S^{'}}[-d].
\end{equation}
\item  Let $f \colon X \to S$ and $g \colon Y \to S$ be preoriented $S$-stacks. A cospan $X \to Z \xleftarrow{}Y$ is preoriented if the square
    \begin{equation}
\begin{tikzcd}
	{h_{*} \mathcal{O}_{Z}} & {f_{*} \mathcal{O}_{X}} \\
	{f_{*} \mathcal{O}_{Y}} & {\mathcal{O}_{S}[-d]}
	\arrow[from=1-1, to=1-2]
	\arrow[from=1-1, to=2-1]
	\arrow[from=1-2, to=2-2]
	\arrow[from=2-1, to=2-2]
\end{tikzcd}
    \end{equation}
commutes with $h \colon Z \to S$.
    \item let $X$ satisfy assumption 1.1 in \cite{naef2023torsion}. Assume that there is an adjunction $p^{*}[-d] \dashv p_{\#}$. 
    Then a \textbf{fundamental class} of degree $d$ is the unit morphism morphism $\mathbb{C} \to p_{\#} \mathcal{O}_{X}[-d]$.
\end{enumerate}
\end{defn}
\begin{remark}
    Note that given a fundamental class we can write the following chain of adjunctions
    \begin{equation}
        p^{*}[-d] \dashv p_{\#} \dashv p^{*} \dashv p_{*}
    \end{equation}
    shifting by $d$ we then get
    \begin{equation}
        p^{*} \dashv p_{\#}[d] \dashv p^{*}[d] \dashv p_{*}[d]
    \end{equation}
    but by uniqueness of adjoints this implies that there is a natural isomorphism
    \begin{equation}
        p_{*} \xrightarrow{\cong} p_{\#}[d]
    \end{equation}
We view this isomorphism as giving us Poincar\'e duality on $X$.
\end{remark}
The following proposition proves that a fundamental class gives an $\mathcal{O}$-orientation.
\begin{prop}\cite[Proposition 1.22]{naef2023torsion}
    Let $X$ be a prestack satisfying  assumption 1.1 in \cite{naef2023torsion}. Then we have the following equivalences
    \begin{enumerate}
        \item the data of an $\mathcal{O}$-orientation $p_{*}\mathcal{O}_{X} \to \mathbb{C}[-d]$
        \item a Poincar\'e duality isomorphism $p_{*} \to p_{\#}[-d]$
        \item a fundamental class $\mathbb{C} \to p_{\#} \mathcal{O}_{X}[-d]$.
    \end{enumerate}
\end{prop}
\begin{remark}
    Note that in \cite{naef2023torsion} the stronger assumption on $X$ is asked. However, to obtain the part of the Proposition 1.22 we are interested in, it is enough to restrict to existence of the left adjoint of the pullback $p^{*}$.
\end{remark}
\begin{ex} \label{fund_classes example}
We give important examples of fundamental classes.
\begin{enumerate}
    \item Let $M$ be a finitely dominated space. Then a fundamental class on $M$ induces a fundamental class on $M_{\mathrm{B}}$. This is equivalent to having usual Poincar\'e duality.
    \item Let $X$ be a smooth and proper scheme. Then a fundamental class in de Rham and respectively Dolbeault cohomology of  $X$ induces a fundamental class on $X_{\mathrm{dR}}$ and $X_{\mathrm{Dol}}$. 
    \item Let $\mathfrak{g}$ a complex Lie algebra of dimension $d$. Then we have Poincar\'e duality \cite[Chapter 5 Section 3]{knapp1988lie}
    \begin{equation}
        \mathrm{H}^{n}(\mathfrak{g},V) \cong \mathrm{H}_{d-n}(\mathfrak{g},V \otimes_{\mathbb{C}} (\wedge^{d} \mathfrak{g})^{\vee})
    \end{equation}
    If $\mathfrak{g}$ acts trivially on $\wedge^{d} \mathfrak{g}$, then we have a Poincar\'e duality on Lie algebra cohomology and thus a fundamental class. See also \cite[Theorem 4.9]{holstein2024koszul}
\end{enumerate}
\end{ex}
We will explicitly use fundamental classes in Section \ref{exp_shifted_tangents}
\subsection{Transgression of forms}
We are now ready to state the main theorem which describes a procedure to induce differential forms from the target $Y$ of a mapping stack to the full space $\Map(X,Y)$. The resulting form on $\Map(X,Y)$ is often called the transgressed form.
\begin{thm}[AKSZ] \cite[Theorem 2.5]{PTVV} \cite[Section 3.2]{calaque2022aksz}\label{aksz}
Let $X$ be an $\mathcal{O}$-compact stack with an orientation of degree $d$ and $Y$ a derived artin stack. 
There is a map
\begin{equation}
    (-)_{\aksz} \colon \mathcal{A}^{p,cl}(Y,n) \xrightarrow{ \int_{[X]}\ev^{*}}  \mathcal{A}^{p,cl}(\Map(X,Y),n-d).
\end{equation}
Here the map
\begin{equation}
    \int_{[X]} \colon \mathcal{A}^{p,cl}(X \times \Map(X,Y),n) \to  \mathcal{A}^{p,cl}(\Map(X,Y),n-d)
\end{equation}
which is induced by the pushforward $p_{\#}$ along the map $p \colon  X \times \Map(X,Y) \to \Map(X,Y)$. We have the following properties of this construction
\begin{enumerate}
    \item If $Y$ is $n$-shifted symplectic then the above map preserves non-degeneracy and there is a canonical $(n-d)$-shifted symplectic structure on $\Map(X,Y)$.
    \item The construction is covariantly functorial in the target. Let $f \colon Y_{1} \to Y_{2}$ be a map of derived artin stacks. Denote the induced map by $f_{X} \colon \Map(X,Y_{1}) \to \Map(X,Y_{2})$. Then there is a commutative diagram
\begin{equation}
\begin{tikzcd}
	{\mathcal{A}^{p,cl}(Y_{2},n)} & {\mathcal{A}^{p,cl}(\Map(X,Y_{2}),n-d)} \\
	{\mathcal{A}^{p,cl}(Y_{1},n)} & {\mathcal{A}^{p,cl}(\Map(X,Y_{2}),n-d)}
	\arrow["{(-)_{\aksz}}", from=1-1, to=1-2]
	\arrow["{f^{*}}"', from=1-1, to=2-1]
	\arrow["{f^{*}_{X}}", from=1-2, to=2-2]
	\arrow["{(-)_{\aksz}}"', from=2-1, to=2-2]
\end{tikzcd}
\end{equation} 
\end{enumerate} 
Furthermore, if $f$ is $n$-symplectic, then the induced map $f_{X}$ is also $n-d$-symplectic.
\end{thm}
\subsection{Formal stacks and shifted symplectic structures on formal stacks}
In this thesis we will need to work with non algebraic maps and behaviour of shifted symplectic structures under such maps we will generally follow \cite{cptvv} and \cite{gaitsgory2017study}. In particular, we will need parts of this theory when working with the exponential map in Chapter \ref{exp_map_chapter} and d-critical structures on formal completions in Section \ref{s0_formal_sect}. We will work with formal derived stacks as in \cite[Definition 2.1.1]{cptvv}. In particular, any derived algebraic stack is a formal stack. We now introduce the de Rham stack and formal completions, which the examples of formal stacks we will consider. \par
Consider the inclusion map $i: \mathrm{alg}^{\mathrm{red}} \to \mathrm{cdga}^{\leq0}$ including reduced classical algebras into cdgas. This functor has a left adjoint 
\begin{equation}
\begin{tikzcd}
	{\mathrm{cdga}^{\leq0}} && {\mathrm{alg}^{\mathrm{red}}}
	\arrow["{(-)_{\mathrm{red}}}", shift left=2, from=1-1, to=1-3]
	\arrow["i", shift left=2, from=1-3, to=1-1]
\end{tikzcd}
\end{equation}
Here $(A)_{\mathrm{red}} = (\mathrm{H}^{0}(A))_{\mathrm{red}}$  \\
This extends to an adjunction, where $i^*$ is left adjoint to $i_{*}$ and $i_{!}$ is left adjoint to $i^{*}$
\begin{equation}
\begin{tikzcd}
	{\mathrm{dStk}} && {\mathrm{Stk}_{\mathrm{red}}}
	\arrow["{i^{*}}", shift right=2, from=1-1, to=1-3]
	\arrow["{i_{*}}", shift left=5, from=1-3, to=1-1]
	\arrow["{i_{!}}"', shift right=4, from=1-3, to=1-1]
\end{tikzcd}
\end{equation}
\begin{defn}[de Rham stack]
Define the de Rham stack functor
\begin{equation}
    (-)_{\mathrm{dR}}: \mathrm{dStk} \to \mathrm{dStk} = i^{*}i_{*}.
\end{equation}
On points we can compute that $(X)_{\mathrm{dR}}(A) = X(\mathrm{H}^{0}(A)_{\mathrm{red}})$. By adjunction, there is a canonical map $X \to X_{\mathrm{dR}}$.
\end{defn}
\begin{defn}[Formal completion as relative de Rham stack]
Let $X \to Y$ be a map. Define the formal completion by the pullback
\begin{center}
\begin{tikzcd}
	{\hat{Y}_{X}} & {X_{\mathrm{dR}}} \\
	Y & {Y_{\mathrm{dR}}}
	\arrow[from=1-1, to=1-2]
	\arrow[from=1-1, to=2-1]
	\arrow[from=1-2, to=2-2]
	\arrow[from=2-1, to=2-2]
\end{tikzcd}
\end{center}
On points we have $\hat{Y}_{X}(A) = Y(A) \times_{Y(A_{\mathrm{red}})} X(A_{\mathrm{red}})$.
\end{defn}
We now check that the more familiar notion of completion of a scheme along a closed subscheme agrees with the above notion of formal completion. Let $X$ be a finite type scheme and $\widehat{X}^{x}$ its completion at a point.

\begin{lem}\cite[Lemma B.1.2]{cptvv} \label{cptvv_formal_lemma}
    Let $R$ be a noetherian classical ring and denote $X = \Spec R$. Let $X_{n} = \Spec (R/I^{n})$, with $I$ being the ideal defining $x$ in $\Spec R$, then we have an equivalence of prestacks
    \begin{equation}
        \widehat{X}^{x} \coloneqq X \times_{X_{\dR}} \pt \cong \colim_{n \in \mathbb{N}} X_{n}
    \end{equation}
\end{lem}
\begin{remark}
 Note that if we consider a classical stack and embed it into derived prestacks it is generally only locally almost of finite type. Similarly formal completions of finite type schemes are locally almost of finite type by \cite[Corollary 6.3.2]{dg_indschemes}. This becomes an important point for us when we consider d-critical structures on formal completions in Section \ref{s0_formal_sect}.  
\end{remark}
\begin{prop} \cite[Proposition 2.1.4]{cptvv}
Let $F$ be a derived stack and $f \colon X \to Y$ be a map of derived stacks.
\begin{enumerate}
\item   $F_{\mathrm{dR}}$ is a formal derived stack and $\widehat{Y}^{X}$ is a formal derived stack.
    \item The map $F \to F_{\mathrm{dR}}$ induces an equivalence $F_{\mathrm{red}} \to (F_{\mathrm{dR}})_{\mathrm{red}}$. The canonical map $X \to \widehat{X}^{f} $ induces an equivalence $X_{\mathrm{red}} \to \widehat{X}^{f}_{\mathrm{red}} $.
    \item Let $j \colon t_{0} F \to F$ be the canonical map from the classical truncation of $F$. Then we have an equivalence $\widehat{F}^{j} \to F$.
\end{enumerate}
\end{prop}

\section{D-critical loci and DT sheaves} \label{mhm_background}
\subsection{Constructible sheaves and  perverse sheaves}
In this section we will collect some aspects of the theory of constructible and perverse sheaves that we will use in the thesis. We assume the reader is familiar with the theory so we will not give a full introduction. We refer to \cite{de2009decomposition} and \cite{achar2021perverse} for great introductions.   \par We work with vector spaces with coefficients in $\mathbb{Q}$, and from now on we will drop the coefficients from the notation.
Let $X$ be a finite type $\mathbb{C}$-scheme or complex analytic space, when considering (constructible) sheaves on these spaces we are implicitly considering them over $X_{\red}$. Every functor written is also implicitly derived.
\begin{remark}
    When talking about local systems as part of the category of perverse sheaves we will sometimes drop the shift by dimension and write $\mathcal{L} \in \mathrm{Perv}(X)$ instead of $\mathcal{L}[\dim X]$.
\end{remark}
\subsection{D-critical loci and DT sheaves}
In this section we will define d-critical structures, which we view as classical truncations of  $(-1)$-shifted symplectic stacks and schemes.   We will briefly recall the construction of the global DT sheaf on a d-critical scheme or stack. 
\subsection{D-critical loci on schemes.}
\begin{prop}\cite[Theorem 6.1]{brav2019darboux} \label{s_sheaf}
Let $X$ be a $\mathbb{C}$-scheme or a complex analytic space. For every open $R \xrightarrow{j} X$ with a closed embedding $R \xrightarrow{i} U$ into a smooth scheme $U$ we have a short exact sequence
    $$0 \to I_{R,U} \xrightarrow{i} i^{-1} \mathcal{O}_{U} \to \mathcal{O}_{R} \to 0.$$
There exists a sheaf $S_{X}$ of $\mathbb{C}$-vector spaces on $X$ such that the following hold 
\begin{enumerate}
    \item there is an exact sequence $0 \to S_{X} |_{R} \xrightarrow{i} i^{-1} \mathcal{O}_{U} / I_{R,U}^{2} \xrightarrow{d} i^{-1} \Omega_{U}/ I_{R,U} \cdot i^{-1} \Omega_{U}$
    \item there is a morphism $\beta_{X} \colon S_{X} \to \mathcal{O}_{X}$ inducing a decomposition $S_{X} = \mathbb{C}_{X} \oplus S^{0}_{X}$, with 
$S^{0}_{X} = \ker(S_{X} \xrightarrow{\beta_{X}}  \mathcal{O}_{X} \to \mathcal{O}_{X,\red})$ and $\mathbb{C}_X$ the constant sheaf.
 \item we have an exact sequence
$$0 \to \HH^{-1}(\mathbb{L}_{X}) \to S_{X} \to \mathcal{O}_{X} \to \Omega_{X} $$
and a similar one for $S^{0}_{X}$.
\end{enumerate}
\end{prop}
\begin{defn}[D-critical locus] \cite[Definition 3.1]{ben2015darboux}
A $d$-critical structure on a scheme $X$ or complex analytic space is the data of
\begin{enumerate}
    \item section $s \in \HH^{0}(S^{0}_{X})$
    \item for each point $p \in X$ a \textbf{critical chart} $ (R, U , f , i) $ with $p \in R$ an open of $X$, $i \colon R \to U$ a closed embedding,  $U$ smooth and $U \xrightarrow{f} \mathbb{C}$ a function satisfying $s|_{R} = i^{-1}f + I^{2}_{R,U}$ and $i(R) = \crit f$.
\end{enumerate}
Let $f \colon X \to Y$ be a morphism of schemes or complex analytic spaces with d-critical loci structures, then there is an induced map $f^{\star} \colon f^{-1} S^{0}_{Y} \to S^{0}_{X}$. We say $f$ is a morphism of d-critical loci if $f^\star s_Y = s_X$.  
\end{defn}
\subsection{Canonical bundle and orientations.}
\begin{defn}[Canonical bundle and orientation]   \cite[Theorem 6.4]{brav2019darboux}
Let $(X,s)$ be a $d$-critical locus. Then there exists a line bundle $K_{X}$ on $X_{\red}$ with the property that for every critical chart $(R,U,f,i)$ there is an isomorphism
$$i \colon K_{X}|_{R_{\red}} \to (\omega_{U}^{\otimes 2}) |_{R_{\red}}.$$
An \textbf{orientation} for a $d$-critical locus is the data $(L,\theta)$ with $L$ a line bundle on $X$ and  an isomorphism $\theta \colon L^{\otimes 2} \to K_{X}$. An isomorphism of two orientations $\psi \colon (L_{1}, \theta_{1}) \to (L_{2}, \theta_{2})$   is given by an isomorphism $\psi \colon L_{1} \to L_{2}$ which satisfies $\theta_{2} \circ \psi^{\otimes 2}  = \theta_{1}$.
\end{defn}
Now we  state the theorem about the passage from $(-1)$-shifted symplectic derived schemes to $d$-critical loci.
\begin{thm} \cite[Theorem 6.6]{brav2019darboux}\label{truncation}
Let $X$ be a $(-1)$-shifted symplectic derived scheme.
Then the truncation $\tt_{0}(X)$ has a natural structure of $d$-critical locus with $(\det \mathbb{L}_{X} )|_{\tt_{0}X_{\red}} = K_{\tt_{0}(X)}$.
\end{thm}
Note that if we can find a volume form for a derived stack $X$, namely an equivalence $\mathcal{O}_{X} \to \det \mathbb{L}_{X}$, then on the truncation of $X$ we get an orientation of the induced d-critical structure.
\subsection{DT sheaf on stacks}
A similar story holds for $(-1)$-shifted symplectic derived Artin stacks. There is a truncation to $d$-critical stacks and similarly a perverse sheaf. We will briefly recall the constructions. First, we need to define sheaves on Artin stacks. As in \cite[Section 2.7]{Joyce_dcrit} we work with the site Lis-\'et$(X)$. A sheaf $\mathcal{F}$ on $X$ will be the data of an \'etale sheaf $\mathcal(F)_{T}$ for every smooth map $f \colon T \to X$ from a scheme $T$ with some compatibility conditions.
\begin{prop}[D-critical structures for stacks]  \cite[Corollary 2.52]{Joyce_dcrit} \label{d_critc_stacks}
Let $X$ be an Artin stack or complex analytic stack.
\begin{enumerate}
    \item We have the following 
\begin{enumerate}
    \item there exists a sheaf $S_{X}$ of $\mathbb{C}$-vector spaces on $X$ such that for each smooth morphism $f \colon T \to X$ we have an isomorphism $\theta_{f} \colon f^{*} S_{X} \to S_{T}$ 
    \item there is a canonical splitting $S_{X} = \mathbb{C}_{X} \oplus S^{0}_{X} $.
\end{enumerate}
\item $X$ has a d-critical structure if there is a section $s_{X} \in S^{0}_{X}$ such that for each smooth morphism $f \colon T \to X$ we have that $f^{*} s_{X}$ defines a d-critical structure on $T$. We call $X$ a d-critical stack.
\item Let $X$ be a d-critical stack. Then there is a canonical line bundle $K_{X}$ on $X_{\red}$. An \textbf{orientation} on a d-critical stack $X$ is the data $(L, \theta)$ of a line bundle $L$ on $X_{\red}$ and an isomorphism $\theta \colon L^{\otimes 2} \to K_{X}$.
\end{enumerate}
\end{prop}
Again we have a truncation theorem
\begin{thm}\cite[Theorem 3.18]{ben2015darboux} \label{truncation_dcrit_stacks}
Let $X$ be a $-1$-shifted symplectic derived artin stack. Then the truncation $\tt_{0}(X)$ has a natural structure of $d$-critical stack and $(\det \mathbb{L}_{X} )|_{(\tt_{0}X)_{\red}} \cong K_{\tt_{0}(X)}$. 
\end{thm}
Then similarly one can define the DT sheaf on stacks.
\begin{thm} \cite[Theorem 4.8]{ben2015darboux}  \label{joyce_sheaf_stacks}
Let $(X,s)$ be a $d$-critical stack with an orientation $(K^{1/2}_{X}, \theta)$. Then there is a perverse sheaf $\varphi_{X}$ on $(X,s)$ such that for each smooth map $f \colon T \to X$ we have $f^{*}[d] \varphi_{X} \cong \varphi_{T} $. Here $d$ is the relative dimension of $f$ and $T$ has the induced $d$-critical structure from $X$. Furthermore, there is an upgrade of $\varphi$ to a mixed Hodge module on $X$.
\end{thm}
\begin{ex}[Products] \cite[Proposition 4.3]{kinjo2024cohomological} \label{dcrit_products}
Let $X$ and $Y$ be $(-1)$-shifted symplectic oriented stacks. Then  the d-critical locus structure on $X \times Y$ is given by $s_{X} \oplus s_{Y}$ and we have that $\varphi_{X \times Y} \cong \varphi_{X} \boxtimes \varphi_{Y} $. Also see \cite[Remark 5.23]{joyce_conje_lino_ben}.
\end{ex}
We will need to be a bit more explicit about d-critical structures on quotient stacks. We recall \cite[Section 3.2]{ben2015darboux}, where it is explained that $d$-critical structures on a quotient stack are the same as $G$-equivariant $d$-critical structures. Let $G$ be an algebraic group acting on a scheme $X$ and denote the action map by $a \colon G \times X \to X$ and the projection map by $\pi \colon G \times X \to X$. Then a $G$-equivariant $d$-critical structure is a section $s \in \Gamma(X, S^{0}_{X})$, with the property that $\pi^{\star} s = a^{\star} s \in \Gamma(G \times X, S^{0}_{G \times X})$. Equivalently for each $g \colon X \to X$ we have that $g^{\star}s = s$, where $g$ is the map induced by the action of $g \in G$.
\section{\texorpdfstring{$S^{0}$}{S 0} sheaves for formal completions} \label{s0_formal_sect}
In this section we establish some folklore results about $S$ sheaves on formal completions as well as comparisons between the $S$ sheaves of an Artin stack $X$, its analytification $X_{\an}$ and its formal completion at a point $\widehat{X}^{x}$. The results in this section are technical in nature and will only be used in the proof of Theorem \ref{main_dcrit_thm}.
\subsection{Closed forms on formal completions }
Note that if we consider a classical stack and embed it into derived prestacks it is generally only locally almost of finite type. Similarly formal completions of finite type schemes are locally almost of finite type by \cite[Corollary 6.3.2]{dg_indschemes}.  First, let $X$ be a finite type scheme and $\widehat{X}^{x}$ its completion at a point. Then we can use Lemma \ref{cptvv_formal_lemma} to identify the de Rham stack definition of completion with the more usual definition of completion for closed immersions of schemes.
This allows us to identify $\QCoh (\widehat{X}^{x}) \cong \lim_{n \in \mathbb{N}} \QCoh (X_{n})$ and the functor $\eta^{*}$ induced by $\eta \colon \widehat{X}^{x} \to X$ is $\eta^{*}(M) = (M / I^{n} M)_{n \in \mathbb{N}}$. We also have the completion $\widehat{M} = \lim M/I^{n}M$. We can identify $\Gamma(\widehat{X}^{x}, \eta^{*} M) = \widehat{M}$.
When considering formal completions of algebraic stacks we can reduce to the case of quotient stacks in the following way. By \cite[Theorem 4.12]{AHR_luna} there is an \'etale map $  Y / G_{x} \to X$ for $Y = \Spec R$.  Therefore, the formal completions of $X$ and $Y/ G_{x}$ will coincide and we have $\widehat{X}^{\B G_{x}} = \widehat{Y}^{x} / G_{x}$ and $\widehat{X}^{x} = \widehat{Y}^{x} / \widehat{G}^{1}_{x}$. Here $\widehat{X}^{\B G_{x}}$ is the completion along the map $\B G_{x} \to X$.
Using Proposition \ref{exact_neg1forms} for a cdga $R$ we can rewrite 
\begin{equation}
    \mathcal{A}^{2,ex}(\Spec R,-1) \cong \cofib(\DR(R)(0) \xrightarrow{\epsilon} \DR(R)(1)[1])[-1]
\end{equation}
the map $\DR(R)(0) \xrightarrow{\epsilon} \DR(R)(1)[1]$ can be further rewritten as the map
\begin{equation}
    \cofib(R \xrightarrow{d_{\dR}} \mathbb{L}_{R})[-1].
\end{equation}
By right Kan extension we can then also write 
\begin{equation}
    \mathcal{A}^{2,\exa}(X) = |\cofib(\DR(X)(0) \xrightarrow{\epsilon} \DR(X)(1)[1])[-1]|.
\end{equation}
In particular, we are interested in $X = \widehat{X}^{x}$ a completion of a finite type scheme at a point. In this case, we can use Theorem \ref{forms_artin} and Proposition \ref{exact_neg1forms} to deduce that
\begin{align*}
    \mathcal{A}^{2,\exa}(\widehat{X}^{x}) & \cong |\cofib(\Gamma(\widehat{X}^{x},\mathcal{O}_{\widehat{X}^{x}} \xrightarrow{d_{\dR}} \mathbb{L}_{\widehat{X}^{x}}))[-1]| \\
    \mathcal{A}^{2,\exa}(X) & \cong |\cofib(\Gamma(X,\mathcal{O}_{X} \xrightarrow{d_{\dR}} \mathbb{L}_{X}))[-1]|.
\end{align*}
The generality of Theorem \ref{forms_artin} is necessary here since we cannot use \cite[Proposition 1.14]{PTVV} since $\widehat{X}^{x}$ is not Artin and $X$ is not of finite type as a derived prestack. 
We will now give a definition of sections of the $S$-sheaf which works for formal completions.
\begin{defn}[S sheaf for formal completions] \label{new_S_sheaf_def}
Let $X$ be a classical Artin stack. If $\widehat{X}^{x}$ is the completion of an Artin stack at a point $x$, we define the vector spaces $S_{\widehat{X}^{x}} = \pi_{0}\mathcal{A}^{2,\exa}(\widehat{X}^{x},-1)$ and $S^{0}_{\widehat{X}^{x}} = \pi_{0}\mathcal{A}^{2,\cl}(\widehat{X}^{x},-1)$. 
\end{defn} 
We can now consider a restriction $\Gamma(X,S_{X}) \to S_{\widehat{X}^{x}}$, which is given by 
\begin{equation}
    (f, \alpha) \mapsto (\widehat{f}, \widehat{\alpha})
\end{equation}
where $\widehat{f}$ and $\widehat{\alpha}$ mean the restrictions  along the map $\widehat{X}^{x} \to X$.
The following proposition now ensures that the definition above is compatible with the original Definition \ref{s_sheaf}.
\begin{prop} \label{closed_forms_ssheaf}
 Let $X$ be a classical scheme or Artin stack. Then we have 
 \begin{equation}
     \pi_{0}\mathcal{A}^{2, \exa }(X,-1) \cong \Gamma(X, S_{X}),
 \end{equation}
 where the sheaf $S_{X}$ is as defined in \ref{s_sheaf}.   \par 
 Let $\eta_{X} \colon \widehat{X}^{x} \to X$ be the formal completion of a classical scheme at a point $x \colon \pt \to X$.
    Fix an open neighbourhood $R \subseteq X$ of $x$ and a closed immersion $R \xrightarrow{i} U$ with ideal $I$. This induces a map $\widehat{\mathcal{O}}^{x}_{U} \to \widehat{\mathcal{O}}^{x}_{R} = \widehat{\mathcal{O}}^{x}_{X}$. Then  $S_{\widehat{X}^{x}}$ fits into the following short exact sequence of vector spaces
    \begin{equation} \label{formal_s_sheaf_ses}
        0 \to S_{\widehat{X}^{x}} \to \widehat{\mathcal{O}}^{x}_{U} / \widehat{I}^{2} \to  \widehat{\Omega}^{x}_{U} / \widehat{I} \widehat{\Omega}^{x}_{U}.
    \end{equation}
\end{prop}
\begin{proof}
The equivalence of the two definitions of sections of the $S$ sheaf follow by \cite[Remark 2.2b]{Joyce_dcrit}. In particular, one can consider the truncation of the cotangent complex of $X$. Given an embedding $i \colon X \to U$ for $U$ smooth we have
\begin{equation} \label{cot_truncated}
    \tau_{\geq -1} \mathbb{L}_{X} = I/I^{2} \to i^{*} \Omega_{U}
\end{equation} 
We can form the following exact sequence of two term complexes
\begin{equation} \label{ses_S_sheaf_cone}
\begin{tikzcd}
	{I/I^{2}} & {\mathcal{O}_{U}/I^{2}} & {\mathcal{O}_{X}} \\
	{i^{*}\Omega_{U}} & {i^{*}\Omega_{U}} & 0
	\arrow[from=1-1, to=1-2]
	\arrow[from=1-1, to=2-1]
	\arrow[from=1-2, to=1-3]
	\arrow[from=1-2, to=2-2]
	\arrow[from=1-3, to=2-3]
	\arrow[from=2-1, to=2-2]
	\arrow[from=2-2, to=2-3]
\end{tikzcd}
\end{equation}
Shifting we can show that $\cone(\mathcal{O}_{X} \to \tau_{\geq -1}\mathbb{L}_{X}) \cong \mathcal{O}_{U}/I^{2} \to i^{*}\Omega_{U}$. From this we can see that \begin{equation*}
     \pi_{0}\mathcal{A}^{2,\exa}(X,-1) \cong \Gamma(X, S_{X}).
 \end{equation*}
We will now repeat the same proof for formal completions. Note for the purposes of formal completion we can work affine locally so we can assume $X = \Spec A$ is affine and we have a closed embedding $X \to U$ with $U= \Spec R$ smooth and affine. Then we have an ideal $I \subseteq R$ such that $A = R/I$. \par
    The second property will follow from a description of the cotangent complex of the formal completion of $X$.  We get induced maps on formal completions that make the following square commute.
    \begin{equation}
\begin{tikzcd}
	X & U \\
	{\widehat{X}^{x}} & {\widehat{U}^{x}}
	\arrow["i"', from=1-1, to=1-2]
	\arrow["{\eta_{X}}", from=2-1, to=1-1]
	\arrow["{\hat{i}}", from=2-1, to=2-2]
	\arrow["{\eta_{U}}", from=2-2, to=1-2]
\end{tikzcd}
    \end{equation}
 Since we are considering $\mathcal{A}^{2,\exa}(\widehat{X}^{x}, -1)$, it is enough to consider the truncation $\tau_{\geq -1} \mathbb{L}_{\widehat{X}^{x}}$. 
The maps $\eta_{X}$ and $\eta_{U}$ are formally \'etale so we get $\tau_{\geq -1}\mathbb{L}_{\widehat{X}^{x}} =  \eta^{*}_{X} \tau_{\geq -1} \mathbb{L}_{X}$. The inverse systems $\eta^{*} I/I^{2}$ and $\eta^{*} i^{*}\Omega_{U}$ are Mittag-Leffler since all the maps in the inverse system are surjective, therefore the limit functor does not have any higher cohomology. Then using Lemma \ref{cptvv_formal_lemma} and equation \eqref{cot_truncated} we can deduce that on global sections on $\widehat{X}^{x}$ we have
\begin{equation} \label{cot_truncated_formal}
    \tau_{\geq -1} \mathbb{L}_{\widehat{X}^{x}} = \widehat{I}/\widehat{I}^{2} \to i^{*} \widehat{\Omega}_{U}.
\end{equation} 
Now we can consider the following short exact sequence of complexes, which comes from completion of the analogous exact sequence \ref{ses_S_sheaf_cone}. Using the Mittag-Leffler condition again we get
\begin{equation}
\begin{tikzcd}
	{\widehat{I}/\widehat{I}^{2}} & {\mathcal{O}_{\widehat{U}}/\widehat{I}^{2}} & {\mathcal{O}_{\widehat{X}}} \\
	{\widehat{i}\Omega_{\widehat{U}}} & {\widehat{i}\Omega_{\widehat{U}}} & 0
	\arrow[from=1-1, to=1-2]
	\arrow[from=1-1, to=2-1]
	\arrow[from=1-2, to=1-3]
	\arrow[from=1-2, to=2-2]
	\arrow[from=1-3, to=2-3]
	\arrow[from=2-1, to=2-2]
	\arrow[from=2-2, to=2-3]
\end{tikzcd}
\end{equation}
By shifting this short exact sequence it follows that we have an isomorphism $\cone(\mathcal{O}_{\widehat{X}} \to \tau_{\geq -1}\mathbb{L}_{\widehat{X}}) \cong (\widehat{\mathcal{O}}^{x}_{U}/ \widehat{I}^{2} \to \widehat{\Omega}^{x}_{U})$. Since we defined the space of exact $-1$ forms to be the cohomology of the cone we get the desired description of $S^{0}_{\widehat{X}^{x}}$. 
\end{proof}
From this point of view we can express the induced $d$-critical structure on $\tt_{0} X$ in Theorem \ref{truncation_dcrit_stacks} as the one induced by the map $\mathcal{A}^{2,\cl}(X,-1) \to \mathcal{A}^{2,\cl}(\tt_{0} X , -1)$. Proposition $\ref{exact_neg1forms}$ now shows that we have a decomposition $S_{\widehat{X}^{x}} = S^{0}_{\widehat{X}^{x}} \oplus \mathbb{C}$.
\begin{remark}
    Note that in the classical or analytic setting it makes sense to also define the $S$ sheaf on $X$ for the \'etale or Zariski topologies on $X$. For formal completions at a point the underlying space is just a point so we only have a vector space. 
\end{remark}
\begin{lem} \label{ss_ssheaf_inj}
    Let $X$ be a finite type scheme. We have
    \begin{enumerate}
        \item a map $\Gamma(X,S^{0}_{X}) \to \Gamma(X_{\an},S^{0}_{X_{\an}})$,
        \item an injective map on stalks $S^{0}_{X,x} \to S^{0}_{X_{\an},x}$,
        \item injective maps $S^{0}_{X,x} \to  S^{0}_{\widehat{X}^{x}}$ and $S^{0}_{X_{\an},x} \to  S^{0}_{\widehat{X}^{x}}$,
        \item a commutative diagram \begin{equation} \label{stalk_alg_formal_an}
\begin{tikzcd}
	& {S^{0}_{X,x}} \\
	{S^{0}_{\widehat{X}^{x}}} & {S^{0}_{X_{\an},x}}
	\arrow[hook, from=1-2, to=2-1]
	\arrow[hook, from=1-2, to=2-2]
	\arrow[hook, from=2-2, to=2-1]
\end{tikzcd}
    \end{equation}
    \end{enumerate}
    \end{lem}
\begin{proof}
Denote by $h \colon X_{\an} \to X$ the inclusion map. To define the map $\Gamma(X,S^{0}_{X}) \to \Gamma(X_{\an},S^{0}_{X_{\an}})$, as in \cite[Section 3.1]{Joyce_dcrit} we can cover $X$ by opens $R$  such that $R \xhookrightarrow{} U$ is a closed embedding into $U$ smooth. We can then use the following diagram
    \begin{equation}
\begin{tikzcd}
	R & U & {\mathcal{O}_{U}/I^{2}_{U}} & {\Omega_{U}/I_{U}\Omega_{U}} \\
	{R_{\an}} & {U_{\an}} & {h_*(\mathcal{O}_{U_{\an}}/I^{2}_{U_{\an}})} & {h_{*}(\Omega_{U_{\an}}/I_{U_{\an}}\Omega_{U_{\an}})}
	\arrow["i", hook, from=1-1, to=1-2]
	\arrow[from=1-3, to=1-4]
	\arrow[from=1-3, to=2-3]
	\arrow[from=1-4, to=2-4]
	\arrow[from=2-1, to=1-1]
	\arrow["{i_{\an}}", hook, from=2-1, to=2-2]
	\arrow[from=2-2, to=1-2]
	\arrow[from=2-3, to=2-4]
\end{tikzcd}
    \end{equation}
    which induces the map of short exact sequences
    \begin{equation}
\begin{tikzcd}
	0 & {S_{X}|_{R}} & {i^{-1}\mathcal{O}_{U}/I^{2}} & {i^{-1} \Omega_{U}/Ii^{-1}\Omega_{U}} \\
	0 & {h_{*}(S_{X_{\an}}|_{R_{\an}})} & {h_{*}i^{-1}_{\an}\mathcal{O}_{U_{\an}}/I^{2}_{\an}} & {h_{*}i^{-1}_{\an} \Omega_{U_{\an}}/I_{\an}i^{-1}_{\an}\Omega_{U_{\an}}}
	\arrow[from=1-1, to=1-2]
	\arrow[from=1-2, to=1-3]
	\arrow[dashed, from=1-2, to=2-2]
	\arrow[from=1-3, to=1-4]
	\arrow[from=1-3, to=2-3]
	\arrow[from=1-4, to=2-4]
	\arrow[from=2-1, to=2-2]
	\arrow[from=2-2, to=2-3]
	\arrow[from=2-3, to=2-4]
\end{tikzcd}
    \end{equation}
    This defines a map $S^{0}_{X} \to h_{*} S^{0}_{X_{\an}}$ and thus a map $\Gamma(X, S^{0}_{X}) \to \Gamma(X_{\an}, S^{0}_{X_{\an}})$. This proves $1$. \par
    If the map $i^{-1} \mathcal{O}_{U} / I^{2} \to i^{-1} \mathcal{O}_{U_{\an}} / I^{2}_{\an}$ is injective, then the map $S_{X}|_{R} \to S_{X_{\an}|_{R_{\an}}}$ is injective. Note that $I = i^{-1}I_{U}$. The  sheaf $O_{U}/I^{2}_{U}$ is coherent so the canonical map to the analytification is injective. This also means that the map on stalks $S^{0}_{X,x} \to S^{0}_{X_{\an},x}$ is injective. This proves $2$. \par 
    The maps $S^{0}_{X,x} \to  S^{0}_{\widehat{X}^{x}}$ and $S^{0}_{X_{\an},x} \to  S^{0}_{\widehat{X}^{x}}$ are defined by taking colimits over analytic opens $U \subseteq X_{\an}$ or Zariski opens $U \subseteq X$ of the maps $\Gamma(U,S^{0}_{X_{\an}}) \to \Gamma(X, S^{0}_{\widehat{X}^{x}})$ or $\Gamma(U,S^{0}_{X}) \to \Gamma(X, S^{0}_{\widehat{X}^{x}})$ respectively. To prove the maps to $S^{0}_{\widehat{X}^{x}}$ are injective we can use the argument in \cite[Proposition 3.12]{ricolfi_savvas}, where it is proven in the algebraic case using the description of $S^{0}_{\widehat{X}^{x}}$ in \eqref{formal_s_sheaf_ses}. The map $S^{0}_{X,x} \to  S^{0}_{\widehat{X}^{x}}$ is then induced by the map $\mathcal{O}_{U,x}/ I^{2}_{x} \to \widehat{\mathcal{O}}_{U,x}/ \widehat{I}^{2}_{x}$.  We can repeat the argument of \cite{ricolfi_savvas} also in the complex analytic case because the map $\mathcal{O}_{U_{\an},x} \to \widehat{\mathcal{O}}_{U_{\an},x} \cong \widehat{\mathcal{O}}_{U,x}$ is still faithfully flat. This follows because $\mathcal{O}_{U_{\an},x}$ is still a noetherian ring despite $\mathcal{O}_{U_{\an}}(U_{\an})$ not being noetherian in general. This proves $3$. \par
    Finally, to prove that the diagram \eqref{stalk_alg_formal_an} commutes we can again consider the local models of $S$ sheaves and the commutative diagram
    \begin{equation}
\begin{tikzcd}
	& {\mathcal{O}_{U,x}/I^{2}_{x}} \\
	{\widehat{\mathcal{O}}_{U,x}/\widehat{I}^{2}_{x}} & {\mathcal{O}_{U_{\an},x}/I^{2}_{an,x}}
	\arrow[from=1-2, to=2-1]
	\arrow[from=1-2, to=2-2]
	\arrow[from=2-2, to=2-1]
\end{tikzcd}
    \end{equation}
This diagram commutes because of the fact that a completion of the algebraic functions and analytic functions at a point is the same. 
\end{proof}
\begin{proof}
Consider the map $S_{X,x} \to S_{\widehat{X}^{x}}$. Now note that an element of $S_{X,x}$ is a system of pairs $(f \in \mathcal{O}(U), \alpha \in \Gamma(U, \mathbb{L}_{U}))$ for each $U$ with $x \in U$ such that
\begin{equation}
    df =0 \quad \text{ and } d_{\mathrm{dR}}f + d \alpha = 0
\end{equation}
Now assume that restricted to $\widehat{X}^{x}$ the pair vanishes. Then the function $f$ vanishes, this will imply that there is a $U^{'}$ such that $f|_{U^{'}} = 0$. Then this implies that  $d \alpha = 0$ so $\alpha \in \HH^{-1}(U^{'}, \mathbb{L}_{U^{'}})$. Then by possibly shrinking $U^{'}$ we also get that $\alpha =0$ on $U^{'}$ since it vanishes on $\widehat{X}^{x}$.
\end{proof}
The constructions of the above lemma are functorial in the sense that for a map of schemes $f \colon X \to Y$ the following diagrams commute
    \begin{equation}\label{an_to_formal}
\begin{tikzcd}
	& {S^{0}_{X_{\an},x}} & {S^{0}_{Y_{\an},x}} \\
	& {S^{0}_{\widehat{X}^{x}}} & {S^{0}_{\widehat{Y}^{y}}} \\
	{S^{0}_{X,x}} &&& {S^{0}_{Y,y}}
	\arrow[hook, from=1-2, to=2-2]
	\arrow["{f^*}"', from=1-3, to=1-2]
	\arrow[hook', from=1-3, to=2-3]
	\arrow["{\widehat{f}}"', from=2-3, to=2-2]
	\arrow[hook, from=3-1, to=1-2]
	\arrow[hook, from=3-1, to=2-2]
	\arrow[hook', from=3-4, to=1-3]
	\arrow[hook, from=3-4, to=2-3]
\end{tikzcd}
\end{equation}
\begin{equation}\label{functoriality}
\begin{tikzcd}
	{S^{0}_{\widehat{X}^{x}}} & {S^{0}_{\widehat{Y}^{x}}} \\
	{\Gamma(X,S^{0}_{X})} & {\Gamma(Y,S^{0}_{Y})} \\
	{S^{0}_{X,x}} & {S^{0}_{Y,y}}
	\arrow[from=1-2, to=1-1]
	\arrow[from=2-1, to=1-1]
	\arrow[from=2-1, to=3-1]
	\arrow[from=2-2, to=1-2]
	\arrow[from=2-2, to=2-1]
	\arrow[from=2-2, to=3-2]
	\arrow[bend left = 60, hook', from=3-1, to=1-1]
	\arrow[bend right = 60, hook, from=3-2, to=1-2]
	\arrow[from=3-2, to=3-1]
\end{tikzcd}
\end{equation}
 The diagram \eqref{functoriality} commutes either for complex analytic spaces or schemes.
We now explain how to upgrade the previous lemma to stacks. For quotient stacks $X/G$ we will by abuse of notation denote $\Gamma(X/G, S^{0}_{X/G})$ by  $\Gamma(X, S^{0}_{X})^{G}$. Strictly speaking $S^{0}_{X}$ is not $G$-equivariant in the usual sense since we are working with the lisse-\'etale site.
\begin{lem} \label{stack_S_sheafcomm}
    Let $X$ be an Artin stack, $U \to X$ be an atlas and $x \in U$ a $\mathbb{C}$ point with stabiliser $G_{x}$.  Then there is a commutative diagram
    \begin{equation}
\begin{tikzcd}
	{S^{0}_{\widehat{X}^{\B G_x}}} & {\Gamma(X,S^{0}_{X})} \\
	{S^{0}_{\widehat{U}^{x}}} & {\Gamma(U,S^{0}_{U})} \\
	& {S^{0}_{U,x}} & {}
	\arrow[hook, from=1-1, to=2-1]
	\arrow[from=1-2, to=1-1]
	\arrow[hook, from=1-2, to=2-2]
	\arrow[from=2-2, to=2-1]
	\arrow[from=2-2, to=3-2]
	\arrow[hook, from=3-2, to=2-1]
\end{tikzcd}
    \end{equation}
    The same diagram commutes for $X$ an analytic stack.
\end{lem}
\begin{proof}
We can reduce to the case of a quotient stack with $U/G$ and a point $x \in U$ with stabiliser $G$. In particular, we have $\widehat{U/G}^{\B G} = \widehat{U}^{x}/G$. Then we immediately get the injectivity of the map $S^{0}_{\widehat{U}^{x}/G} =(S^{0}_{\widehat{U}^{x}})^{G}  \to S^{0}_{\widehat{U}^{x}}$. The action on  $S^{0}_{\widehat{U}^{x}}$ is given by the natural $G$ action on the space of closed $(-1)$-shifted two forms. The commutativity of the square follows from functoriality and the triangle commutes already from the previous lemma. The same argument works for the complex analytic case.
\end{proof}
This lemma says that given two sections $s_1$ and $s_{2} \in \Gamma(X,S^{0}_X)$ we can check if they agree at a point $x \in X$ by checking if they agree on the formal completion at that point. This follows from the lemma by the commutativity of the diagram and the injectivity of the maps. Note that for comparing sections of $\Gamma(X,S^{0}_{X})$ we can work on a subspace of $\Gamma(U,S^{0}_{U})$ by \cite[Proposition 2.54]{Joyce_dcrit}. Furthermore, even though $S^{0}_{U}$ is an \' etale sheaf the global sections are the same as the associated Zariski sheaf.
\section{Formal Exponential map} \label{exp_map_chapter}
In this chapter, we will introduce the exponential map between shifted cotangent spaces and loop spaces. We will use this technology to ultimately compare DT invariants on shifted cotangent bundles with DT invariants on loop stacks. \par
We will in fact work with the exponential map in two different settings. In formal derived algebraic geometry the exponential exists on the level of formal completions giving us a map
\begin{equation}
    \widehat{\To}^{0}[-1] X \to \widehat{\mathcal{L}}^{\mathrm{const}} X
\end{equation}
However, this only sees the formal neighborhood around the point $0$-section and constant loops. We will also work on the classical analytic level where we have a map of analytic stacks induced from the usual exponential $\exp \colon \mathfrak{g} \to G$. This map agrees with the former one on the classical truncations around $0$ and $I$. We will exploit this to show that this global map preserves d-critical structures in Section \ref{d_crit_exp_section}.
\subsection{Loop spaces, shifted tangents and HKR} \label{exp_shifted_tangents}
We first consider different versions of loop spaces of algebraic stacks.
\begin{defn}[Loop spaces]
Let $X$ be a derived stack.
\begin{enumerate}
    \item Define the \textit{loop stack} to be $\mathcal{L} X = \Map(S^1, X) = X \times_{X \times X}X $
    \item  Define the \textit{unipotent loop stack} to be $\mathcal{L}^{u}\mathbf{X} = \Map(\B \mathbb{G}_{a} , X)$.
\end{enumerate}
\end{defn}
There are natural maps
\begin{align*}
    \mathrm{const} \colon  X & \to \mathcal{L}X \\
    f \colon S \to X & \mapsto S^{1} \times S \xrightarrow{\mathrm{pr}_S} S \xrightarrow{f} X 
\end{align*}
and $\mathrm{ev} \colon \mathcal{L}X \to X$ given by evaluating at $1 \in S^{1}$. 
Similarly we have constant unipotent loops
\begin{equation}
    X \to \mathcal{L}^{u}X.
\end{equation}
We can consider the fibres of the evaluation map $\mathcal{L} X \to X$, which will be given by
\begin{equation}
\begin{tikzcd}
	{\Omega_{x}X} & {\mathcal{L}X} \\
	S & X
	\arrow[from=1-1, to=1-2]
	\arrow[from=1-1, to=2-1]
	\arrow[from=1-2, to=2-2]
	\arrow["x", from=2-1, to=2-2]
\end{tikzcd}
\end{equation}
Where $\Omega_{x} X$ is the \textit{based loop space } of $x \colon S \to X$. This is a derived group scheme over $S$. We can view this as the derived version of the stabilizer of $x \in X(\mathbb{C})$. In particular, the classical truncation $t_{0} \Omega_{x}X = G_{x}$ is given by the stabiliser of $x$. Indeed we have $t_{0} \mathcal{L}X = IX$ the inertia stack of $X$. Furthermore, following \cite[Proposition 3.1.13]{chen2018localization} we can identify the  $S$ points of $\mathcal{L}X$ with pairs of $S$-points of $X$ and an $S$-point of $\Omega_{x}X$.
\begin{ex}[Examples of loop spaces]
Let $X = \B G$, then $\mathcal{L}X \cong G/G$ and $\mathcal{L}^{u}X \cong \widehat{G}^{U}/G$ with $U \subseteq G$ the unipotent cone of $G$.  See \cite[Proposition 3.1.28]{chen2018localization}.
\end{ex}
Note that since $S^{1}_{\mathrm{B}}$ is $1$-oriented we get that if $X$ is $n$-symplectic, the stack $\mathcal{L}X$ is naturally $(n-1)$-shifted symplectic.
\subsection{Shifted tangents, HKR and exponential map}
Let $X$ be a derived stack with perfect cotangent complex. It turns out we can define the shifted tangent space as a mapping stack in a similar way as the loop stack.
\begin{prop}[5.12 in \cite{naef2023torsion}]
There is a canonical $\mathbb{G}_{m}$-equivariant isomorphism $\To[-1] X \cong \Map(\B \widehat{\mathbb{G}}_{a}, X)$. The projection to the base $p \colon \To[-1] X \to X$ is identified with the evaluation map $\ev \circ a \colon  \Map(\B \widehat{\mathbb{G}}_{a},X) \to  \B \widehat{\mathbb{G}}_{a} \times \Map(\B \widehat{\mathbb{G}}_{a},X) \to X$. Here $a$ is induced from the map $\pt \to \B \widehat{\mathbb{G}}_{a}$.
\end{prop} \label{shifted_tangents_bgahat}
We then have a similar description of the $S$ points of $\To[-1]X$ as for $\mathcal{L}X$. Indeed an $S$-point of $\To[-1]X$ is given by an $S$-point $x \colon S \to X$ and an $S$-point of $\mathrm{Lie}(\Omega_{X,x})$.
\begin{ex}
    Let $X = \B G$. Then we have
    \begin{equation}
        \To[-1]\B G = \mathfrak{g}/G
    \end{equation}
    and  
    \begin{equation}
        \To^{*}[1] \B G = \mathfrak{g}^{\vee} / G.
    \end{equation}
We can see that the $2$-shifted symplectic structure on $\B G$, in particular the invariant bilinear form $\mathfrak{g}^{\vee} \cong \mathfrak{g}$ gives an identification
\begin{equation}
    \To[-1]\B G \cong \To^{*}[1] \B G.
\end{equation}
\end{ex}
\subsection{Affinization}
To compare  loops and shifted tangents we more explicitly describe the stacks $S^{1}_{\mathrm{B}}$, $\B \widehat{\mathbb{G}}_{a}$ and $\B \mathbb{G}_{a}$ and sheaves on them.
\begin{prop} \label{qcoh_comp}
   Fix a derived affine scheme $X$ We have the following identifications and properties
    \begin{enumerate}
        \item We have an identification $\QCoh (S^{1}_{B} \times X) \cong \Mod \mathbb{C}[x,x^{-1}](\Perf( X))$. The structure sheaf $\mathcal{O}_{S^{1}_{B} \times X}$ is sent to $\mathcal{O}_{X}$ with $x$ acting as identity. \cite[Proposition 5.1]{naef2023torsion}
        \item  We have an identification $\QCoh (\B \widehat{\mathbb{G}}_{a} \times X) \cong \Mod \mathbb{C}[x]  (\Perf (X))$. The structure sheaf $\mathcal{O}_{\B \widehat{\mathbb{G}}_{a} \times X }$ is sent to $\mathcal{O}_{X}$ with $x$ acting as $0$. \cite[Proposition 5.9]{naef2023torsion}
        \item  We have an identification $\QCoh (\B \mathbb{G}_{a} 
 \times X) \cong \Mod^{\mathrm{nilp}} \mathbb{C}[x] (\Perf (X))$ the category of perfect complexes with a nilpotent endomorphism. The structure sheaf $\mathcal{O}_{\B \mathbb{G}_{a} \times X}$ is sent to $\mathcal{O}_{X}$ with $x$ acting as $0$. \cite[Proposition 5.17]{naef2023torsion}
        \item Pullback on quasicoherent sheaves along the map $\B \widehat{\mathbb{G}}_{a} \times X \to \B \mathbb{G}_{a} \times X $ is identified with the inclusion $\Mod^{\mathrm{nilp}} \mathbb{C}[x] (\Perf (X)) \to \Mod \mathbb{C}[x]  (\Perf (X))$.  Pullback on quasicoherent sheaves along the map $S^{1}_{B} \times X \to \B \mathbb{G}_{a} \times X $ is identified with exponentiating the nilpotent endomorphism giving the map $\exp: \Mod^{\mathrm{nilp}} \mathbb{C}[x] (\Perf (X)) \to \Mod \mathbb{C}[x,x^{-1}] (\Perf (X)) $.
    \end{enumerate}
\end{prop}
Recall from \cite[Proposition 3.1]{ben2012loop} that there is a canonical adjunction
\begin{equation}
\begin{tikzcd}
	{\mathrm{dStk}} && {\mathrm{DGAlg}^{\op}}
	\arrow[""{name=0, anchor=center, inner sep=0}, "{\mathcal{O}(-)}", shift left=2, from=1-1, to=1-3]
	\arrow[""{name=1, anchor=center, inner sep=0}, "{\mathrm{Spec}(-)}", shift left=2, from=1-3, to=1-1]
	\arrow["\dashv"{anchor=center, rotate=-90}, draw=none, from=0, to=1]
\end{tikzcd}
\end{equation}
Now given a derived stack $X$ we define the affinization of $X$ to be  $ \mathrm{Aff}(X) \coloneqq \mathrm{Spec}(\mathcal{O}(X))$. The following proposition will be important for us
\begin{prop}
    Consider ${S}^{1}_{\mathrm{B}}$ and $\B \widehat{\mathbb{G}}_{a}$ then we have
    \begin{enumerate}
        \item $\mathcal{O}(S^{1}_{\mathrm{B}}) \cong \mathcal{O}(\B \widehat{\mathbb{G}}_{a}) \cong \mathbb{C}[\eta]$ with $\eta$ in degree $1$.
        \item $\mathrm{Aff}(S^{1}_{\mathrm{B}}) \cong \mathrm{Aff}(\B \widehat{\mathbb{G}}_{a}) \cong \B \mathbb{G}_{a}$.
    \end{enumerate}
\end{prop}
\begin{proof}
    This is already explained in \cite{ben2012loop} for $S^{1}_{\mathrm{B}}$. However, we can compute functions on $\B \widehat{\mathbb{G}}_{a}$ and $S^{1}_{\mathrm{B}}$ using the resolutions
    \begin{align*}
        \mathbb{C}[x] \xrightarrow{\cdot x} \mathbb{C}[x] &\cong \mathbb{C} \\
        \mathbb{C}[x^{\pm 1}] \xrightarrow{\cdot (x-1)} \mathbb{C}[x^{\pm 1}] &\cong \mathbb{C} \\
    \end{align*}
    for $\B \mathbb{G}_{a}$ and $S^{1}_{\mathrm{B}}$ respectively.
\end{proof}
 Note that a map from a stack $X \to \mathbb{G}_{a}$ is a degree $0$ function $\mathcal{O}(X) \to \mathbb{C}$. Shifting, a map $X \to \B \mathbb{G}_{a}$ is a degree $1$ function $\mathcal{O}(X) \to \mathbb{C}$. Let us explicitly describe the maps $S^{1} \to \B \mathbb{G}_{a}$ and $\B \widehat{\mathbb{G}}_{a} \to \B \mathbb{G}_{a}$. In particular, using the concrete models above we see that the maps are given by
 \begin{align*}
     \mathcal{O}(\B \widehat{\mathbb{G}}_{a}) = \mathbb{C}[\eta] & \to \mathbb{C}[1] \\
     \mathcal{O}(S^{1}_{B}) = \mathbb{C}[\eta] & \to \mathbb{C}[1] \\
     \eta & \mapsto 1
 \end{align*}
 Note that this map is exactly the dual of the fundamental class of $S^{1}_{\mathrm{B}}$ and $\B \widehat{\mathbb{G}}_{a}$.
 \begin{prop} 
    The stack $\B \widehat{\mathbb{G}}_{a}$ is $\mathcal{O}$-compact and  $1$-oriented.
\end{prop}
\begin{proof}
The stack $\B \widehat{\mathbb{G}}_{a}$ is $\mathcal{O}$-compact since it satisfies assumption 1.15 in \cite{naef2023torsion}. \par
There is an orientation since for unimodular algebraic groups $G$ of dimension $d$ there is a fundamental class in Lie algebra cohomology of dimension $d$. In particular, $\mathbb{G}_{a}$ is unimodular since it is abelian and hence $\QCoh \B \widehat{ \mathbb{G}}_{a}$ is $1$-oriented. See also Example \ref{fund_classes example}.
\end{proof}
Using this proposition we see that if $X$ is $n$-shifted symplectic, the stack $\To[-1]X$ has a natural $(n-1)$-shifted symplectic structure.
\subsection{Exponential map and shifted symplectic structures}
Using the canonical maps 
\begin{equation} \label{cospan_bg}
    \B \mathbb{Z} \to \B \mathbb{G}_{a} \xleftarrow{} \B \widehat{\mathbb{G}}_{a}
\end{equation}
we get the following correspondence
\begin{equation}
    \To[-1]X \xleftarrow{} \mathcal{L}^{u}X \xrightarrow{} \mathcal{L}X.
\end{equation}
\begin{thm} \cite[Theorem 6.9]{ben2012loop}
    Let $X$ be a derived artin stack. Then there is a formal map between the $-1$ shited tangent stack and loop stack
\begin{align*}
    \exp: \widehat{\To}^{0_{X}}[-1] X & \to \widehat{\mathcal{L}}^{\mathrm{const}}X \\
\end{align*}
 fitting into the commutative diagram
\begin{equation} \label{exp_map_diagram}
\begin{tikzcd}
	{\widehat{\To}^{0_{X}}[-1]X} && {\widehat{\mathcal{L}}^{\mathrm{const}}X} \\
	& {\widehat{\mathcal{L}}^{u,\mathrm{const}}X} \\
	& {\mathcal{L}^{u}X} \\
	{\To[-1]X} && {\mathcal{L}X}
	\arrow["\exp"', from=1-1, to=1-3]
	\arrow["\exp"', from=1-1, to=2-2]
	\arrow[from=1-1, to=4-1]
	\arrow[shift right=2, from=1-3, to=2-2]
	\arrow[from=1-3, to=4-3]
	\arrow["\cong"{description}, shift right=2, from=2-2, to=1-3]
	\arrow[from=2-2, to=3-2]
	\arrow["{q_{a}}", from=3-2, to=4-1]
	\arrow["{q_{m}}"', from=3-2, to=4-3]
\end{tikzcd}
\end{equation}
\end{thm}
\begin{ex}
To make the relationship to the usual exponential map transparent let us consider $X = \B G$ and write out the diagram in equation \eqref{exp_map_diagram}.
\begin{equation}
\begin{tikzcd}
	{\widehat{\mathfrak{g}}^{0}/G} && {\widehat{G}^{1}/G} \\
	& {\widehat{\mathcal{L}}^{u}X} \\
	& {\widehat{G}^{U}/G} \\
	{\mathfrak{g}/G} && {G/G}
	\arrow["\exp"', from=1-1, to=1-3]
	\arrow["\exp"', from=1-1, to=2-2]
	\arrow[from=1-1, to=4-1]
	\arrow[shift right=2, from=1-3, to=2-2]
	\arrow[from=1-3, to=4-3]
	\arrow["\cong"{description}, shift right=2, from=2-2, to=1-3]
	\arrow[from=2-2, to=3-2]
	\arrow[from=3-2, to=4-1]
	\arrow[from=3-2, to=4-3]
\end{tikzcd}
\end{equation}
\end{ex}
 We will want to compare induced closed forms using the diagram \eqref{exp_map_diagram}.
\begin{prop} \label{preorientation_cospan}
    The cospan $\B \widehat{\mathbb{G}}_{a} \xrightarrow{j_a}\B \mathbb{G}_{a} \xleftarrow{j_{m}}S^{1}_{\mathrm{B}} $ is  preoriented.
\end{prop}
\begin{proof}
Recall the definition of preorientations in \ref{orientations_csh}. We get two induced pre-orientations on $\B \mathbb{G}_{a}$. To prove that this cospan is preoriented we need to show that they agree. In other words we need the following diagram to commute
\begin{equation}
\begin{tikzcd}
	{p_{\B \mathbb{G}_{a}*} \mathcal{O}_{\B \mathbb{G}_{a}}} & {p_{\B \mathbb{G}_{a}*}j_{a*}j_{a}^{*} \mathcal{O}_{\B \mathbb{G}_{a}}} & {p_{\B \widehat{\mathbb{G}}_{a}*} \mathcal{O}_{\B \widehat{\mathbb{G}}_{a}}} \\
	{p_{\B \mathbb{G}_{a}*}j_{m*}j^{*}_{m} \mathcal{O}_{\B \mathbb{G}_{a}}} & {p_{S^{1}_{\mathrm{B}}*} \mathcal{O}_{S^{1}_{\mathrm{B}}}} & {\mathbb{C}[-1]}
	\arrow[from=1-1, to=1-2]
	\arrow[from=1-1, to=2-1]
	\arrow["\cong", from=1-2, to=1-3]
	\arrow["{[\B \widehat{\mathbb{G}}_{a}]}", from=1-3, to=2-3]
	\arrow["\cong"', from=2-1, to=2-2]
	\arrow["{[S^{1}_{\mathrm{B}}]}"', from=2-2, to=2-3]
\end{tikzcd}
\end{equation}
Note that in this case the maps $p_{\B \mathbb{G}_{a}*} \mathcal{O}_{\B \mathbb{G}_{a}} \to p_{\B \mathbb{G}_{a}*}j_{a*}j_{a}^{*} \mathcal{O}_{\B \mathbb{G}_{a}}$ and $p_{\B \mathbb{G}_{a}*} \mathcal{O}_{\B \mathbb{G}_{a}}\to p_{\B \mathbb{G}_{a}*}j_{m*}j^{*}_{m} \mathcal{O}_{\B \mathbb{G}_{a}}$ are isomorphisms, since we have $\mathrm{Aff}(S^{1}_{\mathrm{B}}) \cong \mathrm{Aff}(\B \widehat{\mathbb{G}}_{a}) \cong \B \mathbb{G}_{a}$. Therefore it suffices to compare the two maps  $[\B \widehat{\mathbb{G}}_{a}]$ and $[S^{1}_{B}]$ on the explicit model
\begin{equation}
    \mathbb{C}[\eta]  \to \mathbb{C}[-1]
\end{equation}
with $\eta$ a generator in degree $1$. In both cases the maps $[\B \widehat{\mathbb{G}}_{a}]$ and $[S^{1}_{B}]$ send $\eta \to 1$ so we are done.
\end{proof}
Now if we fix a closed $n$-shifted form $\omega$ on $X$. We have the following result.
\begin{thm} \label{exp_general_thm}
Let $X$ be a derived artin stack. Consider the span 
\begin{equation}
    \To[-1]X \xleftarrow{q_{a}} \mathcal{L}^{u}X \xrightarrow{q_{m}} \mathcal{L}X
\end{equation}
given by applying $\Map(-,X)$ to \eqref{cospan_bg}. Denote by $\omega_{a}$ the $(n-1)$-shifted form on $\To[-1] X$ and $\omega_{m}$ the $(n-1)$-shifted form on $\mathcal{L}X$. Then we have
\begin{equation} \label{unip_aksz}
    q_{a}^{*} \omega_{a} \cong q^{*}_{m} \omega_{m}. 
\end{equation}
The map $\exp \colon \widehat{\To}^{0_X}[-1]X \to  \widehat{\mathcal{L}}^{\mathrm{const}} X$
satisfies
\begin{equation}
    \exp^{*} \omega_{m} \cong \omega_{a} 
\end{equation}
where we have restricted $\omega_{-}$ to the respective formal completions.
\end{thm} 
\begin{proof}
We immediately obtain equation \eqref{unip_aksz} by applying functoriality of the AKSZ construction. In particular, apriori we have two induced $(n-1)$-shifted closed forms  $\omega_{u_{a}}$ and $\omega_{u_{m}}$ on $\mathcal{L}^{u} X$ induced by the orientation on $\B \widehat{ \mathbb{G}}_{a}$ and $S^{1}_{\mathrm{B}}$ respectively. However, by Proposition \ref{preorientation_cospan} these two induced orientations agree and therefore we get the equivalences
\begin{equation}
    q^{*}_{a} \omega_{a} \cong \omega_{u_{a}} \cong \omega_{u_{m}} \cong q^{*}_{m} \omega_{m}.
\end{equation}
The second part follows from the first by considering the diagram \eqref{exp_map_diagram}.
 The $-1$-shifted symplectic structures on $\widehat{\To}^{0_{X}}[-1]X$ and $\widehat{\mathcal{L}}^{\mathrm{const}} X$ are pulled back via the maps $i_a= \iota \circ q_a$ and $i_b= \iota \circ q_m$. 
\end{proof}
\subsection{Complex analytic situation}
Let $X$ be an algebraic derived stack which is tannakian. Then we can consider the following diagram
\begin{equation}
\begin{tikzcd}
	& {\Map_{\an}((\B \mathbb{G}_{a})_{\an},X_{\an})} \\
	& {(\Map(\B \mathbb{G}_{a},X))_{an}} \\
	{ \Map_{\an}((\B \widehat{\mathbb{G}}_{a})_{\an},X_{\an})} && { \Map_{\an}((S^{1}_{B})_{\an},X_{\an})} \\
	{(\Map(\B \widehat{\mathbb{G}}_{a},X))_{\an}} && {(\Map(S^{1}_{B},X))_{\an}}
	\arrow[curve={height=12pt}, from=1-2, to=3-1]
	\arrow[curve={height=-12pt}, from=1-2, to=3-3]
	\arrow[from=2-2, to=1-2]
	\arrow[from=2-2, to=3-1]
	\arrow[from=2-2, to=3-3]
	\arrow["\cong", from=4-1, to=3-1]
	\arrow["\cong", from=4-3, to=3-3]
\end{tikzcd}
\end{equation}
However, note that the mapping stack out of $\B \mathbb{G}_{a}$ does not behave well under analytification. Indeed, taking $X = \B \mathbb{G}_{m}$ we get the diagram
\begin{equation}
\begin{tikzcd}
	& {\mathbb{G}_{a, \an} \times \B \mathbb{G}_{m,\an}} \\
	& {\widehat{\mathbb{G}}^{1}_{m, \an} \times \B \mathbb{G}_{m, \an}} \\
	{\mathbb{G}_{a,\an} \times \B \mathbb{G}_{m,\an}} && {\mathbb{G}_{m,\an} \times \B \mathbb{G}_{m,\an}}
	\arrow[curve={height=12pt}, from=1-2, to=3-1]
	\arrow[curve={height=-12pt}, from=1-2, to=3-3]
	\arrow[from=2-2, to=1-2]
	\arrow[from=2-2, to=3-1]
	\arrow[from=2-2, to=3-3]
\end{tikzcd}
\end{equation}
However, considering that a closed $\mathbb{C}$-point of $\Map(\B \widehat{\mathbb{G}}_{a},X)$ is a map $\B \widehat{\mathbb{G}}_{a} \to X$, which is equivalent to a point of $X$ and an element $a \in \mathfrak{g}_{x}$. Since we are working in the complex analytic setting, we can exponentiate $a$ to obtain a map $\mathbb{G}_{a, \an} \to G_{x, an}$. This should then lift to a map
\begin{equation}
    \widetilde{\exp} \colon \Map_{\an}(\B (\widehat{\mathbb{G}}_{a})_{\an}, X_{\an}) \to \Map_{\an}(\B (\mathbb{G}_{a})_{\an}, X_{\an})
\end{equation}
Therefore, we could define the derived exponential as the composition
\begin{equation}
    \Map_{\an}(\B (\widehat{\mathbb{G}}_{a})_{\an}, X_{\an}) \xrightarrow{\widetilde{\exp}} \Map_{\an}(\B (\mathbb{G}_{a})_{\an}, X_{\an}) = \mathcal{L}^{u,\an}X \to \mathcal{L}X.
\end{equation}
However, since the above construction would not simplify the arguments of the following section significantly, we do not pursue this construction further. However, see Proposition \ref{derived_complex exponential} where we can construct a derived exponential in some cases using a trick.

\section{Complex analytic exponential in examples} \label{exp_examples}
In this section we will explain the construction of a complex analytic exponential map in our cases of interest.
\subsection{Base case} \label{base_case}
We establish here basic facts about the exponential map $\exp \colon \mathfrak{g} \to G$ we will use later. We define the \'etale locus.
\begin{lem}[Derivative of the exponential map] 
Let $\exp  \colon  \mathfrak{gl}_{n} \to \GL_{n}$ be the exponential map. The derivative of $\exp$ is 
\begin{equation}\label{der_exp}
    d \exp_{X} Y = \exp(X) \frac{1- \exp(-\ad_{X})}{\ad_{X}}Y 
\end{equation}
then $\exp$ is \'etale when $\ad_{X}$ is invertible. Equivalently the exponential map is \'etale for all $X \in \mathfrak{gl}_{n}$ that satisfy
\begin{equation}\label{etale_cond}
   \lambda_{i} - \lambda_j \neq 2 \pi i k  \text{ for any two eigenvalues of } X \text{ and }  k \in \mathbb{Z} \setminus 0. 
\end{equation}
\end{lem}
\begin{notation} \label{base_etale}
    We will call the condition in equation \eqref{etale_cond} the \'etale condition for an endomorphism. We write $\mathfrak{gl}^{\et}_n$ for the space of matrices in $\mathfrak{gl}_n$ that satisfy condition \eqref{etale_cond}. \par
    Taking the exponential of another reductive group $G$ we can similarly define the \'etale locus of $\mathfrak{g}$, also denoted $\mathfrak{g}^{\et}$, by taking away the locus where the derivative of $\exp$ vanishes.
\end{notation}
\begin{lem} \label{stabiliser_exp_diag}
     The exponential map $\exp \colon \mathfrak{gl}^{\et}_{n} \to \GL_n$ preserves stabilisers of diagonalisable matrices in $\mathfrak{gl}^{\et}_{n}$ under the conjugation action of $\GL_n$.
 \end{lem}
 \begin{proof}
     It is enough to check the claim on any representative of a diagonalizable matrix in the orbit under conjugation since the stabilisers are isomorphic under conjugation. Therefore, we may assume that $D \in \mathfrak{gl}^{\et}$ is diagonal with $D = \diag(\gamma_{1} I_{\lambda_{1}}, \dots , \gamma_{l} I_{\lambda_{l}})$ corresponding to some partition $\lambda$ of $n$ of length $l$ and $\gamma_{i} \neq \gamma_{j}$ for $i \neq j$. Now the stabiliser of $D$ only depends on the partition $\lambda$ and not on the values $\gamma_{i}$ so the only way it can change is if $\exp \gamma_{i} = \exp \gamma_{j}$ but this implies that
     \begin{equation}
         \gamma_{i} - \gamma_{j} = 2 \pi i k 
     \end{equation}
     for $k \neq 0$. This would contradict the condition in equation \eqref{etale_cond} so the exponential preserves stabilisers.
 \end{proof}
 \begin{lem} \label{one_variab_exp_lem}
     The commutative square
     \begin{equation}
\begin{tikzcd}
	{\mathfrak{gl}^{\et}_{n}} & {\GL_n} \\
	{\So^{n} \mathbb{G}^{\et}_{a}} & {\So^{n} \mathbb{G}_{m}}
	\arrow["\exp", from=1-1, to=1-2]
	\arrow[from=1-1, to=2-1]
	\arrow[from=1-2, to=2-2]
	\arrow["\exp"', from=2-1, to=2-2]
\end{tikzcd}
     \end{equation}
     is a pullback of complex analytic spaces. 
 \end{lem}
 \begin{proof}
     Recall that we can view $\So^{n} \mathbb{G}^{\et}_{a}$ as the space of diagonal $n \times n$ matrices up to permutation that satisfy the condition in equation \eqref{etale_cond} and $\So^{n} \mathbb{G}_{m}$ as the space of diagonal matrices up to permutation with non-zero entries. The pullback $P = \So^{n} \mathbb{G}^{\et}_{a} \times_{\So^{n} \mathbb{G}_{m}} \GL_{n}$ has points
     \begin{equation}
         \{ (D, A) \in \So^{n} \mathbb{G}^{\et}_{a} \times \GL_n \mid \exp D = A_{ss} \}
     \end{equation}
     where $A_{ss}$ is the diagonal $n \times n$ matrix that contains the eigenvalues of $A$ up to multiplicity. We can define a map
     \begin{align*}
         f \colon \mathfrak{gl}^{\et}_{n} & \to P \\
         x & \mapsto (x_{ss}, \exp x)
     \end{align*}
     where $x_{ss}$ is the diagonal matrix containing eigenvalues of $x$ up to multiplicity.  The map $f$ is surjective since the exponential map is still surjective once restricted to $\mathfrak{gl}^{\et}_{n}$. We will now show that $f$ is injective. Assume that we have $x, y \in \mathfrak{gl}^{\et}_{n}$ such that
     \begin{equation}
         (x_{ss}, \exp x) = (y_{ss}, \exp y)
     \end{equation}
Now consider Jordan-Chevalley decompositions of $x$ and $y$
\begin{equation}
    x = s_{x} + n_{x} \quad y = s_{y} + n_{y}
\end{equation}
with $s_{-}$ diagonalisable $n_{-}$ nilpotent and 
\begin{equation*}
    s_{x}n_{x} = n_{x} s_{x} \quad s_{y}n_{y} = n_{y}s_{y}.
\end{equation*}
Then we can take the exponential to get the equation
\begin{equation}
    \exp(s_{x}) \exp(n_{x})= \exp x = \exp y = \exp(s_{y}) \exp(n_{y})  
\end{equation}
here $\exp(s_{x})$ is still diagonalisable and $\exp(n_{x})$ is unipotent so we can use the uniqueness of the multiplicative Jordan-Chevalley decompositions of $\exp x = \exp y$  to conclude that 
\begin{equation} \label{ss_part_exp}
    \exp s_{x} = \exp s_{y} \quad \exp n_{x} = \exp n_{y}.
\end{equation}
The exponential map defines a bijection between the nilpotent cone and the unipotent cone which implies that $n_{x} = n_{y}$. Now because $x_{ss} = y_{ss}$ we have that there exist invertible matrices $S_{1}, S_{2}$ such that
\begin{equation*}
    S^{-1}_{1} s_{x} S_{1} = x_{ss} = y_{ss} = S^{-1}_{2} s_{y} S_{2}
\end{equation*}
Then exponentiating, using the fact that the exponential commutes with conjugation and equation \eqref{ss_part_exp} we get
\begin{equation*}
    \exp s_{x} = S_{1} S^{-1}_{2} \exp s_{x} S_{2} S^{-1}_{1}.
\end{equation*}
However, we know from Lemma \ref{stabiliser_exp_diag} that the exponential preserves stabilisers, which implies that 
\begin{align*}
    s_{x} & = S_{1} S^{-1}_{2}  s_{x} S_{2} S^{-1}_{1} \\
    S^{-1}_{1}s_{x}S_{1} & =  S^{-1}_{2}  s_{x} S_{2}  \\
    S^{-1}_{2}  s_{x} S_{2} & = S^{-1}_{2}  s_{y} S_{2}
\end{align*}
hence $s_{x} = s_{y}$ and therefore $x = y$.
\end{proof}
\subsection{(Twisted) local systems and preprojective representations} \label{twisted_preprojective}
Consider the following examples of algebras $A$
\begin{enumerate}
    \item (deformed) preprojective algebra of a quiver $Q$
    \item (twisted) fundamental group algbera $\mathbb{C} \pi_{1}(\Sigma_{g}, \omega)$ of a surface.
\end{enumerate}
Then we can consider the derived moduli stack of 
 $\mathfrak{M}_{A,d}$ of representations and $\mathrm{M}_{A,d}$ the good moduli space representing semisimple representations. It is shown in \cite[Section 6.2]{davison2022bps} that the above algebras are $2$-Calabi-Yau, therefore the stack of representations has a $0$-shifted symplectic structure. Then $\To[-1] \mathfrak{M}_{A,d}$ is the moduli space of representations $\mathfrak{M}_{A[x],d}$ equipped with an endomorphism and $\mathcal{L} \mathfrak{M}_{A,d}$ is the moduli space of representations $\mathfrak{M}_{A[x,x^{-1}],d}$ equipped with an automorphism. The good moduli spaces $\mathrm{M}_{A[x],d}$ and $\mathrm{M}(A[x,x^{-1}],d)$ parametrise semisimple representations. Note in this special case we can explicitly define an exponential map over the base
\begin{equation}
\begin{tikzcd}
	{\mathfrak{M}_{A,d}} \\
	{\mathfrak{M}_{A[x],d,\an}} & {\mathfrak{M}_{A[x,x^{-1}],d,an}} \\
	{\mathrm{M}_{A[x],d,\an}} & {\mathrm{M}_{A[x,x^{-1}],d,an}} \\
	{\mathrm{M}_{A,d,\an}}
	\arrow[from=2-1, to=1-1]
	\arrow["\exp", from=2-1, to=2-2]
	\arrow[from=2-1, to=3-1]
	\arrow[from=2-2, to=1-1]
	\arrow[from=2-2, to=3-2]
	\arrow["\exp"', from=3-1, to=3-2]
	\arrow[from=3-1, to=4-1]
	\arrow[from=3-2, to=4-1]
\end{tikzcd}
\end{equation}
by defining the exponential on representation spaces and taking quotient by $\GL_n$. This exponential is then compatible with direct sum since the exponential is compatible with block matrices.
\subsection{Derived complex analytic exponential} \label{exp_for_higgs_loc}
In this section we show that for certain moduli spaces of interest we can actually define a derived complex exponential map.
\begin{prop} \label{derived_complex exponential}
Let $X$ be a smooth projective variety and $M$ a finite homotopy type. Let $G$ be a connected reductive algebraic group.
Let $Y$ be $X$,$X_{\mathrm{dR}}$, $X_{\B}$,$X_{\mathrm{Dol}}$ or $M_{\B}$ .
    We can define a derived exponential map
\begin{equation}
    \exp \colon (\To[-1] \Map(Y, \B G))_{\mathrm{an}} \to (\mathcal{L} \Map(Y,\B G))_{\mathrm{an}}
\end{equation}
over the base $(\Map(Y,\B G))_{\an}$
Furthermore, we can obtain an \'etale map
\begin{equation}
    \exp \colon (\To[-1] \Map(Y, \B G))^{\mathrm{et}}_{\mathrm{an}} \xhookrightarrow{}(\To[-1] \Map(Y, \B G))_{\mathrm{an}} \to (\mathcal{L} \Map(Y,\B G))_{\mathrm{an}}
\end{equation}
Here $(\To[-1] \Map(Y, \B G))^{\mathrm{et}}_{\mathrm{an}}$ is defined by taking $\To[-1]$ of $\mathfrak{g}^{\mathrm{et}}/G$. 
\end{prop}
\begin{proof}
We start with the map
\begin{equation}
    \exp \colon \Map_{\an}((\B \widehat{\mathbb{G}}_{a})_{\an}, (\B G)_{\an}) \to  \Map_{\an}((S^{1}_{\B})_{\an}, (\B G)_{\an})
\end{equation}
Now we use that $\B \widehat{\mathbb{G}}_{a}$ and $\B G$ fit into the requirements for \cite[Theorem 6.13 ]{holstein2018analytification} we get that this is identified with the map
\begin{equation}
    \exp \colon (\Map(\B \widehat{\mathbb{G}}_{a}, \B G)_{\an} \to  (\Map(S^{1}_{\B}, \B G)_{\an}
\end{equation}
This defines the exponential map
\begin{equation}
    \exp \colon (\mathfrak{g}/G)_{\an} \xrightarrow{} (G/G)_{\an}.
\end{equation}
We now apply $\mathrm{Map}_{\mathrm{an}}(Y_{\mathrm{an}},-)$ to get
\begin{equation}
    \exp \colon \Map_{\mathrm{an}}(Y_{\mathrm{an}},(\mathfrak{g}/G)_{\an}) \to  \Map_{\mathrm{an}}(Y_{\mathrm{an},}(G/G)_{\an})
\end{equation}
Now note that $\mathfrak{g}/G$ and $G/G$ are Tannakian and the spaces in the condition of the proposition, satisfy the Universal GAGA property \cite[Definition 5.1, Example 5.3,Section 5.2]{holstein2018analytification}. Therefore we have that 
\begin{align*}
    \Map_{\mathrm{an}}(Y_{\mathrm{an}},(\mathfrak{g}/G)_{\an}) & \cong \Map(Y, \mathfrak{g}/G) \cong  \To[-1]\Map(Y,\B G) \\
    \Map_{\mathrm{an}}(Y_{\mathrm{an}},(G/G)_{\an}) & \cong \Map(Y, G/G) \cong  \mathcal{L}\Map(Y,\B G)
\end{align*}
 Furthermore, if we restrict to the \'etale locus as in \ref{base_etale} then taking mapping stacks preserves \'etale maps so we also obtain an \'etale map.
\end{proof}
\begin{remark}
    This proposition constructs a derived complex exponential in two of our main cases of interest namely local systems, by taking $M_{\B}$, and Higgs bundles, by taking $X_{\Dol}$. By restricting to the semistable locus we also obtain a derived complex analytic exponential for semistable Higgs bundles. 
\end{remark}
Despite being able to construct a derived complex analytic exponential in these cases, we cannot immediately deduce compatibility with good moduli spaces. Therefore, we now turn to more explicit constructions to define an exponential on good moduli spaces of the moduli problems we are interested in.
\subsection{Construction of exponential map for coherent sheaves with an endomorphism} \label{exp_for_coh}
We will now construct an exponential map  on good moduli spaces. For this, we will need some more explicit constructions, which we carry out now. \par
We first construct an exponential map for an endomorphism of coherent sheaves.
\begin{lem} \label{exp_for_coh_explicit}
    Fix an endomorphism $a$ of a coherent sheaf on $\mathcal{F}$ then there exists an automorphism $\exp a$ of $\mathcal{F}$ defined by
    \begin{equation}
        \exp a = \sum^{\infty}_{i =0} \frac{a^{n}}{n!}
    \end{equation}
    with $a^{0} = I_{\mathcal{F}}$.
\end{lem}
\begin{proof}
    We first prove that $\exp a$ converges to an automorphism of the associated analytic coherent sheaf $\mathcal{F}^{\an}$ on $X^{\an}$. We can consider $\mathcal{F}$ as as an analytic coherent sheaf on $X_{\an}$ then we have locally an exact sequence
    \begin{equation} \label{local_coherence}
        \mathcal{O}^{p}_{X^{\an}}|_{U} \to \mathcal{O}^{q}_{X^{\an}}|_{U}  \to \mathcal{F}_{\an}|_{U} \to 0
    \end{equation}
    Now note that for any complex analytic space $Y$ and a function $f \colon Y \to \mathbb{C}$ we can define an invertible function $\exp f$ by the composition
    \begin{equation}
        \exp f \colon Y \xrightarrow{f} \mathbb{C} \xrightarrow{\exp} \mathbb{C}^{*} \xhookrightarrow{} \mathbb{C}
    \end{equation}
    therefore for a morphism of a trivial bundle $\mathcal{O}^{n}_{Y} \xrightarrow{x} \mathcal{O}^{n}_{Y} $ we can define an exponential $\exp x \colon \mathcal{O}^{n}_{Y} \to \mathcal{O}^{n}_{Y}$. Now given an endomorphism $x$ of $\mathcal{F}$ we can lift the local presentation of equation \eqref{local_coherence} to a square
    \begin{equation}
\begin{tikzcd}
	{\mathcal{O}^{p}_{X_{\an}}|_{U}} & {\mathcal{O}^{q}_{X_{\an}}|_{U}} & {\mathcal{F}_{\an}|_{U}} & 0 \\
	{\mathcal{O}^{p}_{X_{\an}}|_{U}} & {\mathcal{O}^{q}_{X_{\an}}|_{U}} & {\mathcal{F}_{\an}|_{U}} & 0
	\arrow[from=1-1, to=1-2]
	\arrow["{x^{2}|_{U}}"{description}, from=1-1, to=2-1]
	\arrow[from=1-2, to=1-3]
	\arrow["{x^{1}|_{U}}", from=1-2, to=2-2]
	\arrow[from=1-3, to=1-4]
	\arrow["{x|_{U}}", from=1-3, to=2-3]
	\arrow[from=2-1, to=2-2]
	\arrow[from=2-2, to=2-3]
	\arrow[from=2-3, to=2-4]
\end{tikzcd}
    \end{equation}
    now we can use the construction of exp on trivial bundles and the usual properties of the exponential map to get an exact sequence
    \begin{equation}
\begin{tikzcd}
	{\mathcal{O}^{p}_{X_{\an}}|_{U}} & {\mathcal{O}^{q}_{X_{\an}}|_{U}} & {\mathcal{F}_{\an}|_{U}} & 0 \\
	{\mathcal{O}^{p}_{X_{\an}}|_{U}} & {\mathcal{O}^{q}_{X_{\an}}|_{U}} & {\mathcal{F}_{\an}|_{U}} & 0
	\arrow[from=1-1, to=1-2]
	\arrow["{\exp (x^{2}|_{U})}"', from=1-1, to=2-1]
	\arrow[from=1-2, to=1-3]
	\arrow["{\exp(x^{1}|_{U})}", from=1-2, to=2-2]
	\arrow[from=1-3, to=1-4]
	\arrow["{\exp(x|_{U})}", from=1-3, to=2-3]
	\arrow[from=2-1, to=2-2]
	\arrow[from=2-2, to=2-3]
	\arrow[from=2-3, to=2-4]
\end{tikzcd}
    \end{equation}
    defining the exponential of $x$ locally. Again by usual properties of the exponential this description glues to a global map $\exp x \colon \mathcal{F}_{\an} \to \mathcal{F}_{\an}$.
    Then by the GAGA correspondence, since $X$ is projective  we have that $\mathrm{End}(\mathcal{F}) \cong \mathrm{End}(\mathcal{F}^{\an})$ so that there is a corresponding algebraic morphism $\exp \colon \mathcal{F} \to \mathcal{F}$.
\end{proof}
We now introduce notation for the moduli spaces we will consider
\begin{notation}[Good moduli assumptions]
\label{moduli_space_notation}
    From now on we will use the notation $M$ for one of our moduli spaces in consideration. In particular, $M$ is one of
    \begin{enumerate}
        \item Higgs sheaves $\mathcal{H}\mathrm{iggs}^{ss}_{P}(X)$ or Higgs bundles $\mathrm{Higgs}^{ss}_{P}(X)$.
        \item Coherent sheaves $\mathrm{Coh}^{ss}_{P}(X)$ on a smooth projective variety $X$.
        \item Preprojective algebra representations $\mathcal{M}_{\Pi_{Q},d}$
        \item Twisted local systems $\mathrm{Loc}^{\omega}_{\GL_n}( \Sigma_{g})$.
    \end{enumerate}
We will then also use the following notation
\begin{enumerate}
    \item $X$ for the corresponding good moduli spaces
    \item $M_{a}$ for the classical truncation of $\To[-1]M$ and $X_{a}$ for the corresponding good moduli space
    \item $M_{m}$ for the classical truncation of $\mathcal{L}M$ and $X_{m}$ for the corresponding good moduli space.
\end{enumerate}
\end{notation}
\subsection{Quotient description of (Higgs)-coherent sheaves on a smooth projective variety}
Let $X$ be a smooth projective variety over $\mathbb{C}$. Here we follow \cite{simpson1994moduli}. I would like to thank Yukinobu Toda for pointing out the constructions of these moduli spaces in the paper \cite{simpson1994moduli}. For a coherent sheaf $\mathcal{F}$, we denote the Hilbert polynomial of $\mathcal{F}$ by $p_{\mathcal{F}}$. A $\Lambda$-module $\mathcal{F}$ is $p$-semistable if it is of pure dimension $d$ and for any sub $\Lambda$-module $\mathcal{F}^{'} \subseteq \mathcal{F}$ with $0 < r(\mathcal{F}^{'}) < r(\mathcal{F})$ we have an $N$ such that for $n \geq N$ 
\begin{equation}
    \frac{p_{\mathcal{F}(n)}}{r(\mathcal{F})} \leq \frac{p_{\mathcal{F}(n)}}{r(\mathcal{F})}
\end{equation}
\begin{remark}
    $p$-semistability implies $\mu$-semistability.
\end{remark}
In particular we consider $\Lambda$-modules for
\begin{enumerate}
\item $\Lambda = \mathcal{O}_{X}$, modules are coherent sheaves on $X$. 
    \item $\Lambda_{a} = \mathrm{Sym}(\mathcal{O}_{X})$, modules are coherent sheaves on $X$ with an endomorphism.
    \item $\Lambda_{H} = \mathrm{Sym}( \omega_{X})$, modules are Higgs sheaves on $X$.
    \item $\Lambda_{Ha} = \mathrm{Sym}(\mathcal{O}_{X} \oplus \omega_{X})$, modules are Higgs sheaves on $X$ with a compatible endomorphism.
\end{enumerate}
Fix a polynomial $P$ and let $N \geq N_{0}$ with $N_{0}$ as in \cite[Corollary 3.6 Lemma 4.2]{simpson1994moduli}. By \cite[Theorem 3.8]{simpson1994moduli}, for any $\Lambda$ we can build a quasiprojective
$Q_{\Lambda,P}$ which parametrizes
\begin{equation}
    \{ (\mathcal{F}, \alpha) \mid \alpha \colon \mathbb{C}^{p(N)} \cong \mathrm{H}^{0}(X, \mathcal{F}(N)) \}
\end{equation}
with $\mathcal{F}$ semistable, having Hilbert polynomial $P$. 
 Consider the spaces defined as 
\begin{align*}
     Q_{\Lambda_{a}, P} & \text{  space of coherent sheaf with an endomorphism} \\
     Q^{m}_{\Lambda_{a}, P} &\subseteq Q_{\Lambda_{a}, P} \text{ as the open subscheme where the endomorphism is invertible} \\
     Q_{\Lambda_{Ha}, P} & \text{  space of coherent sheaf with a Higgs field and compatible endomorphism} \\
     Q^{m}_{\Lambda_{Ha}, P} &\subseteq Q_{\Lambda_{Ha}, P} \text{ as the open subscheme where the endomorphism is invertible} \\
     Q^{\mathrm{bun}}_{\Lambda_{Ha}, P} & \text{  space of vector bundles with a Higgs field and compatible endomorphism} \\
     Q^{m,\mathrm{bun}}_{\Lambda_{Ha}, P} &\subseteq Q_{\Lambda_{Ha}, P} \text{ as the open subscheme where the endomorphism is invertible}.
\end{align*}
Then by using \cite[Theorem 4.7]{simpson1994moduli} we get the following proposition.
\begin{prop} \label{quotient_stack_coh}
    Let $X$ be a smooth projective variety. Then we have
    \begin{align*}
        Q_{\Lambda_{a}, P } / \SL_{P(N)} & \cong  t_{0}\To[-1]\mathrm{Coh}^{ss}_{P}(X) \text{ with GIT quotient } Q_{\Lambda, P} / \! \! / \SL_{P(N)} \\
        Q^{m}_{\Lambda_{a}, P } / \SL_{P(N)} & \cong  t_{0}\mathcal{L}\mathrm{Coh}^{ss}_{P}(X) \text{ with GIT quotient } Q^{m}_{\Lambda, P} / \! \! / \SL_{P(N)} \\
        Q_{\Lambda_{Ha}, P } / \SL_{P(N)} & \cong  t_{0}\To[-1]\mathcal{H}\mathrm{iggs}^{ss}_{P}(X) \text{ with GIT quotient } Q_{\Lambda_{Ha}, P} / \! \! / \SL_{P(N)} \\
        Q^{m}_{\Lambda_{Ha}, P } / \SL_{P(N)} & \cong  t_{0}\mathcal{L}\mathcal{H}\mathrm{iggs}^{ss}_{P}(X) \text{ with GIT quotient } Q^{m}_{\Lambda_{Ha}, P} / \! \! / \SL_{P(N)} \\
        Q^{bun}_{\Lambda_{Ha}, P } / \SL_{P(N)} & \cong  t_{0}\To[-1]\mathrm{Higgs}^{ss}_{P}(X) \text{ with GIT quotient } Q^{bun}_{\Lambda_{Ha}, P} / \! \! / \SL_{P(N)} \\
        Q^{m,bun}_{\Lambda_{Ha}, P } / \SL_{P(N)} & \cong  t_{0}\mathcal{L}\mathrm{Higgs}^{ss}_{P}(X) \text{ with GIT quotient } Q^{m,bun}_{\Lambda_{Ha}, P} / \! \! / \SL_{P(N)} \\
    \end{align*}
    in all cases the GIT quotient is the good moduli space of the corresponding quotient stack.
\end{prop}
\begin{remark}[Constructions for Higgs sheaves]
  We could alternatively  consider the projectivization $X = \mathbb{P}(K_{C})$ of the total space of the canonical $K_{C}$ on a smooth projective curve $C$. Then by taking the construction above and considering sheaves with support that does not cross the divisor at infinity we obtain a description of the moduli stack of semistable Higgs sheaves equipped with endomorphism on $C$.  
\end{remark}
\begin{remark}
    For the locally free case we could alternatively use the description in terms of representation spaces
    \begin{equation}
    R_{\Lambda, P}  = \{ (\mathcal{F}, \beta) \mid \mathcal{F}, \quad p \text{-semistable } \Lambda\text{-module }, \mathcal{F}_{x} \cong^{\beta} \mathbb{C}^{n}\},
\end{equation}
where the $\beta$ is viewed as a framing. This is closer to the quotient stack description of local systems. \cite[4.10]{simpson1994moduli} then constructs the desired moduli space.
\end{remark}
\subsection{Construction of exponential map on moduli spaces}
We now have a quotient presentation for all the moduli spaces in Notation \ref{moduli_space_notation}. We now can construct an exponential map on the complex analytification of the moduli spaces in Notation \ref{moduli_space_notation}.
\begin{prop}\label{exp_rep_spaces}
Let $M$ be one of the moduli spaces in Notation \ref{moduli_space_notation}.
We can define a map
\begin{equation}
   \exp \colon M_{a} \to  M_{m}.
\end{equation}
Furthermore, 
\begin{enumerate}
    \item the map is over the base $M$
    \item $\exp$ induces a commuting diagram
\begin{equation}
\begin{tikzcd}
	{(M_{a})_{\an}} & {(M_{m})_{\an}} \\
	{(X_{a})_{\an}} & {(X_{m})_{\an}}
	\arrow["\exp", from=1-1, to=1-2]
	\arrow[from=1-1, to=2-1]
	\arrow[from=1-2, to=2-2]
	\arrow["\exp", from=2-1, to=2-2]
\end{tikzcd}
\end{equation}
\item the induced map on good moduli spaces is also over the base $X$.
\item $\exp$ is compatible with direct sums.
\end{enumerate}
\end{prop}
\begin{proof}
If $M$ is the moduli space of preprojective algebra representations or twisted local systems. We can use the standard representation space descriptions to define a map in the same way as the base case in subsection \ref{base_case}. \par
    We now describe how to define the map for coherent sheaves. We  can define a map $\exp \colon  Q^{a}_{\Lambda_{\#}, P, \an } \to Q^{m}_{\Lambda_{\#}, P, \an }$ by using the exponential map construction of Lemma \ref{exp_for_coh_explicit} since the exponential preserves the underlying sheaf and framing. It is also $\SL_{P(N)}$-equivariant and therefore descends. The direct sum compatibility again follows as for matrices. The commuting diagram comes from the fact that the map
    \begin{equation}
        Q_{\Lambda_{\#}, P, \an } / \SL_{P(N),\an} \to Q_{\Lambda_{\#}, P, \an } / \! \! / \SL_{P(N),\an}
    \end{equation}
     satisfies a universal property for analytic quotients. Indeed, since these are quotients of quasiprojective schemes by linear $\SL_{P(N)}$ actions we can use \cite[Lemma 2.9, Lemma 2.13]{toda2018moduli}.
\end{proof}
\subsection{Spectral correspondence}
We can give a spectral description of the stacks of the stacks $\To[-1]$ and $\mathcal{L}$ in this case. Denote by $\mathrm{Coh}Y$ the derived stack of coherent sheaves on a variety $Y$. Now assume that $X$ is smooth projective of dimension $d$. Then by the compatibility with loops and moduli of objects proved in \cite[Theorem 9.3]{kinjo2024cohomological} \footnote{We do not give a formal proof of this statement on the derived level since we will not use it. The classical statement is more standard, see \cite{tanaka2017vafa} for an exposition.} and a similar one for $\To[-1]$ we can write
\begin{align*}
     \To[-1] \mathrm{Coh}(X)  \cong  \mathrm{Coh}_{c}(X \times \mathbb{A}^{1}) \\
     \mathcal{L} \mathrm{Coh}(X)  \cong  \mathrm{Coh}_{c}(X \times \mathbb{G}_{m}). 
\end{align*}
In the above the subscript $c$ refers to compactly supported sheaves. Similarly, for Higgs sheaves on a curve $C$
\begin{align*}
    \To[-1] \mathcal{H}\mathrm{iggs}(C)  \cong  \mathrm{Coh}_{ \leq 1}(\mathrm{Tot}_{C}(\Omega_{C}) \times \mathbb{A}^{1}) \\
     \mathcal{L} \mathcal{H}\mathrm{iggs}(C)  \cong  \mathrm{Coh}_{\leq 1}(\mathrm{Tot}_{X}(\Omega_{C}) \times \mathbb{G}_{m}) 
\end{align*}
In the above the subscript $\leq1$ refers to $1$ or $0$-dimensional compactly supported sheaves. Similarly, restricting to bundles on the left we restrict to compactly supported pure $1$ dimensional sheaves on the right.
\subsection{Jordan-Chevalley decomposition for coherent sheaves}
Using the spectral correspondence above we can view a coherent sheaf $\mathcal{F}$ on $X$ with endomorphism $x$ as a compactly supported sheaf $\widetilde{\mathcal{F}}$ on $X \times \mathbb{A}^{1}$. Then pushing forward $\widetilde{F}$ to $\mathbb{A}^{1}$ we obtain a compactly supported sheaf on $\mathbb{A}^{1}$, such a data is equivalent to a vector space $V$ with endomorphism $\tilde{x}$. This gives a way to define the characteristic polynomial of the pair $\mathcal{F},x$ on $X$ defined as the characteristic polynomial of $\tilde{x}$.
\begin{lem}[Jordan-Chevalley decomposition for sheaves] \label{jc_decomp}
Let $\mathcal{F}$ be a coherent sheaf on a smooth projective variety $X$. Then we have
\begin{enumerate}
    \item Let $x$ be an endomorphism of $\mathcal{F}$, then we can decompose $x$ uniquely as 
    \begin{equation}
        x = x_{s} + x_{n}
    \end{equation}
    such that $x_{s}$ is semisimple and $x_{n}$ is nilpotent. Furthermore, there is a decomposition
    \begin{equation}
        \mathcal{F} = \bigoplus \mathcal{F}_{\lambda_{i}}
    \end{equation}
    with $\mathcal{F}_{\lambda_{i}}$-eigensubsheaves of $\mathcal{F}$. If $x$ is invertible than we can also decompose $x = x_{s}x_{u}$ with $x_{u}$ unipotent in the sense that $I-x_{u}$ is nilpotent.
    \item If $\mathcal{F}$ is a locally free sheaf then $\mathcal{F}_{\lambda}$ are also locally free
\end{enumerate}
\end{lem}
\begin{proof}
    We can mimic the proof in linear algebra.
We can write the characteristic polynomial as
\begin{equation}
    \chi = \prod^{n}_{i=1} (t-\lambda_{i})^{c_{i}}
\end{equation}
where $\lambda_{i}$ are distinct. We can consider the products $f_{k} = \prod^{n}_{i \neq k}(t-\lambda_{i})^{c_{i}}$, then using the Chinese remainder theorem we can find polynomials $Q_{i}$ for $i=1 , \dots , n$ such that
\begin{equation}
    \id = \sum^{n}_{i =1}(Q_{i} f_{i}).
\end{equation}
Note we can deduce a global Cayley-Hamilton theorem by restrcting to the fibre at each point.
Write $P_{i}$ for $Q_{i}f_{i}$. Then we can define
\begin{equation}
    s = \sum^{n}_{i=1} \lambda_{i} P_{i}.
\end{equation}
and 
\begin{equation}
    n = x- s.
\end{equation}
We define
    \begin{equation}
        \mathcal{F}_{\lambda_{i}} = \mathrm{im} P_{i}
    \end{equation}
Note that we have
\begin{align*}
    P_{i}P_{j} & = 0 \, i \neq j \\
    P^{2}_{i} & = P_{i}.
\end{align*}
Then we can decompose 
\begin{equation}
    \mathcal{F} \cong \bigoplus^{n}_{i=1} \mathcal{F}_{\lambda_{i}}
\end{equation}
and see that $n$ is nilpotent and $s$ semisimple as in the linear algebra case. If $\mathcal{F}$ is itself a vector bundle, since the $P_i$ are idempotent, each  sheaf $\mathcal{F}_{\lambda_{i}}$ will also be a vector bundle.
\end{proof}
\subsection{Good moduli spaces and \'etale loci}
In this subsection we will establish a pullback diagram which we will use to compute pullback by the exponential map of the pushforward to the good moduli space of the loop stack. We begin by defining \'etale loci on the good moduli space level.
\begin{defn}[\'Etale locus for good moduli spaces] \label{etale_loci}
Recall Notation \ref{moduli_space_notation}. We define moduli spaces
\begin{equation}
    X^{\et}_{a} = \{ (E, F) \in (X_{a})_{\an} \mid F= \bigoplus_{i} \lambda_{i} \id_{E_{i}} \text{ satisfies  condition } \eqref{etale_cond} \} 
\end{equation}
We define the stack version 
\begin{equation}
    M^{\et}_{a} = (M_{a})_{\an} \times_{(X_{a})_{\an}} X^{\et}_{a}
\end{equation}
by pullback.
\end{defn}
\subsection{Surjectivity of exponential map}
We will now prove that the exponential map is surjective even when restricted to the \'etale locus. We start with a Lemma.
\begin{lem} \label{preserve_stabilisers}
    The map $\mathrm{M}^{\et}_{a} \xrightarrow{j^{'}} M_{a} \xrightarrow{\exp} M_{m}$ preserves stabilisers at $\mathbb{C}$-points.
\end{lem}
\begin{proof}
    Consider the description of the $\mathbb{C}$ points of $\mathrm{M}^{\et}_{a}$. Then since the exponential is a map over the base $M$ it preserves stabilisers of the underlying object. Therefore, it is enough to check that it preserves stabilisers of the endomorphism. However, using the Jordan-Chevalley decomposition  (for vector spaces or coherent sheaves in Lemma \ref{jc_decomp}) this then follows as in Lemma \ref{stabiliser_exp_diag}.
\end{proof}
\begin{lem}  \label{surjective_etale_locus}
    Consider the situation of Proposition \ref{exp_rep_spaces}. The maps
\begin{enumerate}
    \item $M^{\et}_{a} \xrightarrow{j^{'}} M_{a} \xrightarrow{\exp} M_{m}$
    \item $X^{\et}_{a} \xrightarrow{j^{'}} X_{a} \xrightarrow{\exp} X_{m}$
\end{enumerate}
are \'etale and surjective. 
\end{lem}
\begin{proof}
We begin by noting that via the exponential we have an isomorphism of formal stacks
\begin{equation}
    \exp \colon \widehat{\To}^{\mathrm{nilp}}[-1]M \cong \widehat{\mathcal{L}}^{\mathrm{unip}}M
\end{equation}
Therefore, the exponential map is surjective on unipotent automorphism. Now if the endomorphism is semisimple, we can also easily construct a preimage for the exponential. We can then use the Jordan-Chevalley decomposition from Lemma \ref{jc_decomp} to conclude that the analytic exponential map is surjective. Note we can restrict to the \'etale locus and still have surjectivity since the eigenvalues we have removed can all be replaced by the $0$ eigenvalue to get the same automorphism.\par
We now turn to proving that the maps are \'etale. We prove that $\exp$ is \'etale separately in Proposition \ref{d_crit_equivalence}. \par
Now we turn to good moduli spaces. Pick a point $x \in M^{\et}_{a}$ and write $\overline{x} \in X^{\et}_{a}$.
We begin by noting that the map $\exp \colon M^{\et}_{a} \to M_{m}$ is \'etale, and preserves stabilizers by Lemma \ref{preserve_stabilisers}. Therefore, there is an isomorphism of formal stacks
\begin{equation} \label{formal_equivalence_alp}
    \widehat{M^{\et}_{a}}^{\B G_{x}} \xrightarrow{ \cong }  \widehat{M^{\et}_{m}}^{\B G_{\exp x}}
\end{equation}
Now by \cite[Theorem 1.1]{alper2012local} or the general \'etale local structure of stacks we have that 
\begin{equation}
    \widehat{X^{\et}_{a}}^{\overline{x}} \cong (\widehat{M^{\et}_{a}}^{\B G_{x}} \times_{\B G_{x}} \pt)/ \! \! / G_{x} \text{ and } \widehat{h^{\et}_{m}}^{\overline{\exp x}} = (\widehat{M^{\et}_{m}}^{\B G_{\exp x}} \times_{\B G_{\exp x}} \pt) / \! \! / G_{\exp x}
\end{equation}
but then by preservation of stabilizers and equation \eqref{formal_equivalence_alp} we get
\begin{equation}
    \widehat{X^{\et}_{a}}^{\overline{x}} \cong \widehat{X^{\et}_{m}}^{\overline{\exp x}}.
\end{equation}
Therefore, $\exp$ is \'etale at $\overline{x}$.
\end{proof}
\begin{remark}
    The presentation of the proof that $\exp$ is \'etale is slighly out of order because we use the compatibility with graded points. However, the main point is that we have that exp is compatible with direct sum and the direct sum map is \'etale at specific points. This allows us to reduce to the case of the formal exponetial along the zero section and constant loops. 
\end{remark}
\begin{prop} \label{etale_pullback_diagram}
    Consider the \'etale loci in \ref{etale_loci}. We have a pullback diagram of complex analytic stacks and spaces with the horizontal maps \'etale.
    \begin{equation}
\begin{tikzcd}
	{\mathrm{M}^{\et}_{a}} & {\mathrm{M}_{a}} & {\mathrm{M}_{m}} \\
	{\mathrm{X}^{\et}_{a}} & {\mathrm{X}_{a}} & {\mathrm{X}_{m}}
	\arrow["{j^{'}}", hook, from=1-1, to=1-2]
	\arrow["{\pi^{\et}_{a}}"', from=1-1, to=2-1]
	\arrow["\exp", from=1-2, to=1-3]
	\arrow["{\pi_{a}}", from=1-2, to=2-2]
	\arrow["{\pi_m}", from=1-3, to=2-3]
	\arrow["j", hook, from=2-1, to=2-2]
	\arrow["\exp", from=2-2, to=2-3]
\end{tikzcd}
    \end{equation}
Here $j$ and $j^{'}$ are open immersions and the leftmost square and the whole square are pullbacks.
\end{prop}
\begin{remark}
    Note that the rightmost square is not a pullback.
\end{remark}
\begin{proof}
By definition the leftmost square is a pullback. Now, consider the diagram
\begin{equation}
\begin{tikzcd}
	{\mathrm{M}^{\et}_{a}} \\
	& {\mathrm{M}^{'}_{a}} & {\mathrm{M}_{a}} & {\mathrm{M}_{m}} \\
	& {\mathrm{X}^{\et}_{a}} & {\mathrm{X}_{a}} & {\mathrm{X}_{m}}
	\arrow["f"', from=1-1, to=2-2]
	\arrow["h", curve={height=-18pt}, from=1-1, to=2-4]
	\arrow["g"', curve={height=18pt}, from=1-1, to=3-2]
	\arrow["{j^{'}}", hook, from=2-2, to=2-3]
	\arrow["{\pi^{\et}_{a}}"', from=2-2, to=3-2]
	\arrow["\exp", from=2-3, to=2-4]
	\arrow["{\pi_{a}}", from=2-3, to=3-3]
	\arrow["{\pi_m}", from=2-4, to=3-4]
	\arrow["j", hook, from=3-2, to=3-3]
	\arrow["\exp", from=3-3, to=3-4]
\end{tikzcd}
\end{equation}
where $\mathrm{M}^{'}_{a}$ is by definition the pullback of the outer square. Now since the diagram in the proposition commutes there is a map $f \colon \mathrm{H}^{\et}_{a} \to \mathrm{H}^{'}_{a}$. Now since $j \circ \exp$ is \'etale the pullback map $h$ is \'etale. Then since $j^{'} \circ \exp$ is also \'etale by the $2$ out of $3$ property $f$ is \'etale. Now we consider the map on $\mathbb{C}$-points
\begin{align*}
    f \colon \mathrm{M}^{\et}_{a}(\mathbb{C}) & \to \mathrm{M}^{'}_{a}(\mathbb{C}) \\
    (E,a) & \mapsto ((E, a_{\mathrm{ss}}),(E, \exp a))
\end{align*}
We claim this is an equivalence since the exponential map is surjective by Lemma \ref{surjective_etale_locus} and preserves stabilisers by Lemma \ref{preserve_stabilisers}. 
\end{proof}

\section{D-critical structures and orientability} \label{d_crit_exp_section}
In this section we will prove the compatibility of d critical structures under the exponential map in certain cases. We give these assumptions now
\begin{assume} \label{exp_assumptions}
    Let $X$ be a derived artin stack equipped with a $0$-shifted symplectic structure. Further assume that there are good moduli spaces for $X$, $\mathcal{L}X$ and $\To[-1]X$. Assume that there exists a map $\exp \colon (t_{0} \To[-1]X)_{\an} \to (t_{0} \mathcal{L}X)_{\an}$ that satisfies
    \begin{enumerate}
        \item completion of the map $\exp$ agrees with the truncation of the formal exponential map $\exp \colon \widehat{\To}^{0_{X}}[-1]X \to \widehat{\mathcal{L}}^{\mathrm{const}}X$.
        \item the map $\exp$ extends to graded points and there is a commutative diagram
        \begin{equation}
\begin{tikzcd}
	{(\Grad t_{0} \To[-1]X)_{\an}} & {(\Grad t_{0} \mathcal{L}X)_{\an}} \\
	{(t_{0} \To[-1]X)_{\an}} & {(t_{0} \mathcal{L}X)_{\an}}
	\arrow["\exp", from=1-1, to=1-2]
	\arrow[from=1-1, to=2-1]
	\arrow[from=1-2, to=2-2]
	\arrow["\exp", from=2-1, to=2-2]
\end{tikzcd}
        \end{equation}
        \item Given a central torus $T$ inside the stabilizer of every point of $X$. We assume an induced compatibility of exp with the induced actions of $\To[-1] \B T$ and $\mathcal{L} \B T$ via the following diagram 
            \begin{equation}
\begin{tikzcd}
	{(\B T \times \mathfrak{t} \times \To[-1]X_{\cl \an}} & {(\B T \times T \times \mathcal{L}X)_{\cl \an}} \\
	{(\To[-1]X)_{\cl \an}} & {(\mathcal{L}X)_{\cl \an}}
	\arrow["{\exp \times \exp}", from=1-1, to=1-2]
	\arrow[from=1-1, to=2-1]
	\arrow[from=1-2, to=2-2]
	\arrow["\exp", from=2-1, to=2-2]
\end{tikzcd}
    \end{equation}
    \end{enumerate}
\end{assume}
\begin{ex}[Stacks satisfying Assumption \ref{exp_assumptions}] \label{exp_assumptions_examples}
 We give a list of moduli problems satisfying the assumptions above for which our results hold. These follow from the discussion in Section \ref{exp_examples}. We also have the corresponding constructions for good moduli spaces \ref{moduli_space_notation}.
\begin{enumerate}
    \item [\ref{exp_for_higgs_loc}]The stack $\Loc_{G}(\Sigma_{g})$ for $G$ a reductive group. The mapping stack construction of the exponential map gives the required compatibilities. The \'etale locus can be described by Proposition \ref{derived_complex exponential}.
    \item [\ref{exp_for_higgs_loc}]The stack $\Higgs^{\mathrm{ss}}_{G}(C)$. The mapping stack construction of the exponential map gives the required compatibilities. The \'etale locus can be described by Proposition \ref{derived_complex exponential}.
    \item [\ref{twisted_preprojective}]The stack $\mathfrak{M}_{\Pi^{\lambda}_{Q}}$ of deformed preprojective algebra representations or the stack of twisted representations $\mathrm{Loc}^{\mathrm{tw}}_{\GL_n}(\Sigma)$. The required compatibilities follow from compatibility of the exponential map with direct sum of vector spaces.
    \item [\ref{exp_for_coh}]The stack $\mathrm{Coh}^{\mathrm{ss}}_{P}(X)$ for $X$ a symplectic surface. The compatibility of direct sum and exponential gives the required compatibilities.
\end{enumerate}
Furthermore, if given stability conditions the properties also hold for the open substacks of semistable objects. If we are viewing the moduli problems with structure group $\GL_n$, viewed as substacks of the moduli of objects of a category, we will call them of linear type. The \'etale locus in the linear type case is given in \ref{etale_loci}. In all these cases, the action of the central torus commutes with the exponential via the usual formulas of the exponential relating the addition on the Lie algebra of commuting elements and multiplication in the group.
\end{ex}
\subsection{Main theorem}
In this section, our goal is to prove the following theorem.  
\begin{thm} \label{main_dcrit_thm}
Let $X$ satisfy Assumption \ref{exp_assumptions}. Then we have that the map between classical analytic stacks
    \begin{equation}
        \exp \colon (\To[-1]X)^{\et}_{\cl\an} \to (\mathcal{L}X)_{\cl\an}
    \end{equation}
\begin{enumerate}
    \item preserves $d$-critical structures.
    \item For $X$ as in Example \ref{exp_assumptions_examples} the orientation differs by an explicit function $j_{X}$.
    \item The function $j_{X}$ has a square root for all the examples in \ref{exp_assumptions_examples} apart from coherent sheaves.
    \item For $X$ as in Example \ref{exp_assumptions_examples} and of linear type, $\exp$ is surjective.
\end{enumerate}
\end{thm}
\begin{proof}
    This is proved in two mains parts, which form their own subsections. In Proposition \ref{d_crit_equivalence} we check compatibility of $d$-critical structures and in Proposition \ref{exp_is_oriented} we will prove that the map is oriented. Surjectivity follows from Lemma \ref{surjective_etale_locus}.
\end{proof}
We now get the following corollary, after applying 
\begin{prop}\cite[Proposition 4.5]{ben2015darboux}
Let $f \colon X \to Y$ be a smooth map of oriented $d$-critical loci of relative dimension $n$. Then we have the natural isomorphism $\varphi_X \cong f^{*}[n] \varphi_{Y}$.
\end{prop}
\begin{cor}
    The above theorem implies that for $X$ satisfying Assumption \ref{exp_assumptions} such that $\exp$ is oriented. Write $\varphi_{\To^{\et}[-1]X}$ for the DT sheaf on $(\To[-1] X)_{\cl \an}$ and $\varphi_{\mathcal{LX}}$ for the DT sheaf on $\mathcal{LX}$. We have
    \begin{equation}
         \varphi_{\To^{\et}[-1]X} \cong \exp^{*}\varphi_{\mathcal{L}X}.
    \end{equation}
    By push and pull we then get a morphism
    \begin{equation}
        \mathrm{H}^{*}(\mathcal{L}X,\varphi_{\mathcal{L}X}) \to \mathrm{H}^{*}(\To^{\et}[-1]X,\varphi_{\To^{\et}[-1]X}).
    \end{equation}
\end{cor}
We start with some facts about graded points which we will need for the proof of the main Theorem \ref{main_dcrit_thm}.
\begin{remark}
    Note that in all our cases of interest Theorem \ref{main_dcrit_thm} could be proved by considering the compatibility of the direct sum map compatibility (or map from Levi version of the stack). We use the language of graded points to give a unified proof in all the cases. Let us give some remarks about the general steps of our proof
    \begin{enumerate}
        \item We can check that $d$-critical structures agree on the $0$-section using Theorem \ref{exp_general_thm}.
        \item We now can consider a point away from this locus, then there is a graded point that maps to our chosen point. By picking a graded point we have increased central rank in the terminology \cite{bu2025cohomology}. This is the step which is a generalization of compatibility with direct sums
        \item Using the higher central rank we can use translations by a torus action, to translate the problem to the zero section of some component of the stack of graded points and then again use step 1.
    \end{enumerate}
    We use a similar idea to check orientations as well.
\end{remark}
We will now give a brief reminder on graded points and prove some Lemmas about the behaviour of graded points, loops and shifted tangents. We will use these in  Subsections \ref{fully} and \ref{orientab_fully}.
\subsection{Recollection on graded points}
Let $X$ be a derived artin stack. 
Recall the stack of $n$-graded points, defined as the mapping stack
\begin{equation}
    \mathrm{Grad}^{n} X = \Map(\B \mathbb{G}^{n}_{m}, X).
\end{equation}
Consider a point $x \in X$ with a reductive stabilizer $G_{x} = t_{0} \Omega_{x}$. Then given a semisimple element $a \in \mathfrak{g}_{x} = \mathrm{Lie}(G_x)$ we can consider the following composition
\begin{equation}
    \B \mathrm{Z}^{0}(L_{x,a}) \to \B \mathrm{Z}(L_{x,a}) \to  \B L_{x,a} \to \B G_{x} \to X
\end{equation}
where $L_{x,a} = \mathrm{C}_{G_{x}}(a)$ is the Levi subgroup given by the centralizer of $a$. There is also a cocharacter $\chi_{a} \colon \mathbb{G}_{m} \to G_{x}$. Such that
\begin{equation}
    L_{x,a} = L_{\chi_{a}}.
\end{equation}
Note that $\mathrm{Z}^{0}(L_{x,a})$ is a torus of some dimension $n$. Then we also get a point in $\Grad^{n}(X)$.
\subsection{Graded points of shifted tangents and loops} \label{shifted_graded_setup}
Then similarly we have a map
\begin{align*}
    \gr_{\B \widehat{\mathbb{G}}_{a}} \colon \To[-1] \Grad^{n}X = \Grad^{n} \To[-1]X & \to \To[-1]X \\
    \gr_{S^{1}} \colon \mathcal{L} \Grad^{n}X = \Grad^{n} \mathcal{L}X & \to \mathcal{L}X
\end{align*}
We now work in the following set up let 
\begin{align*}
  (x,a) & \in \To[-1]X \\
  (x,A) & \in \mathcal{L}X
\end{align*}
for $a \in \mathfrak{t}_{x} \subseteq \mathfrak{g}_{x}$ and $A \in G_{x}$ be closed points. Here $T_{x}$ is the maximal torus of $G_{x}$. Then we have
\begin{align*}
    G_{(x,a)} &= \mathrm{C}_{G_{x}} (a)\subseteq G_{x} \\
    G_{(x,A)} &= \mathrm{C}_{G_{x}} (A)\subseteq G_{x}.
\end{align*}
Now $\mathrm{C}_{G_{x}}(a)$ is a Levi subgroup so we can consider a cocharacter $\chi_{a} \colon \mathbb{G}_{m} \to G_{x}$ such that the derivative gives the element $a \in \mathfrak{g}_{x}$.
\begin{equation}
    \mathrm{C}_{G_{x}}(a) = \mathrm{C}_{G_{x}}(\chi_{a})
\end{equation}
This cocharacter defines a point $\widetilde{\chi_a} = (\chi_{a},a) \in \To[-1] \Grad X$ with 
\begin{equation}
    \gr_{\B \widehat{\mathbb{G}}_{a}}(\widetilde{\chi_a}) = (x,a).
\end{equation}
Note this makes sense since the stabilizer of the map $\B \mathbb{G}_{m} \to \To[-1]X$ is the same as the stabilizer of the point $(x,a)$. Furthermore, we get a non-degenerate $n$-graded point from the inclusion of centre of $\mathrm{Z}^{0}( \mathrm{C}_{G_{x}}(a))$ into $\mathrm{C}_{G_{x}}(a)$. This also defines a graded point in $\Grad^{n} X$ via the same inclusion. Now taking exponential of $a$ we have inclusions
\begin{equation}
    \mathrm{Z}^{0}( \mathrm{C}_{G_{x}}(a)) \to  \mathrm{C}_{G_{x}}(\exp a)
\end{equation}
Then we get an induced point $\widetilde{\exp{a}}$ of $\mathcal{L} \mathrm{Grad} X$.
Similarly, this then defines an $n$-graded point of $\mathcal{L} \Grad^{n} X$. 
\begin{lem} \label{grad_lemma_2}
We have
    \begin{enumerate}
        \item There is an open and closed substack $(\To[-1] \Grad^{n})_{\widetilde{\chi_a}} X$ in $\To[-1] \mathrm{Grad}^{n}X$ that contains $\widetilde{\chi_a}$. 
        \item There is a corresponding open and closed substack $\Grad_{\chi_a} X$ of $\Grad X$ that contains $\chi_a$ and every stabilizer contains a copy of the torus $\mathrm{Z}^{0}( \mathrm{C}_{G_{x}}(\chi_{a}))$.
        \item There is an action of $\To[-1] \B \mathrm{Z}^{0}( \mathrm{C}_{G_{x}}(\chi_{a}))$ on $(\To[-1] \Grad)_{\widetilde{\chi_{a}}} X$ that preserves closed forms. There is a corresponding component $\mathcal{L} \Grad_{\chi_{a}} X$ of $\mathcal{L} \Grad X$ and an action of $\mathcal{L} \B \mathrm{Z}^{0}(C_{G_{x}}(\chi_{a}))$
    \end{enumerate}
\end{lem}
\begin{proof}
Note that the graded point in $\To[-1] \Grad X$ and $\Grad X$ has the same stabilizer.
This graded point is non-degenerate so by the arguments in \cite[Proposition 1.3.9]{halpernleistner2022} there is a component $(\To[-1] \Grad)_{z} X$ and a corresponding component $\Grad_{\chi_{a}} X$ such that 
\begin{equation}
    (\To[-1] \Grad)_{\widetilde{\chi}_{a}} X = \To[-1] (\Grad_{ \chi_{a}}X).
\end{equation}
now every point of $\Grad_{\chi_{a}}X$ contains a copy of $\mathrm{Z}^{0}(C_{G_{x}}(a))$ so there is an action
\begin{equation}
    \B \mathrm{Z}(C_{G_{x}}(a)) \times \Grad_{\chi_{a}}X \to \Grad_{\chi_a}X
\end{equation}
by taking $\To[-1]$ we get an induced action
\begin{equation}
    (\B \mathrm{Z}(C_{G_{x}}(a)) \times \mathrm{Lie}(\mathrm{Z}(C_{G_{x}}(a))) ) \times \To[-1]\Grad_{\chi_a}X \to \To[-1]\Grad_{\chi_a}X.
\end{equation}
The fact that the action preserves closed forms is given by Lemma \ref{center_symplectic_act}. If $X$ satisfies \ref{exp_assumptions} and $(x,a)$ is in the \'etale locus, then the component $\mathcal{L} \Grad_{\chi_{a}}$ corresponds to $(\mathcal{L} \Grad )_{\widetilde{\exp{a}}}X$.
\end{proof}
\begin{remark}
    Note that in general it is not true that there is an $\mathbb{A}^{1}$-deformation retraction from $\Grad \mathcal{L}X$ to $X$ and so we only consider a particular component corresponding to $\exp (a)$. However, every component of $\Grad X$ induces a component of $\mathcal{L} \Grad X$.
\end{remark}
\begin{lem} \label{center_symplectic_act}
    Assume that we have a derived artin stack $X$ with a $0$-shifted symplectic structure and a $\B T$ action for a torus $T$. 
    The induced actions $\B T \times \mathfrak{t} \times \To[-1]X$  and $\B T \times T \times \mathcal{L}X \to \mathcal{L}X$, restricted to a point in $a \in \mathfrak{t}$ or $A \in T$ give translations $\To[-1]X \xrightarrow{a} \To[-1]X$ and $\mathcal{L}X \xrightarrow{A} \mathcal{L}X$. The translations preserve the induced $(-1)$-shifted closed forms.
\end{lem}
\begin{proof}
For part $1$ since we have an inclusion $T \to G_{x}$ for every point we obtain an action $\B T \times X \to X$ taking mapping stacks $\Map(\B \widehat{\mathbb{G}}_{a},-)$ and $\Map(S^{1},-)$ we obtain the desired actions. Explicitly the translations are defined by
\begin{equation}
\begin{tikzcd}
	{\To[-1]X} & {\mathfrak{t} \times \To[-1]X} & {\B T \times \mathfrak{t} \times \To[-1] X} & {\To[-1]X} \\
	& X & {\B T \times X} & X
	\arrow["{a \times \id}"', from=1-1, to=1-2]
	\arrow["a", curve={height=-24pt}, from=1-1, to=1-4]
	\arrow["\ev"', from=1-1, to=2-2]
	\arrow["{1 \times \id \times \id}", from=1-2, to=1-3]
	\arrow["{\ev \circ\pr_{\To[-1]X}}", from=1-2, to=2-2]
	\arrow["\act", from=1-3, to=1-4]
	\arrow["{\ev \times \ev}", from=1-3, to=2-3]
	\arrow["\ev", from=1-4, to=2-4]
	\arrow["{1 \times \id}", from=2-2, to=2-3]
	\arrow["\id"', curve={height=18pt}, from=2-2, to=2-4]
	\arrow["\act", from=2-3, to=2-4]
\end{tikzcd}
\end{equation}
The same diagram holds for $\mathcal{L}X$ replacting $\mathfrak{t}$ with $T$.
\par
We now turn to part $2$. The proof is the same for $\To[-1]X$ and $\mathcal{L}X$ so we only write it for the former. Consider the commutative diagram
    \begin{equation} \label{diagram_translation}
\begin{tikzcd}
	{\To[-1]X} & {\B \widehat{\mathbb{G}}_{a}\times \To[-1]X} & X \\
	{\To[-1]X} & {\B \widehat{\mathbb{G}}_{a} \times \To[-1]X}
	\arrow["{ a_{x}}"', from=1-1, to=2-1]
	\arrow["\pi"', from=1-2, to=1-1]
	\arrow["\ev", from=1-2, to=1-3]
	\arrow["{\id \times a_{x}}", from=1-2, to=2-2]
	\arrow["\ev"', from=2-2, to=1-3]
	\arrow["\pi"', from=2-2, to=2-1]
\end{tikzcd}
    \end{equation}
    Then since the AKSZ symplectic form on $\To[-1]X$ is defined by pullback by $\ev^{*}$ from $X$ and then pushing forward along $\pi$ the claim that
    \begin{equation}
        a^{*}_{x} \omega_{\To[-1]X} \sim \omega_{\To[-1]X}
    \end{equation}
    follows from the fact that the triangle on the right of equation \eqref{diagram_translation} commutes and the functoriality of the pushforward of differential forms as in \cite[Remark 3.1.3]{calaque2022aksz}.
\end{proof}
\begin{lem}\label{fiber_lemma}
Let $X$ be a derived stack. Then consider  $(x,a)$ and $(x,A)$ as above. We have the following equivalences
\begin{align*}
    \mathbb{T}_{\To[-1]X,(x,a)} & = (\mathbb{T}_{X,x})^{a} \\
    \mathbb{T}_{\mathcal{L}X,(x,A)} & = (\mathbb{T}_{X,x})^{A}
\end{align*}
\end{lem}
\begin{proof}
    Consider the following diagram
    \begin{equation} \label{base_change_fiber_lemma}
\begin{tikzcd}
	{\B \widehat{\mathbb{G}}_{a}} & {\B \widehat{\mathbb{G}}_{a} \times \To[-1]X} \\
	{(x,a)} & {\To[-1]X} & X
	\arrow["{\id \times(x,a)}", from=1-1, to=1-2]
	\arrow["\pi"{description}, from=1-1, to=2-1]
	\arrow["\pi"', from=1-2, to=2-2]
	\arrow["\ev", from=1-2, to=2-3]
	\arrow[from=2-1, to=2-2]
\end{tikzcd}
    \end{equation}
    we know that
    \begin{align*}
        \mathbb{T}_{\To[-1]X,(x,a)} & = \iota^{*} \pi_{*} \ev^{*} \mathbb{T}_{X} \\
        & = \pi_{*}(\id \times \iota)^{*}\ev^{*} \mathbb{T}_{X} \text{ using base change along} \eqref{base_change_fiber_lemma} \\
        & = \pi_{*} \mathbb{T}_{X,x} \\
        & = (\mathbb{T}_{X,x})^{a}.
    \end{align*}
    In the last two steps above we view $\mathbb{T}_{X,x}$ a sheaf on $\B \widehat{\mathbb{G}}_{a}$ with the action of $a$. Then pushing forward to the point is equivalent to taking invariants of this action. The $\mathcal{L}X$ case follows in the same way replacing $\B \widehat{\mathbb{G}}_{a}$ with $S^{1}$.
\end{proof}
\begin{lem} \label{etale_grad}
Let $X$ be a derived artin stack and consider the set up in subsection \ref{shifted_graded_setup}.
\begin{enumerate}
    \item the map $\gr_{\B \widehat{\mathbb{G}}_{a}} \colon \To[-1] \Grad X \to \To[-1] X$ is formally \'etale at $\widetilde{\chi}_{a}$.
    \item Now assume that $X$ is as in \ref{exp_assumptions} and $(x,a)$ is in the \'etale locus. Then, the map $\gr_{S^{1}} \colon \mathcal{L} \Grad X \to \mathcal{L}X$ is formally \'etale at $\widetilde{\exp{a} }$, $\gr_{S^{1}}(\exp(\widetilde{\chi_{a}})) = (x, \exp a)$. 
\end{enumerate}
\end{lem}
\begin{proof}
We have a map
\begin{equation}
\begin{tikzcd}
	{\B \mathbb{G}_{m}} & {\B C_{G_{x}}(a)} & {\To[-1]X} \\
	& {\B G_{x}} & X
	\arrow["{\chi_{a}}", from=1-1, to=1-2]
	\arrow[from=1-2, to=1-3]
	\arrow[from=1-2, to=2-2]
	\arrow[from=1-3, to=2-3]
	\arrow[from=2-2, to=2-3]
\end{tikzcd}
\end{equation}
and by restricting $\mathbb{T}_{\To[-1]X}$ we obtain a $\mathbb{G}_{m}$ action on the perfect complex of vector spaces $\mathbb{T}_{\To[-1]X,(x,a)}$. Note that we can view $\mathbb{T}_{X,x}$ as a perfect complex over $ \B \mathbb{G}_{m}$ via the action of $\chi_{a}$. In particular, $\mathbb{T}_{X,x}$ splits into weight spaces
\begin{equation}
    \mathbb{T}_{X,x} \cong \bigoplus_{\omega \in \mathrm{X}^{*}(\mathbb{G}_{m})} \mathbb{T}_{X,x,\omega}
\end{equation}
where $\mathbb{T}_{X,x,\omega}$ is the complex acted on by weight $\omega$. Now we can consider the map 
\begin{equation}
    \B \widehat{\mathbb{G}}_{a} \to \B \mathbb{G}_{m},
\end{equation}
where we view $\B \widehat{\mathbb{G}}_{a}$ as $\B \exp \mathrm{Lie}(\mathrm{G}_m)$ as in \cite{gaitsgory2017study}. Because $a$ is semisimple, and we chose $\chi_{a}$ to induce $a$ by taking derivatives, pulling back the complex $\mathbb{T}_{X,x}$ on $\B \mathbb{G}_{m}$ we obtain the action of $a$ on $\mathbb{T}_{X,x}$ on $\B \widehat{\mathbb{G}}_{a}$ and we get a decomposition
\begin{equation}
    \mathbb{T}_{X,x} \cong \bigoplus_{\omega \in \mathfrak{t}^{*}} \mathbb{T}_{X,x,\omega}.
\end{equation}
If $\omega \neq 0$, then we can explicitly compute that $(\mathbb{T}_{X,x,\omega})^{a} \cong 0$. Therefore, only the weight $0$-part contributes to the invariants of $\mathbb{T}_{X,x}$. In particular, we have
\begin{equation}
    (\mathbb{T}_{X,x})^{a} \cong (\mathbb{T}_{X,x,0})^{a}
\end{equation}
 By Lemma \ref{fiber_lemma}, $\mathbb{T}_{\To[-1]X,(x.a)}$ is given by the invariants of the $a$ action on $\mathbb{T}_{X,x}$
\begin{equation} \label{weight_zeroonly}
    \mathbb{T}_{\To[-1]X,(x,a)} \cong (\mathbb{T}_{X,x,0})^{a}.
\end{equation}
Now recall that
\begin{equation}
    \mathbb{T}_{\To[-1]\Grad X,\widetilde{\chi_a}} \cong (\mathbb{T}_{X,x})^{\chi_{a},a} \cong  (\mathbb{T}_{\To[-1]X,(x,a)})^{\chi_{a}}
\end{equation}
By $\chi_{a,},a$ we mean the invariants under the action of both $\chi_{a}$ and $a$, viewing $\Grad \To[-1]X$ as the mapping stack $\Map(\B \mathbb{G}_{m} \times \B \widehat{\mathbb{G}}_{a},x)$. By $\chi_{a}$ we mean invariants of the character $\chi_{a}$, in other words, the weight $0$ part.
Now we can use that
\begin{align*}
    \mathbb{T}_{\To[-1] \mathrm{\Grad}X,\widetilde{\chi_a}}  & = (\mathbb{T}_{\To[-1]X,(x,a)})^{\chi_{a}} \cong ((\mathbb{T}_{X,x})^{a})^{\chi_{a}} \\
    & \cong ((\mathbb{T}_{X,x})^{ \chi_{a}})^{a} \text{ since the two invariants commute} \\
      & \cong  ((\mathbb{T}_{X,x,0}))^{a} \, \text{ since the } \chi_{a} \text{ and } a \text{ weights are the same} \\
      & = \mathbb{T}_{\To[-1]X,(x,a)} \, \text{ by equation} \eqref{weight_zeroonly}.
\end{align*}
Note that the commutation of the two invariants holds only because we picked $\chi_a$ and $a$ in our specific way in the set up \ref{shifted_graded_setup}. Therefore, since we have a quasiisomorphism of tangent complexes we can conclude that the map is formally \'etale at this point. For the second part, since we have chosen a point in the \'etale locus, the action weights of the action of $\exp a$ on $\mathbb{T}_{X,x}$ will still be the same as $\chi_{a}$ and then we can repeat the same argument but now for $S^1$.
\end{proof}
\begin{lem} \label{grad_symplectic_dcrit}
Let $X$ be equipped with a $2$-form $\omega_{X}$. Then there is an induced form $\omega_{\Grad X}$ given by AKSZ. The map $\gr \colon (\Grad X,\omega_{\Grad X}) \to (X, \omega_{X})$ preserves the closed two forms.
\end{lem}
\begin{proof}
    This is proven in \cite[Lemma 6.12]{descombes2025hyperbolic}.
\end{proof}
\subsection{Volume forms on Grad} \label{volume_form_grad}
We now prove some lemmas about orientations on $\Grad X$ and the special cases of $\Grad \To[-1]X$ and $\Grad \mathcal{L}X$.
We can write down an equivalence
\begin{equation}
    \det \mathbb{L}_{\Grad X} = \gr^{*}\det \mathbb{L}_{X} \otimes \det \mathbb{L}_{\Grad X / X}
\end{equation}
Now we have a formula
\begin{equation}
    \mathbb{L}_{\Grad X / X} = (\gr ^{*}\mathbb{L}_{X})^{> 0} \oplus (\gr ^{*}\mathbb{L}_{X})^{< 0} 
\end{equation}
Then
\begin{equation}
    \det \gr^{*} \mathbb{L}_{X} = \det \gr^{*} \mathbb{L}^{0}_{X} \otimes  \det \gr^{*} \mathbb{L}^{ < 0}_{X} \times \det \gr^{*} \mathbb{L}^{>0}_{X}
\end{equation}
Then if $X$ is $(-1)$-symplectic we can identify
\begin{equation}
    \gr^{*} \mathbb{L}^{>0}_{X} \cong \gr^{*} \mathbb{L}^{<0}_{X}[1]
\end{equation}
So we can write
\begin{equation}
    \det \mathbb{L}_{\Grad X} = \det \gr^{*} \mathbb{L}^{0}_{X} = \det \gr^{*} \mathbb{L}_{X} \otimes \det^{\otimes -2} \gr^{*} \mathbb{L}^{< 0}_{X}
\end{equation}
Then given an orientation or volume form on $X$ we get an induced orientation on $\Grad X$. We call this the orientation $o_{\Grad X}$. See also \cite[Lemma 6.13]{descombes2025hyperbolic} for this construction. \par
Let us denote $\Map(Z,X)$ by $X^{Z}$ for $Z = S^{1}, \B \widehat{\mathbb{G}}_{a}$. And by $p$ the natural map $\Map(Z,X) \to X$.
Now if we take $X^{Z}$ we have a volume form on $X^{Z}$ and also on $(\Grad X)^{Z}$. Call this orientation $o^{Z}_{\Grad X^{Z}}$. So we have two orientations $o^{Z}_{\Grad X^{Z}}$ and $o_{\Grad X^{Z}}$. We will show they are the same. Now we can write the sequence
\begin{equation} \label{p_cot_sequence}
    p^{*} \mathbb{L}_{X}  \to \mathbb{L}_{X^{Z}} \to \mathbb{L}_{p} 
\end{equation}
\subsection{Comparing volume forms} \label{compargin_volume_forms}
We now compare the orientation on $\Grad X$ for a $(-1)$-symplectic stack $X$ and the volume form orientations on $\Grad \To[-1]Y = \To[-1] \Grad Y$  and $\Grad \mathcal{L}Y = \mathcal{L} \Grad Y$. I thank Pierre Descombes for providing the arguments of this subsection.
We first note we have functorialities
\begin{equation} \label{pullback_volume}
\begin{tikzcd}
	{X^{Z}} \\
	& {X \times_{Y} Y^{Z}} & {Y^{Z}} \\
	& X & Y
	\arrow["{p_{f}}"', from=1-1, to=2-2]
	\arrow["{f_{Z}}", curve={height=-12pt}, from=1-1, to=2-3]
	\arrow["{p_{X}}"', curve={height=18pt}, from=1-1, to=3-2]
	\arrow[from=2-2, to=2-3]
	\arrow[from=2-2, to=3-2]
	\arrow["{p_{Y}}", from=2-3, to=3-3]
	\arrow["f"', from=3-2, to=3-3]
\end{tikzcd}
\end{equation}
Then using the invariant trivialization
\begin{equation} \label{left_inv_triv}
    \mathbb{L}_{p} = p^{*} \mathbb{L}_{X}[1]
\end{equation}
we can write the following commutative diagrams of fiber sequences. 
\begin{equation} \label{f_funct}
\begin{tikzcd}
	{f^{*}_{Z} \mathbb{L}_{p_{Y}}} & {\mathbb{L}_{p_{X}}} & {\mathbb{L}_{p_{f}}} \\
	{f^{*}_{Z} p^{*}_{Y}\mathbb{L}_{Y}[1]} & {p^{*}_{X}\mathbb{L}_{X}[1]} & {p^{*}_{X}\mathbb{L}_{f}[1]}
	\arrow[from=1-1, to=1-2]
	\arrow["\cong"', from=1-1, to=2-1]
	\arrow[from=1-2, to=1-3]
	\arrow["\cong"', from=1-2, to=2-2]
	\arrow["\cong", from=1-3, to=2-3]
	\arrow[from=2-1, to=2-2]
	\arrow[from=2-2, to=2-3]
\end{tikzcd}
\end{equation}
The upper horizontal fibre sequence comes from the cotangent sequence of the map $X^{Z} \xrightarrow{p_{f}} X \times_{Y} T^{Z} \to X$ and base change. The lower horizontal fibre sequence comes from shiting by $[1]$, pulling back by $p_{X}$ and commutativity of the diagram.
We also have
\begin{equation}
\begin{tikzcd}
	{\mathbb{T}_{p}} & {\mathbb{T}_{X^{Z}}} & {p^{*}\mathbb{T}_{X}} \\
	{p^{*}\mathbb{T}_{X}[-1]} && {p^{*}\mathbb{L}_{X}} \\
	{p^{*}\mathbb{L}_{X}[-1]} & {\mathbb{L}_{X^{Z}}[-1]} & {\mathbb{L}_{p}[-1]}
	\arrow[from=1-1, to=1-2]
	\arrow["\cong"', from=1-1, to=2-1]
	\arrow[from=1-2, to=1-3]
	\arrow["{\omega_{\mathrm{aksz}}}", from=1-2, to=3-2]
	\arrow["{p^{*} \omega}", from=1-3, to=2-3]
	\arrow["{p^{*} \omega[-1]}"', from=2-1, to=3-1]
	\arrow[from=2-3, to=3-3]
	\arrow[from=3-1, to=3-2]
	\arrow[from=3-2, to=3-3]
\end{tikzcd}
\end{equation}
Now apply diagram \eqref{f_funct} to the map $\gr \colon \Grad X \to X$ to get 
\begin{equation}
\begin{tikzcd}
	{\gr^{*}_{Z} \mathbb{L}_{p_{X}}} & {\mathbb{L}_{p_{\Grad X}}} & {\mathbb{L}_{p_{\gr}}} \\
	{\gr^{*}_{Z} p^{*}_{X}\mathbb{L}_{X}[1]} & {p^{*}_{\Grad X}\mathbb{L}_{\Grad X}[1]} & {p^{*}_{\Grad X}\mathbb{L}_{\gr}[1]}
	\arrow[from=1-1, to=1-2]
	\arrow["\cong"', from=1-1, to=2-1]
	\arrow[from=1-2, to=1-3]
	\arrow["\cong"', from=1-2, to=2-2]
	\arrow["\cong", from=1-3, to=2-3]
	\arrow[from=2-1, to=2-2]
	\arrow[from=2-2, to=2-3]
\end{tikzcd}
\end{equation}
Now note that $\gr_{Z} = \gr_{X^{Z}} $.
Using the cotangent sequence for the map $(\Grad X)^{Z} \to \Grad X \times_{X} X^{Z} \to  X^{Z}$ we get a fiber sequence of cotangent complexes 
\begin{equation}
     p_{\gr}^{*} \mathbb{L}_{\pi} \to \mathbb{L}_{\gr_{Z}} \to \mathbb{L}_{p_{\gr}}
\end{equation}
and the base change \eqref{pullback_volume}
\begin{equation}
    e^{*} \mathbb{L}_{\gr} \cong \mathbb{L}_{\pi} 
\end{equation}
Putting these together we get the sequence
\begin{equation}
    p^{*}_{\Grad X}\mathbb{L}_{\gr} \to \mathbb{L}_{\gr_{X^{Z}}} \to \mathbb{L}_{p_{\gr}}
\end{equation}
Now we get a commutative square
\begin{equation} \label{comm_sq_1}
\begin{tikzcd}
	{\det \mathbb{L}_{\Grad X^{Z} / X^{Z}}} & {\det(\gr^{*}_{Z}\mathbb{L}_{X^{Z}}[1])^{\neq 0}} \\
	{\det p^{*}_{\Grad X} \mathbb{L}_{\Grad X / X} \otimes \det \mathbb{L}_{p_{\gr X}}} & {\det(\gr^{*}_{Z} p^{*}_{X}\mathbb{L}_{X}[1])^{\neq 0} \otimes \det(\gr^{*}_{Z}\mathbb{L}^{\neq 0}_{p_{X}}[1])} \\
	{\det p^{*}_{\Grad X} \mathbb{L}_{\Grad X / X} \otimes \det (p^{*}_{\Grad X} \mathbb{L}_{\Grad X / X}[1])} & {\det(\gr^{*}_{Z} p^{*}_{X}\mathbb{L}_{X}[1])^{\neq 0} \otimes \det(\gr^{*}_{Z}p^{*}_{X}\mathbb{L}^{\neq 0}_{X}[2])} \\
	{\mathcal{O}_{\Grad X^{Z}}} & {\gr^*_{X^{Z}} \mathcal{O}_{X^{Z}}}
	\arrow["\cong"', from=1-1, to=1-2]
	\arrow["{\cong \eqref{p_cot_sequence}}", from=1-1, to=2-1]
	\arrow["{\cong\eqref{p_cot_sequence}}", from=1-2, to=2-2]
	\arrow["\cong"', from=2-1, to=2-2]
	\arrow["{\cong \eqref{left_inv_triv}}", from=2-1, to=3-1]
	\arrow["{\cong \eqref{left_inv_triv}}", from=2-2, to=3-2]
	\arrow["\cong"', from=3-1, to=3-2]
	\arrow["\cong", from=3-1, to=4-1]
	\arrow["\cong"', from=3-2, to=4-2]
	\arrow["\cong", from=4-1, to=4-2]
\end{tikzcd}
\end{equation}
\begin{equation}\label{comm_sq_2}
\begin{tikzcd}
	{\det ((\gr^{*}_{Z} \mathbb{L}_{X^{Z}})^{>0}[1])} & {\det (((\gr^{*}_{Z} \mathbb{L}_{X^{Z}})^{<0})^{\vee}[2])} \\
	{\det ((\gr^{*}_{Z} p^{*}_{X}\mathbb{L}_{X})^{>0}[1]) \otimes \det ((\gr^{*}_{Z}\mathbb{L}_{p_{X}})^{>0}[1]) } & \begin{array}{c} \substack{\det ((\gr^{*}_{Z} p^{*}_{X}\mathbb{L}_{X})^{<0})^{\vee}[2]) \otimes\\\det ((\gr^{*}_{Z}\mathbb{L}_{p_{X}}))^{<0})^{\vee}[2])} \end{array} \\
	{\det (\gr^{*}_{Z} p^{*}_{X}\mathbb{L}_{X})^{>0}[1]) \otimes \det ((\gr^{*}_{Z}p^{*}_{X}\mathbb{L}_{X}))^{>0}[2])} & \begin{array}{c} \substack{\det ((\gr^{*}_{Z} p^{*}_{X}\mathbb{L}_{X})^{<0})^{\vee}[2]) \otimes\\\det ((\gr^{*}_{Z}p^{*}_{X}\mathbb{L}_{X}))^{<0})^{\vee}[3])} \end{array} \\
	{\gr^{*}_{Z}\mathcal{O}_{X^{Z}}} & {\gr^{*}_{Z}\mathcal{O}_{X^{Z}}}
	\arrow["{\gr^{*}_{Z} \omega_{\mathrm{aksz}}}"', from=1-1, to=1-2]
	\arrow["{\cong \eqref{p_cot_sequence}}", from=1-1, to=2-1]
	\arrow["{\cong \eqref{p_cot_sequence}}", from=1-2, to=2-2]
	\arrow["{\cong \eqref{left_inv_triv}}", from=2-1, to=3-1]
	\arrow["{\cong \eqref{left_inv_triv}}", from=2-2, to=3-2]
	\arrow["{\gr^{*}_{Z} p^{*}_{X}\omega_{X}}"', from=3-1, to=3-2]
	\arrow["\cong", from=3-1, to=4-1]
	\arrow["\cong", from=3-2, to=4-2]
	\arrow["\cong", from=4-1, to=4-2]
\end{tikzcd}
\end{equation}
putting together the last two diagrams we get
\begin{equation}
\begin{tikzcd}
	{\det \mathbb{L}_{\Grad X^{Z} / X^{Z}}} & {\mathcal{O}_{\Grad X^{Z}}} \\
	{\det((\gr^{*}_{Z}\mathbb{L}_{X^{Z}})^{<0}[1])\otimes \det ((\gr^{*}_{Z}\mathbb{L}_{X^{Z}})^{>0}[1])} & {\gr^*_{Z} \mathcal{O}_{X^{Z}}} \\
	{\det((\gr^{*}_{Z}\mathbb{L}_{X^{Z}})^{<0}[1])\otimes \det (((\gr^{*}_{Z}\mathbb{L}_{X^{Z}})^{\vee})^{<0}[2])} & {\gr^*_{Z} \mathcal{O}_{X^{Z}}} \\
	{\det((\gr^{*}_{Z}\mathbb{L}_{X^{Z}})^{<0}[1])^{\otimes 2}} & {\gr^*_{Z} \mathcal{O}_{X^{Z}}}
	\arrow["\cong", from=1-1, to=1-2]
	\arrow["\cong"', from=1-1, to=2-1]
	\arrow["{\eqref{comm_sq_1}}"{description}, draw=none, from=1-2, to=2-1]
	\arrow["\cong"', shift right, from=1-2, to=2-2]
	\arrow["\cong"', from=2-1, to=2-2]
	\arrow["{\gr^{*}_{Z} \omega_{\mathrm{aksz}}}", from=2-1, to=3-1]
	\arrow["{\eqref{comm_sq_2}}"{description}, draw=none, from=2-2, to=3-1]
	\arrow["\cong"', from=2-2, to=3-2]
	\arrow["\cong", from=3-1, to=3-2]
	\arrow["\cong", from=3-1, to=4-1]
	\arrow["\cong", from=3-2, to=4-2]
	\arrow["\cong", from=4-1, to=4-2]
\end{tikzcd}
\end{equation}
The above computations imply the following lemma.
\begin{lem} \label{grad_orientation_pd}
 Let $X$ be a $0$-shifted symplectic stack. The orientation $o^{Z}_{\Grad X^{Z}}$ coming from the group structure and the induced orientation $o_{\Grad X^{Z}}$, induced from the volume form on $X^{Z}$, coincide.
\end{lem}
\subsection{Commutative orientation data} \label{comm_orient_data}
We recall the definition of commutative orientation data in \cite[Section 10.2.8]{bu2025cohomology}. We will use this in subsection \ref{loop_naht_proof}.
\begin{lem}
    Let $N$ be a linear moduli stack as in \cite[10.2]{bu2025cohomology} and assume it is the moduli of objects of a $2$-CY category. Then the stacks $N^{Z} = \Map(Z, N)$ for $Z = \B \widehat{\mathbb{G}}_{a}$ or $S^{1}$ have commutative orientation data given by the canonical volume forms on $N^{Z}_{\gamma}$ for $\gamma \in \pi_{0}(N)$.
\end{lem}
\begin{proof}
    We follow the proof of \cite[Lemma 10.2.9]{bu2025cohomology} and also only show the lemma for the case $\gamma$,$\gamma$. To prove that we have commutative orientation data we have to construct an isomorphism
    \begin{equation} \label{lemma_equation_comm_ori}
        \sigma^{\star} o_{N^{Z}_{\gamma}} \cong o_{N^{Z}_{\gamma}} \boxtimes o_{N^{Z}_{\gamma}} \cong o_{N^{Z}_{ \gamma} \times N^{Z}_{\gamma}}
    \end{equation}
    and then show it satisfies a $\mathbb{Z}/ 2 \mathbb{Z}$-equivariance. To see the first part, note that the construction of the orientation $\sigma^{\star} o_{N^{Z}_{\gamma}}$ on $N^{Z}_{2 \gamma}$ in \cite[6.1.9]{bu2025cohomology} is analogous to the construction of the induced volume form on $\Grad X$ from $X$ in \ref{volume_form_grad}. Therefore, the same argument as in Subsection \ref{compargin_volume_forms} works to give the desired isomorphism \eqref{lemma_equation_comm_ori}. Now note that as in \cite[Lemma 10.2.9]{bu2025cohomology} the difference when swapping two cotangent chambers is a sign corresponding to the Euler pairing. Because we are working in the $2$-CY case, this is just $1$. Therefore, this matches with the swap on volume forms via the isomorphism $N^{Z}_{\gamma} \times N^{Z}_{\gamma} \cong (N_{\gamma} \times N_{\gamma})^{Z}$.
\end{proof}
\subsection{Checking d-critical structures over the zero section}
We first check $d$-critical structures along the locus of $0$-section and constant loops.
\begin{lem} \label{const_comp_lemma}
   We have that
\begin{equation} \label{s_sect_equation_thm7}
    (\exp^{\star}s_{m})_{0} = s_{a,0} ,
\end{equation}
where by $0$ is any point in $\To[-1]X$ in the zero section and $1$ is any point in the image of the constant loops of  $\mathcal{L}X$.
\end{lem}
\begin{proof}
Write $U_{a}$ for the atlas of $\To[-1]X$ and $U_{m}$ for the atlas of $\mathcal{L}X$. Using Corollary \ref{exp_general_thm} we can see that the exponential preserves the $(-1)$-shifted closed $2$-forms. We will now deduce the statement by using the results of Section \ref{s0_formal_sect} to perform a chase along the following diagram. In particular, we repeatedly use diagrams \eqref{an_to_formal} and \eqref{functoriality}   \par
\begin{equation}
\begin{tikzcd}
	{\mathcal{A}^{2,cl}(\To[-1]X,-1)} & {\Gamma(t_{0}\To[-1]X,S^{0}_{{t_{0}\To[-1]X}})} & {\Gamma(U_{a},S^{0}_{U_{a}})} & {S^{0}_{U_{a,\an},0}} \\
	{\mathcal{A}^{2,cl}(\widehat{\To}^{\B G_{0}}[-1]X ,-1)} & {S^{0}_{t_{0}\widehat{\To}^{\B G_{0}}[-1]X}} & {S^{0}_{\widehat{U}_{a}}} & {S^{0}_{U_{a,\an},0}} \\
	{\mathcal{A}^{2,cl}(\widehat{\mathcal{L}}^{\B G_1}X,-1)} & {S^{0}_{t_{0}\widehat{\mathcal{L}}^{\B G_1}X}} & {S^{0}_{\widehat{U}_{m}^{1}}} & {S^{0}_{U_{m,\an},1}} \\
	{\mathcal{A}^{2,cl}(\mathcal{L}X,-1)} & {\Gamma(t_{0}\mathcal{L}X,S^{0}_{{t_{0}\mathcal{L}X}})} & {\Gamma(U_{m},S^{0}_{U_{m}})} & {S^{0}_{U_{m},1}}
	\arrow[from=1-1, to=1-2]
	\arrow[from=1-1, to=2-1]
	\arrow[hook, from=1-2, to=1-3]
	\arrow[from=1-2, to=2-2]
	\arrow[from=1-3, to=1-4]
	\arrow[from=1-3, to=2-3]
	\arrow[hook', from=1-4, to=2-4]
	\arrow[from=2-1, to=2-2]
	\arrow[hook, from=2-2, to=2-3]
	\arrow[hook', from=2-4, to=2-3]
	\arrow["\exp", from=3-1, to=2-1]
	\arrow[from=3-1, to=3-2]
	\arrow["\exp", from=3-2, to=2-2]
	\arrow[hook, from=3-2, to=3-3]
	\arrow["\exp", from=3-3, to=2-3]
	\arrow["\exp", from=3-4, to=2-4]
	\arrow[hook', from=3-4, to=3-3]
	\arrow[from=4-1, to=3-1]
	\arrow[from=4-1, to=4-2]
	\arrow[from=4-2, to=3-2]
	\arrow[hook, from=4-2, to=4-3]
	\arrow[from=4-3, to=3-3]
	\arrow[from=4-3, to=4-4]
	\arrow[hook, from=4-4, to=3-4]
\end{tikzcd}
\end{equation}
In more detail: in the first column of the above diagram we work with the derived stacks $\To[-1]X$ and $\mathcal{L}X$. We use Theorem \ref{exp_general_thm} to get that the exponential map pulls back the additive closed form to the multiplicative one. Then by commutativity of the first column of squares we can also deduce that the equation \eqref{s_sect_equation_thm7} also holds for the $S^{0}$ sheaves of the formal completions of the classical truncations at $0$ and $1$. The $d$-critical structures $s_a$ and $s_m$ on the analytifications of $U_{a}$ and $U_{m}$ are induced from algebraic ones via the vertical maps in the rightmost column. By Lemma \ref{ss_ssheaf_inj} the map 
$$S^{0}_{U_{a, \an},0} \to S^{0}_{\widehat{U}^{0}_{a}}$$ is injective and by Lemma \ref{stack_S_sheafcomm} it is then enough to check that they are the same under the exponential map by first embedding into $S^{0}_{\widehat{U}^{0}_{a}}$.
\end{proof}
\subsection{Checking d-critical structures fully} \label{fully}
We will now prove under certain assumptions that the complex analytic exponential map preserves $d$-critical structures
\begin{prop} \label{d_crit_equivalence}
Let $X$ be as in Example \ref{exp_assumptions_examples}. Then for the \'etale locus as given in \ref{etale_loci} $\exp \colon (\To[-1]X)^{\et}_{\cl \an} \to (\mathcal{L}_{X})_{\an}$ is \'etale. \par
Now let $X$ satisfy Assumption \ref{exp_assumptions}. Then we have that the map between classical analytic stacks
    \begin{equation}
        \exp \colon (\To[-1]X)^{\et}_{\cl\an} \to (\mathcal{L}X)_{\cl\an}
    \end{equation}
    preserves $d$-critical structures.
\end{prop}
\begin{proof}
    We check both statements at the level of points on formal completions. Fix a point $(x,a) \in (\To[-1]X)^{\et}_{\cl \an}$. By the Assumption \ref{exp_assumptions} and Lemma \ref{grad_lemma_2}
     we can write a diagram of  formal completions.  
        \begin{equation} \label{translation_diagram}
\begin{tikzcd}
	{\widehat{\To}^{(\chi_{a},0)}[-1] \Grad X} & {\widehat{\mathcal{L}}^{(\chi_{a},1)} \Grad X} \\
	{\widehat{\To}^{(\chi_{a},a)}[-1] \Grad X} & {\widehat{\mathcal{L}}^{(\chi_{a},\exp a)} \Grad X} \\
	{\widehat{\To}^{(x,a)}[-1]X} & {\widehat{\mathcal{L}}^{(x,\exp (a))}X}
	\arrow["\exp", from=1-1, to=1-2]
	\arrow["a", from=1-1, to=2-1]
	\arrow["{\cdot \exp a}", from=1-2, to=2-2]
	\arrow["\exp", from=2-1, to=2-2]
	\arrow["{\gr_{\B \widehat{\mathbb{G}}_{a}}}", from=2-1, to=3-1]
	\arrow["{\gr_{S^{1}}}"', from=2-2, to=3-2]
	\arrow["\exp", from=3-1, to=3-2]
\end{tikzcd}
\end{equation}
Furthermore, the vertical maps are isomorphisms and the top horizontal map is an isomorphism by the construction of the formal exponential map on the constant loops/ zero section. Therefore, the bottom horizontal map is an isomorphism as well. This shows that exp is \'etale. \par
We now turn to proving the statement about d-critical structures. We  have by Lemma \ref{grad_symplectic_dcrit} that 
\begin{equation}
    \gr_{ S^{1}}^{*}s_{\mathcal{L}X,(x,\exp a)} = s_{\mathcal{L} \Grad X, (\chi_{a}, \exp a)}
\end{equation}
 We can then use the above diagram \eqref{translation_diagram} to compute
\begin{align*}
    a^{* }\exp^{*} \gr^{*}_{S^{1}} s_{\mathcal{L}X,(x,\exp a)} & = a^{*}\exp^{*} s_{\mathcal{L} \Grad X, (\chi_{ a}, \exp a)} \\
    & = \exp^{*} (\exp a)^{*} s_{\mathcal{L} \Grad X, (\chi_{a}, \exp a)}  \\
    & = \exp^{*} s_{\mathcal{L} \Grad X, (\chi_{a}, 1)} \quad \text{Lemma \ref{center_symplectic_act}} \\
    & = s_{\To[-1] \Grad X, (\chi_{ a}, 0)} \quad \text{Lemma \ref{const_comp_lemma}}.
\end{align*}
Using the above equations we can conclude that
\begin{align*}
    a^{*}\gr^{*}_{\B \widehat{\mathbb{G}}_{a}}s_{\To[-1]X,(x,a)}  & = a^{*} s_{\To[-1] \Grad X,(\chi_{a}, a)} \quad \text{ Lemma \ref{grad_symplectic_dcrit}} \\
    & = s_{\To[-1] \Grad X,(\chi_{a}, 0)} \quad \text{ Lemma \ref{center_symplectic_act} } \\
    & = a^{* }\exp^{*} \gr^{*}_{S^{1}} s_{\mathcal{L}X,(x,\exp a)} \\
    & = a^{*} \gr^{*}_{\B \widehat{\mathbb{G}}_{a}} \exp^{*} s_{\mathcal{L}X,(x,\exp a)}.
\end{align*}
Now since the maps $\gr_{\B \widehat{\mathbb{G}_{a}}}$ and $a$ are isomorphisms they induce isomorphisms on stalks of S-sheaves so finally we have  that
\begin{equation}
    \exp^{*} s_{\mathcal{L} X, (x, \exp a)} = s_{\To[-1]X,(x,a)}.
\end{equation}
This finishes the proof.
\end{proof}
\subsection{Recollection on Atiyah classes} 
\label{exp_volume_forms_subsect}
We now begin to prove that the exponential map is oriented. We first describe how the exponential map behaves with respect to volume forms as in \cite[Section 5.1]{naef2023torsion}. Firstly, by \cite[Proposition 5.9,5.17]{naef2023torsion} we can view perfect complexes on $\B \widehat{\mathbb{G}}_{a} \times X$ and $\B \mathbb{G}_{a} \times X$ as perfect complexes on $X$ with an endomorphism and a nilpotent endomorphism respectively. There is a $\mathbb{G}_{m}$ action on $\B \widehat{\mathbb{G}}_{a}$ induced by the one on $\widehat{\mathbb{G}}_{a}$. Recall the isomorphism in Proposition \ref{shifted_tangents_bgahat}
\begin{equation}
    \To[-1]X \cong  \Map(\B \widehat{\mathbb{G}}_{a},X)  
\end{equation}
Using the isomorphism in \eqref{qcoh_comp} given a complex $F \in \QCoh(X)$ the complex $\ev^{*}F \in \QCoh (\B \widehat{\mathbb{G}}_{a}) \times \To[-1] X$ corresponds to a weight $1$ endomorphism on the complex $p^{*} F \in \QCoh(\To[-1] X)$ for $p: \To[-1]X \to X$. If $F$ is bounded above the endomorphism $p^{*}F \to p^{*}F$ is equivalent to a map $F \to F \otimes \mathbb{L}_{X}[1]$ by \cite[Theorem 2.5]{monier2021notelinearstacks}.
\begin{prop}[Atiyah class in terms of shifted tangent ]
Consider $p \colon \To[-1]X \to X$. Then for every perfect complex $F$ on $X$ the Atiyah class is equivalent to weight $1$ endomorphism
\begin{equation}
    p^{*} F \to p^{*}F.
\end{equation}
This definition then extends to bounded complexes $F \in \QCoh^{-}(X)$ be defining the Atiyah class to be a weight $1$ endomorphism
\begin{equation}
    \mathrm{at} \colon p^{*}F \to p^{*}F.
\end{equation}
\end{prop}
\begin{ex}[Classifying stack]
Consider $X = \B G$, then $\mathbb{T}_{X} = \mathfrak{g}[1]$ and $\mathbb{T}_{X}[-1] = \mathfrak{g}$ then the Atiyah class action on a perfect complex $E \in \Perf \B G \cong  \Rep G$ is the action of $\mathfrak{g}$ on the representation $E$.
\end{ex}
\subsection{Functoriality of Atiyah class} \label{funct_of_atiyah}
We record three types of functorialities of Atiyah classes. These follow from the definition of the Atiyah class we have given but see also \cite{kuhn2024atiyah}. Let $X$ and $Y$ be derived stacks.
\begin{enumerate}
    \item Given a morphism $f \colon F_{1} \to F_{2}$ of perfect complexes on $X$. Then we have a commutative square on $\To[-1]X$
    \begin{equation}
\begin{tikzcd}
	{p^{*}F_{1}} & {p^{*}F_{1}} \\
	{p^{*}F_{2}} & {p^{*}F_{2}}
	\arrow["{\mathrm{at}_{F_{1}}}", from=1-1, to=1-2]
	\arrow["{p^{*}f}"', from=1-1, to=2-1]
	\arrow["{p^{*}f}", from=1-2, to=2-2]
	\arrow["{\mathrm{at}_{F_2}}"', from=2-1, to=2-2]
\end{tikzcd}
    \end{equation}
    \item We have that $\mathrm{at}_{F^{\vee}} = - \mathrm{at}^{\vee}_{F}$.
    \item Given a map $f \colon X \to Y$ and a sheaf $F$ on $Y$, we have a commutative diagram
    \begin{equation}
\begin{tikzcd}
	{f^{*}F} & {f^{*}F \otimes\mathbb{L}_{X}[1]} \\
	& {f^{*}F \otimes f^{*}\mathbb{L}_{Y}[1]}
	\arrow["{\mathrm{at}_{f^{*}F}}", from=1-1, to=1-2]
	\arrow["{f^{*} \mathrm{at}_{F}}"', from=1-1, to=2-2]
	\arrow["{\mathrm{id} \otimes (df)^{\vee}}"', from=2-2, to=1-2]
\end{tikzcd}
    \end{equation}
\end{enumerate}
\subsection{Volume forms on loops and shifted cotangents}
Recall that we have canonical volume forms $\vol_{a}$ on $\To[-1]X$ and $\vol_{m}$ on $\mathcal{L}X$ induced by the natural group structures over $X$. See \cite[Section 5]{naef2023torsion} for the definitions.
\subsection{Comparison of volume forms on shifted tangent and cotangent} \label{orientation_remark}
In \cite[Proposition 5.16]{naef2023torsion} a volume form $\omega_{a}$ on $\To[-1] X$ is defined using the using the abelian group structure of $\To^{*}[-1] X$ relative to $X$ to get $\mathbb{L}_{\To[-1]X / X} \cong p^{*}_{X} \mathbb{L}_{X}[1]$ with $p_{X} \colon \To[-1]X \to X$. Comparing this volume form to the orientation defined in \cite[Example 2.15]{dim_red_kinjo} we see that they are defined in the same way hence the volume form on $\To^{*}[-1]X$ induces this canonical orientation as in \cite{dim_red_kinjo}.
\subsection{Difference in volume forms in terms of Atiyah classes}
Now consider the correspondence $\B \widehat{\mathbb{G}}_{a} \to \B \mathbb{G}_{a} \xleftarrow{} \B \mathbb{Z}$, which by taking mapping stacks $\Map(-,X)$ induces a correspondence 
\begin{equation} \label{unipotent_loop_corresp}
    \To[-1] X \xleftarrow{q^{a}} \mathcal{L}^{u} X \xrightarrow{q^{m}} \mathcal{L} X.
\end{equation} 
Here $L^{u}$ is called the unipotent loop space. Now we can restrict the Atiyah class to the unipotent loop space to get a nilpotent endomorphism $p^* F \to p^* F$. Therefore, we can evaluate on any invertible power series $f \in \mathbb{C}[[x]]^{\times}$, to get an automorphism $f(\at_{F})$ and an invertible function $\det f(\at_{F}) \in \mathcal{O}^{\times}_{\mathcal{L}^{u} X}(\mathcal{L}^{u} X)$. Then we have the following theorem
\begin{thm}\cite[Theorem 5.23]{naef2023torsion} \label{naef_saf}
Let $\vol_{a}$ be the natural volume form on $\To[-1] X$ and $\vol_{m}$ the natural volume form on $\mathcal{L} X$. Denote by
\begin{equation}
    j_{X} = \det (\frac{\at_{\mathbb{L}_{X}}}{\exp(\at_{\mathbb{L}_{X}}) - 1 })
\end{equation}
the function on $\mathcal{L}^{u}X$.
We have an equality of volume forms on $\mathcal{L}^{u} X$
\begin{equation}
    q^{*}_{m} \vol_{m} = q^{*}_{a} \vol_{a} \cdot j_{X}.
\end{equation}
\end{thm}
\subsection{Duflo type functions} \label{duflo}
We now give some properties of the function $j_{X}$. Ultimately we will give another proof of existence of square root in our cases of interest as given in Proposition \ref{square_root_j}. However, we believe this extra structure will be relevant for the study of this function in more general examples. Consider the case $X = \B G$. Our general case is analogous to this, but this discussion is not necessary for the proof. Let us first define the functions
\begin{align*}
    j^{L}(z) &= \frac{1-\exp(-z)}{z} \\
    j^{R}(z)  & = \frac{\exp(z)-1}{z} \\
    j(z) & = \frac{\sinh{z/2}}{z/2} \\
    \sinh{z} &= \frac{\exp(z)-\exp(-z)}{2}
\end{align*}
We have that for an endomorphism $x$
\begin{equation}
    (j^{L}(x))^{T} = j^{R}(x) \quad (j(x))^{T}=j(x)
\end{equation}
We also have the relation
\begin{equation} \label{sinh_jl}
    j^{L}(z) \exp(z) = j^{R}(z) , \quad j^{R}(z) = \exp(z/2) \frac{\sinh(z/2)}{z/2} 
\end{equation}
$j(x)$,$j^{L}(x)$,$j^{R}(x)$ are invertible if and only if $x$ has no eigenvalues of the form $2 \pi i k$ for $k \in \mathbb{Z} \setminus 0$. Now note that
\begin{equation}
    d \exp (x) = j^{R}(\mathrm{ad}_{x})
\end{equation}
The zeros of the function $\det j^{R}(\mathrm{ad}_{x})$ are exactly given by the condition \eqref{etale_cond}. Therefore, the function
\begin{equation}
    \det^{-1 }j^{R}(\mathrm{ad}_{x}),
\end{equation}
known as the Duflo function, is well defined on $\mathfrak{g}^{\et}$. Now assume that there is a nondegenerate bilinear form on $\mathfrak{g}$ then we have that
\begin{equation} \label{duflo_eq_quad}
    \det j^{R}(\mathrm{ad}_{x}) = \det j(\mathrm{ad}_{X})
\end{equation}
Now in \cite[Appendix C]{meinrenken2013clifford} it is shown that $j (\mathrm{ad}_{x})$ has a global square root. We can view the adjoint action as the Atiyah class of $\B G$ as an endomorphism on $\mathfrak{g}/G$. The nondegenerate bilinear form giving equation \eqref{duflo_eq_quad} is a consequence of the $2$-shifted symplectic structure on $\B G$. \par
Now let $X$ be a $0$-shifted symplectic stack. Write $\omega^{b}$ for the quasi-isomorphism $\mathbb{T}_{X} \to \mathbb{L}_{X}$ The Atiyah class, viewed as a map $p^{*}\mathbb{L}_{X} \to p^{*} \mathbb{L}_{X}$ must satisfy 
\begin{equation}
    \omega^{b} \mathrm{at}_{\mathbb{T}_{X}} = -\mathrm{at}_{\mathbb{L}_{X}} \omega^{b} 
\end{equation}
via the functorialities in Subsection \ref{funct_of_atiyah}. Therefore, the Atiyah class is skew with respect to the symplectic pairing.
Now on the unipotent locus we will have the equality of automorphisms via the relation in equation \eqref{sinh_jl}
\begin{equation}
    j^{R}(\mathrm{at}_{X}) = \exp(\mathrm{at}_{X}/2) \frac{\sinh(\mathrm{at}_{X}/2)}{\mathrm{at}_{X}/2} 
\end{equation}
taking determinants this becomes
\begin{equation}
    \det j^{R}(\mathrm{at}_{X}) = \exp(\mathrm{tr}(\mathrm{at}_{X}/2)) \det \frac{\sinh(\mathrm{at}_{X}/2)}{\mathrm{at}_{X}/2} 
\end{equation}
The fact that $X$ is $0$-symplectic then implies that the trace of $\mathrm{at}_{X}$ vanishes and therefore we have
\begin{equation} \label{j_L_to_sinh}
    \det j^{R}(\mathrm{at}_X) = \det \frac{\sinh(\mathrm{at}_{X}/2)}{\mathrm{at}_{X}/2}
\end{equation}
Now on the unipotent locus we can always define
\begin{equation}
    j_{X} = \det^{-1} j^{R}(\mathrm{at}_X)
\end{equation}
We will now define an analytic version of $j_{X}$ on $\To^{\et}[-1]X$.
Assume that we have a presentation of the cotangent complex of $X$ as a three term complex
\begin{equation}
    \mathbb{L}_{X} \cong  E^{-1} \to E^{0} \to E^{1}
\end{equation}
then the shifted symplectic structure gives us a pairing $E^{-1} \cong (E^{1})^{\vee}$ and a symplectic pairing on $E^{0}$.
\begin{remark}
    We consider this condition on the cotangent complex to avoid technicalities of defining an analytic $j^{R}$ function for the Atiyah class for an analytic perfect complex.
\end{remark}
\begin{remark}
    Morally, what we are doing is a derived analogue of the Duflo case \ref{duflo}. In particular, if we had a derived exponential $\exp \colon \To[-1]X \to \mathcal{L}X$, for any analytic derived stack $X$, then one should express the derivative $d \exp \colon \mathbb{T}_{\To[-1]X} \to \exp^{*}\mathbb{T}_{\mathcal{L}X}$ in terms of $j^{L}(\at_{X})$ or $j^{R}(\at_{X})$ depending on which trivialization (left or right) we are using on $\To[-1]X$. Then the determinant of $j^{L}(\at_X)$ defines a function on $\To[-1]X$ with poles exactly given by the condition \ref{etale_cond}. Considering $\To^{\et}[-1]X$ we get an invertible function. Now if $X$ is $0$-symplectic then as in \cite[Appendix C.3]{meinrenken2013clifford}
    we should get
    \begin{equation}
        \det j^{L}_{\at_X} = \det j^{R}_{\at_{X}}
    \end{equation}
    because of the compatibility of the Atiyah class and the shifted symplectic structure we hope to use  this structure to produce a global analytic square root more generally. 
\end{remark}
\subsection{Construction of square root}
We view the Atiyah class an endomorphism of the pullback of this three term complex to $\To[-1]X$. Then on the etale locus we can define an analytic automorphism of the $3$-term complex
\begin{equation}
    j^{R,\et}(\mathrm{at}_{X})
\end{equation}
then we can take the function 
\begin{equation}
    \mathrm{det} j^{R, \et}_{\mathrm{at}_{X}} = \mathrm{det} j^{R, \et}(\mathrm{at}(E^{0}))  \cdot (\mathrm{det} j^{R, \et}(\mathrm{at}(E^{1})))^{2}.
\end{equation}
The square power appears from the pairing $E^{1} \cong (E^{-1})^{\vee}$. Therefore, to show this function has a square root it is enough to show that 
$\mathrm{det} j^{R, \et}_{\mathrm{at}(E^{0})}$ has a square root.
\begin{prop}\label{square_root_j}
Let $X$ be one of the examples in \ref{exp_assumptions_examples} apart from coherent sheaves. Then
\begin{equation}
    j_{X} = \mathrm{det}^{-1} j^{R, \et}_{\mathrm{at}_{X}}
\end{equation}
on $(\To[-1]X)^{\et}_{\cl \an}$ has a square root.
\end{prop}
\begin{proof}
Note that in the cases in the theorem The fact that the function $\det^{-1}j^{R}(\mathrm{at}E^{0})$ has a square root follows from the fact that $\det E^{0}$ will already be a square of line bundles. In the case of stacks of the form $\To^{*}X$ this follows from the fibre sequence of the projection $\To^{*}X \to X$. For $\Loc_{G}(\Sigma)$ this follows by explicit computation of the cotangent complex for mapping stacks. For (deformed) preprojective algebras we can use the explicit description of the cotangent complex in \cite[Section 7]{davison2022bps}.
\end{proof}
\begin{remark} \label{coherent_sheaves_orient}
    For coherent sheaves or in general moduli of objects of $2$-Calabi-Yau categories we do not have such a global description. Therefore, we do not know how to produce a global analytic square root.
\end{remark}
\subsection{Variant: volume forms on formal completions} \label{volume_formal}
Consider a stack $X$ with a volume form $\vol$. Then there is an induced volume form $\vol_{x}$ on $\widehat{X}^{x}$ since the map $\widehat{X}^{x} \to X$ is formally \'etale. Now given a map of stacks $g \colon X \to Y$, which is \'etale at $x \in X$ we get an isomorphism on formal completions $\widehat{g} \colon \widehat{X}^{x} \to \widehat{Y}^{f(x)}$ then there is a function $f_{g,x}$ on $\widehat{X}^{x}$ such that 
\begin{equation}
    \widehat{g}^{*}\vol_{y} = f_{g,x} \vol_{x}.
\end{equation}
\subsection{Checking orientatibility in general} \label{orientab_fully}
The following Lemma follows from the fact that we have defined the volume forms on $\To[-1]X$ and $\mathcal{L}X$ using the left invariant trivialization, which uses the group structures over $X$. If there is a torus inside the stabilizer of each of the points of $X$, then we get induced translations by the same group structure which preserve these trivializations.
\begin{lem} \label{volume_form_and_translation}
As in Lemma \ref{center_symplectic_act} consider a stack $X$ with a torus $T$ in the stabilizer of each point. Then the translations $a \colon \To[-1]X \to \To[-1]X$ and $A \colon \mathcal{L}X \to \mathcal{L}X$ preserve the canonical volume forms on $\To[-1]X$ and $\mathcal{L}X$ respectively.
\end{lem}
\begin{prop} \label{exp_is_oriented}
    Let $X$ be as in Assumption \ref{exp_assumptions} and also have a global description of the cotangent complex as a three term complex of vector bundles.
    The map
    \begin{equation}
        \exp \colon (\To[-1]X)^{\et}_{\cl \an} \to (\mathcal{L}X)_{\cl \an}
    \end{equation}
    is oriented up to the function
    \begin{equation}
        j_{X} = \mathrm{det}^{-1} j^{R, \et}_{\mathrm{at}_{X}}.
    \end{equation}
\end{prop}
\begin{proof}
    Note that both $(\To[-1]X)^{\et}_{\cl \an}$ and $(\mathcal{L}X)_{\cl \an}$ are oriented $d$-critical loci with orientations $o_{a}$ and $o_{m}$ respectively. Then we have that
    \begin{equation}
        \exp^{*} o_{m} = f o_{a}
    \end{equation}
    for some invertible function $f$ on $(\To[-1]X)^{\et}_{\cl \an}$. We will show that 
    \begin{equation}
        f_{\exp} = j_{X}.
    \end{equation}
    We will prove this equality of functions by pulling back to formal completions at closed points. For a formally \'etale map of prestacks $g \colon X \to Y$, equipped with volume forms $\vol_{X}$ and $\vol_{Y}$ we denote by $f_{g,x}$ the restriction of the function $f_{g}$ that satisfies $g^{*}\vol_{Y} = f_{g} \vol_{X}$ to the formal completion $\widehat{X}^{x}$. See subsection \ref{volume_formal}.
    Consider again the diagram \eqref{translation_diagram}.
Our goal is to show that
\begin{equation}
    a^{*}\gr^{*}_{ \B \widehat{\mathbb{G}}_{a}}f_{\exp,(x,a)} =a^{*}\gr^{*}_{\B \widehat{\mathbb{G}}_{a}} j_{X,(x,a)}.
\end{equation}
This will then imply that 
\begin{equation}
    f_{\exp,(x,a)} = j_{X,(x,a)}.
\end{equation}
Because $\gr_{\B \widehat{G}_{a}}$ is oriented by Lemma \ref{grad_orientation_pd}, so it does not introduce a multiplication by an invertible function, and the diagram \eqref{translation_diagram} commutes we have
\begin{equation}
    \gr^{*}f_{\exp,(x,a)} = f_{\exp,(\chi_{a},a)} 
\end{equation}
Now because the top square in the diagram \eqref{translation_diagram} commutes and translations preserve volume forms by Lemma \ref{volume_form_and_translation} we also have
\begin{equation}
    a^{*}f_{\exp,(\chi_{a},a)} = f_{\exp,(\chi_{a},0)}
\end{equation}
Now we can use Lemma \ref{translation_orientations} to conclude that 
\begin{equation}
    f_{\exp,(\chi_{a},0)} = j_{\Grad X,(\chi_{a},0)}
\end{equation}
The right hand side is well defined since $j_{\Grad X , (\chi_{a},0)}$ is restricted from the unipotent locus.
Furthermore, we can prove that
\begin{equation} \label{atiyah_funct_jfunct}
    a^{*}\gr^{*}_{\B \widehat{\mathbb{G}}_{a}}j_{X,(x,a)} = j_{\Grad X,(\chi_{a},0)}
\end{equation}
 Equation \eqref{atiyah_funct_jfunct} then follows by functoriality of Atiyah classes as in Subsection \ref{funct_of_atiyah}. In particular, since translations and $\gr$ are oriented the induced function by taking determinants in the functoriality of the Atiyah class with respect to the morphism $a \circ \gr$ is just $1$.
\end{proof}
\begin{cor}
Let $X$ be in Example \ref{exp_assumptions_examples} except for coherent sheaves, then $\exp$ is oriented.
\end{cor}
\begin{proof}
    Having shown that $f_{\exp} = \mathrm{det}^{-1} j^{R, \et}_{\mathrm{at}_{X}}$ we can conclude that $f_{\exp}$ has a square root by Proposition \ref{square_root_j}. 
\end{proof}
\begin{lem} \label{translation_orientations}
    Let $X$ be a derived stack with an action of $\B T$ and a compatible exponential map as in Lemma \ref{center_symplectic_act}. Then fix $a \in \mathfrak{t}$ and write $A = \exp a$. We get the following equation on volume forms on formal completions
    \begin{equation}
    \exp^{*}\vol_{(x,A)} = (a^{-1*}j_{X,(x,0)}) \vol_{x,a}.
\end{equation}
Here $j_{X,(x,0)}$ is the restriction of $j_{X} \in \mathcal{O}(\mathcal{L}^{u}X)$ to $\mathcal{O}(\widehat{\To}^{(x,0)}[-1]X)$.
\end{lem}
\begin{proof}
We consider the diagram
    \begin{equation}
\begin{tikzcd}
	{\widehat{\To}^{(x,0)}[-1]X} & {\widehat{\mathcal{L}}^{(x, 1)}X} \\
	{\widehat{\To}^{(x,a)}[-1]X} & {\widehat{\mathcal{L}}^{(x, A)}X}
	\arrow["\exp"', from=1-1, to=1-2]
	\arrow["a", from=1-1, to=2-1]
	\arrow["A", from=1-2, to=2-2]
	\arrow["\exp", from=2-1, to=2-2]
\end{tikzcd}
    \end{equation}
Then we can compute 
\begin{equation}
    \exp^{*} \vol_{(x,A)}  = f_{\exp,(x,a)} \vol_{(x,a)}
\end{equation}
\begin{align*}
    (a^{*}f_{\exp,(x,a)}) \vol_{(x,0)} &= a^{*}(f \vol_{(x,a)}) \\
    & = a^{*} \exp^{*} \vol_{(x,A)} \\
    & = \exp^{*} A^{*} \vol_{(x,A)}  \\
    & = \exp^{*} \vol_{(x,1)} \quad \text{Lemma } \ref{volume_form_and_translation} \\
    & = j_{X,(x,0)} \vol_{(x,0)} \quad \text{Theorem } \ref{naef_saf}
\end{align*}
Therefore, applying $a^{-1*}$ we can conclude that 
\begin{equation}
    f_{\exp, (x,a)} = a^{-1*}j_{X,(x,0)}.
\end{equation}
\end{proof}

\section{Loop dimensional reduction and loop nonabelian Hodge}
We now turn to applying the exponential map to compute DT invariants via a loop dimensional reduction. In this section, using Theorem \ref{main_dcrit_thm}  and the constructions in Section \ref{exp_examples}, we will compute the BPS sheaves of the multiplicative moduli problems $\mathcal{L}X$ for the $\GL_n$ version of the moduli problems $X$ in \ref{exp_assumptions_examples}. For local systems and Higgs bundles, this will lead to the loop nonabelian Hodge theorem.
\subsection{Action of loops and shifted tangents}
\begin{notation} \label{dashed}
    In this section to reduce clutter we will write dashed arrows to mean maps that exist on the complex analytification of the truncations.
\end{notation}
Let $X$ be as in \ref{exp_assumptions_examples} of linear type. We have an action of $\B \mathbb{G}_{m}$ on $X$
\begin{equation}
    \B \mathbb{G}_{m} \times X \to X
\end{equation}
Taking $\Map(\B \widehat{\mathbb{G}}_{a},-)$ and $\Map(S^{1}_{\B},-)$ we get actions, that are compatible with the exponential.
\begin{equation}
\begin{tikzcd}
	{\mathbb{G}_{a} \times \B \mathbb{G}_{m} \times \To[-1]X} & {\mathbb{G}_{m} \times \B \mathbb{G}_{m} \times \mathcal{L}X} \\
	{\To[-1]X} & {\mathcal{L}X}
	\arrow["{\exp \times \id \times \exp}"', dashed, from=1-1, to=1-2]
	\arrow[from=1-1, to=2-1]
	\arrow[from=1-2, to=2-2]
	\arrow["\exp", dashed, from=2-1, to=2-2]
\end{tikzcd}
\end{equation}
In particular, we get induced maps
\begin{equation}
\begin{tikzcd}
	{\mathbb{G}_{a} \times X} & {\mathbb{G}_{a} \times M} \\
	{\mathbb{G}_{a} \times \To[-1]X} & {\mathbb{G}_{a} \times M_{a}} \\
	{\To[-1]X} & {M_{a}}
	\arrow[from=1-1, to=1-2]
	\arrow["{\id \times \mathrm{const}}", from=1-1, to=2-1]
	\arrow[from=1-2, to=2-2]
	\arrow[from=2-1, to=2-2]
	\arrow["{\mathrm{act}_{a}}", from=2-1, to=3-1]
	\arrow["{\iota_{a}}"', from=2-2, to=3-2]
	\arrow[from=3-1, to=3-2]
\end{tikzcd}
\end{equation}
We get the same type of map
\begin{equation}
    \mathbb{G}_{m} \times M \xrightarrow{\iota_{m}} M_{m}
\end{equation}
by replacing $\To[-1]$ with $\mathcal{L}$. 
\begin{prop}
We have that
\begin{align*}
    \BPS_{a} & \text{ is equivariant with respect to } \mathbb{G}_{a} \\
    \BPS_{m} & \text{ is equivariant with respect to } \mathbb{G}_{m}.
\end{align*}
\end{prop}
\begin{proof}
This follows from the fact that the $\mathbb{G}_{a}$ and $\mathbb{G}_{m}$ actions preserve the $(-1)$-shifted symplectic structures and are oriented as in Lemma \ref{center_symplectic_act}. Therefore, the action descends to the DT sheaves and therefore to the first perverse pieces of the pushforward to the good moduli spaces.
\end{proof}
\subsection{Loop dimensional reduction}
Recall notation \ref{dashed}. We have the following diagram of spaces 
\begin{equation} \label{exp_dimensional_support}
\begin{tikzcd}
	{(\To[-1]X)^{\et}} & {\To[-1]X} & {\mathcal{L}X} \\
	{M^{\et}_{a}} & {M_{a}} & {M_{m}} \\
	& {M \times \mathbb{G}_{a}} & {M \times \mathbb{G}_{m}}
	\arrow[dashed, from=1-1, to=1-2]
	\arrow["{\pi^{\et}_{a}}"{description}, dashed, from=1-1, to=2-1]
	\arrow["\exp", dashed, from=1-2, to=1-3]
	\arrow["{\pi_{a}}"', from=1-2, to=2-2]
	\arrow["{\pi_m}", from=1-3, to=2-3]
	\arrow[dashed, from=2-1, to=2-2]
	\arrow["\exp"', dashed, from=2-2, to=2-3]
	\arrow[dashed, from=3-2, to=2-1]
	\arrow["{\iota_{a}}"', from=3-2, to=2-2]
	\arrow["\exp", dashed, from=3-2, to=3-3]
	\arrow["{\iota_{m}}", from=3-3, to=2-3]
\end{tikzcd}
\end{equation}
 Now by
\cite[Theorem 7.2.15 ]{bu2025cohomology} there exists a perverse sheaf $\BPS_{M}$ on $M$ such that we have an isomorphism of perverse sheaves
\begin{equation} \label{support_lemma}
    \BPS_{M_{a}} \cong \iota_{a*}(\BPS_{M} \boxtimes \IC(\mathbb{G}_{a})).
\end{equation}
\begin{remark}
 Note that technically we have to make sure that there is a symplectic equivalence $\To^{*}[-1]X$ and $\To[-1]X$. In our cases this follows ultimately from the statement for $\B G$ or $\mathrm{Perf}$, where this can be checked explicitly, as our moduli spaces are obtained from various other constructions from these stacks. An approach to the general case would be to define a contraction operator for vector fields, along the lines of \cite[Definition 2.4]{brav2019darboux} and use a tautological $(-1)$-vector field on $\To[-1]X$ to contract the pullback of the closed $2$-form on $X$. In this way one should be able to construct a primitive for the AKSZ form on $\To[-1]X$, which one can compare with the Liouville form, pulled back from $\To^{*}[-1]X$.
\end{remark}
We now have the following theorem
\begin{thm}[Loop dimensional reduction] \label{loop_dim_red}
Let $X$ be the examples in subsection \ref{exp_assumptions_examples} of linear type with good moduli spaces as constructed in Section \ref{exp_examples}. Then we have an isomorphism of perverse sheaves
\begin{equation}
    \BPS_{M_m} \cong \iota_{m*}(\BPS_{M} \boxtimes \IC(\mathbb{G}_{m})).
\end{equation}
\end{thm}
\begin{proof}
We begin by computing the support of the BPS sheaf $\BPS_{M_{m}}$. Using Proposition \ref{etale_pullback_diagram}, Theorem \ref{main_dcrit_thm} and diagram \ref{exp_dimensional_support} we can write
\begin{equation}
    \exp^{*}\pi_{m*} \varphi_{m} \cong \pi^{\et}_{a*}\exp^{*}\varphi_{m} \cong \pi_{a*}|_{M^{\et}_{a}} \varphi^{\et}_{a}.
\end{equation}
Note that since the support $M \times \mathbb{G}_{a}$ lives inside the open $M^{\et}_{a}$ we also get 
\begin{equation}
    \BPS_{M^{\et}_{a}} = \BPS_{M_{a}}|_{M^{\et}_{a}} \cong \iota_{a*}(\BPS_{M} \boxtimes \IC(\mathbb{G}_{a})).
\end{equation}
Write $\exp^{\et}$ to be the restriction of $\exp$ to $M^{\et}_{a}$. Then since $\exp^{\et}$ is an \'etale map we get
\begin{equation}
    \exp^{*\et}\BPS_{M_m} \cong \BPS_{M^{\et}_{a}} 
\end{equation}
From which we can deduce that
\begin{equation}
    \BPS_{M_{m}} \cong \iota_{m*}(\BPS_{M} \boxtimes \mathcal{L})
\end{equation}
where $\mathcal{L}$ is a (shifted) local system on $\mathbb{G}_{m}$ of rank $1$. Now note that since $\mathbb{G}_{m}$ acts on $M \times \mathbb{G}_{m}$ by scaling the $\mathbb{G}_{m}$ factor and $\BPS_{M_{m}}$ is $\mathbb{G}_{m}$-equivariant we get that the local system $\mathcal{L}$ is $\mathbb{G}_{m}$-equivariant. Therefore, the local system $\mathcal{L}$ is the shifted constant sheaf. Note that for the case of coherent sheaves we have to apply a kind of sign twisted version of \cite[Theorem 7.2.15]{bu2025cohomology} since we cannot deduce that the usual orientation on $\To[-1] \mathrm{Coh}S$ is equivalent to the one pulled back by the exponential from $\mathcal{L}X$. However, the arguments of the theorem go through  so we can obtain a similar support lemma.
\end{proof}
\begin{remark}[Warning] \label{coh_warning}
    For all the cases apart from coherent sheaves in Example \ref{exp_assumptions_examples} we identify the $2d$-BPS sheaf with the $2d$-BPS sheaf studied in \cite{davison2023bps}. For coherent sheaves we cannot immediately deduce this due to the mismatch with orientations.
\end{remark}

\subsection{Nonabelian Hodge theory}
Let $g\geq 2$. We recall parts of the setup in \cite[Section 14]{davison2022bps}. 
Denote by $\mathrm{Higgs}^{\mathrm{ss}}_{g,r,d}$ the stack of semistable Higgs bundles of degree $d$ and rank $r$ and $\Loc^{\omega_{d}}_{\GL_{r}}(\Sigma_{g})$ the stack of twisted local systems, where $\omega_{d} = \exp(\frac{2 \pi d i}{r})$. Take
\begin{align} \label{coprime_ab}
    a \geq 1 ,b \in \mathbb{Z} \setminus 0 \text{ with } (a,b) =1 \text{ or } \\
    a= 1, b =0 \nonumber
\end{align}
For a slope $\theta = \frac{b}{a} \in \mathbb{Q}$, we can define
\begin{align*}
    \mathrm{BPS}^{\mathrm{B}}_{\theta} = \bigoplus_{n} \BPSo^{\B}_{g , na, nb} = \mathrm{Free}_{\mathrm{Lie}}(\bigoplus_{n \geq 0} \mathrm{IH}^{*}(\mathrm{M}^{\mathrm{B}}_{g,na,nb})) \\
    \mathrm{BPS}^{\mathrm{Dol}}_{\theta} = \bigoplus_{n} \BPSo^{\Dol}_{g , na, nb} = \mathrm{Free}_{\mathrm{Lie}}(\bigoplus_{n \geq 0} \mathrm{IH}^{*}(\mathrm{M}^{\Dol}_{g,na,nb}))
\end{align*}
where $\mathrm{M}^{\Dol}_{na,nb}$ is the good moduli space of the rank $na$ and degree $nb$ Higgs bundles and $\mathrm{M}^{\B}_{na,nb}$ is the good moduli space of rank $na$ and twist $nb$ local systems respectively. The free Lie algebra on intersection cohomology has a double grading and
\begin{equation}
    \BPSo^{\#}_{g,na,nb}
\end{equation}
are the degree $nb$ and rank $na$ parts. In \cite{davison2022bps} for $g \geq 2$  it is shown that the above cohomology is the cohomology of the $2d$ BPS sheaf of $\mathrm{Higgs}^{\mathrm{ss}}_{g,r,d}$ and $\Loc^{\omega_{d}}_{\GL_{r}}(\Sigma_{g})$  as in \cite[Theorem 7.2.15]{bu2025cohomology}. 
\subsection{Computations of BPS cohomology of $\Sigma_g \times S^{1}$}
\begin{prop} \label{sigma_g_bps}
Fix $g \geq 2$. Denote by $\BPSo^{\B,m}_{g,r,0}$ the BPS cohomology of $\Loc_{\GL_r}(\Sigma_{g} \times S^{1})$
We have an isomorphism
\begin{equation}
   \BPSo^{\B,m}_{g,r,0} =  \HHf(\mathrm{M}^{\mathrm{B}}_{g,r,1},\mathbb{Q}_{\mathrm{vir}}) \otimes \mathrm{H}^{*}(\mathbb{G}_{m}). 
\end{equation}
\end{prop}
\begin{proof}
By loop dimensional reduction we get that
\begin{equation}
    \BPSo^{\B,m}_{g,r,0} = \BPSo^{\B}_{g,r,0} \otimes \mathrm{H}^{*}(\mathbb{G}_{m}).
\end{equation}
    Now using classical nonabelian Hodge theory we have
    \begin{equation}
         \BPSo^{\B}_{g , r, 0} \cong \BPSo^{\mathrm{Dol}}_{g ,r,0}.
    \end{equation}
    Then by \cite[Corollary 5.15]{kinjokoseki2024cohomological} we have that 
    \begin{equation}
        \BPS^{\Dol}_{g, r,0} \cong \BPS^{\Dol}_{g,r,1}
    \end{equation}
    Now using the classical NAHT we have
\begin{equation}
    \BPSo^{\Dol}_{g,r,1} \cong \HHf(\mathrm{M}^{\mathrm{B}}_{g,r,1}, \mathbb{Q}_{\vir})
\end{equation}
This follows since the pair $(r,1)$ is coprime and the space is smooth. This is discussed in \cite[Example 5.18]{kinjokoseki2024cohomological}.
\end{proof}
To obtain this theorem we also could have just directly used the Betti $\chi$ independence in \cite[Theorem 14.10]{davison2022bps}.
\begin{remark}
Note that the cohomology of the twisted smooth character varieties has been computed in \cite{hausel_mixed}. The $g=1$ case has been computed in \cite{kaubrys2024cohomological}.
\end{remark}

\subsection{Loop nonabelian Hodge theory for $\GL_n$} \label{loop_naht_proof}
Take $a$,$b$ as in equation \eqref{coprime_ab}. We define
\begin{align*}
    \mathcal{L} \mathrm{Higgs}^{\mathrm{ss}}_{\theta} & = \coprod_{n \geq 0} \mathcal{L} \mathrm{Higgs}^{\mathrm{ss}}_{g,na,nb} \\
    \mathcal{L} \mathrm{Loc}_{\theta} & = \coprod_{n \geq 0} \mathcal{L} \mathrm{Loc}^{\omega_{nb}}_{\GL_{na}}(\Sigma_g) 
\end{align*}
similarly we write
\begin{align*}
   M^{\mathrm{Dol},m}_{\theta} & = \coprod_{n \geq 0}  M^{\mathrm{Dol},m}_{g,na,nb} \\
    M^{\mathrm{B},m}_{\theta} & = \coprod_{n \geq 0}  M^{\mathrm{B},m}_{g,na,nb} .
\end{align*}
Denote the respective good moduli space maps by $\pi_{\mathrm{Dol}}$ and $\pi_{\mathrm{B}}$.
Both of these moduli spaces have induced monoid structures given by direct sum. These induce a convolution monoidal structure $\boxtimes_{\oplus}$ on $\mathrm{D}^{+}_{c}(M^{\mathrm{Dol},m}_{\theta})$ and $\mathrm{D}^{+}_{c}(M^{\mathrm{B},m}_{\theta})$.
We now use the results on cohomological integrality for $3$-Calabi-Yau categories from \cite{bu2025cohomology}.
\begin{prop}[Cohomological integrality for fixed slope]
Let $g \geq 1$. Denote the DT sheaves on $\mathcal{L}\mathrm{Higgs}^{\mathrm{ss}}_{\theta}$ and $\mathcal{L}\mathrm{Loc}_{\theta}$ by $\varphi^{\theta}_{\mathrm{Dol}}$ and $\varphi^{\theta}_{\mathrm{B}}$ respectively. 
We have the following equivalences 
\begin{align} \label{coh_int_slope}
    \pi_{\mathrm{Dol}*} \varphi^{\theta}_{\mathrm{Dol}} \cong \mathrm{Sym}_{\boxtimes_{\oplus}}(\bigoplus_{n} \BPS^{\mathrm{Dol},m}_{g,nr,nd} \otimes \mathrm{H}^{*}(\B \mathbb{G}_{m})_{\mathrm{vir}}),  \text{ in } \mathrm{D}^{+}_{c}(M^{\mathrm{Dol},m}_{\theta}) \\
    \pi_{\mathrm{B}*} \varphi^{\theta}_{\mathrm{B}} \cong \mathrm{Sym}_{\boxtimes_{\oplus}}(\bigoplus_{n} \BPS^{\mathrm{B},m}_{g,nr,nd} \otimes \mathrm{H}^{*}(\B \mathbb{G}_{m})_{\mathrm{vir}}), \text{ in } \mathrm{D}^{+}_{c}(M^{\mathrm{B},m}_{\theta}).
\end{align}
From the above equivalences we then get equivalences in $\mathrm{Vect}^{\mathbb{N} \times \mathbb{Z}}$
\begin{align}\label{coh_int_slope_cohomology}
    \mathrm{H}^{*}(\mathcal{L} \mathrm{Higgs}^{\mathrm{ss}}_{\theta},\varphi^{\theta}_{\mathrm{Dol}}) & \cong \mathrm{Sym}(\bigoplus_{n} \BPSo^{\mathrm{Dol},m}_{g,nr,nd} \otimes \mathrm{H}^{*}(\B \mathbb{G}_{m})_{\mathrm{vir}}) \\
    \mathrm{H}^{*}(\mathcal{L} \mathrm{Loc}_{\theta},\varphi^{\theta}_{\mathrm{B}}) & \cong \mathrm{Sym}(\bigoplus_{n} \BPSo^{\mathrm{B},m}_{g,nr,nd} \otimes \mathrm{H}^{*}(\B \mathbb{G}_{m})_{\mathrm{vir}}).
\end{align}
\end{prop}
\begin{proof}
This follows from \cite[Section 10.2]{bu2025cohomology} after we show that $\mathrm{L}\mathrm{Higgs}^{\mathrm{ss}}_{\theta}$ and $\mathrm{L}\mathrm{Loc}_{\theta}$  satisfy the assumptions of a linear moduli stack with commutative orientation. The stacks $\mathcal{L}\mathrm{Higgs}^{\mathrm{ss}}_{\theta}$ and $\mathcal{L}\Loc_{\theta}$ are linear moduli stacks since we can view them as open substacks of the moduli of objects of the abelian categories of Higgs sheaves equipped with an automorphism and  the category of twisted fundamental group representations equipped with endomorphism. Furthermore, these abelian categories are themselves subcategories of the dg categories $\mathrm{Higgs} \otimes \mathrm{Loc}(S^{1})$ and  $\Loc(\Sigma_{G}) \otimes \mathrm{Loc}(S^{1})$. The symmetry condition on ext groups required follows from the ortogonality of loops stacks which is proved in \cite{bu2025cohomology}. Finally, the commutative orientation is constructed in subsection \ref{comm_orient_data}.
\end{proof}
Now fix a 
then we have the following theorem
\begin{thm}[Loop NAHT] \label{loop_naht} 
Let $g \geq 1$, $r \geq 1$ $d \in \mathbb{Z}$. Then we have an equivalence of BPS cohomology
\begin{equation}
     \BPSo^{\mathrm{Dol},m}_{g,r,d} \cong \BPSo^{\mathrm{B},m}_{g,r,d}.
\end{equation}
Furthermore, we have an isomorphism of graded vector spaces
\begin{equation}
    \mathrm{H}^{*}(\mathcal{L}\mathrm{Higgs}^{\mathrm{ss}}_{g,r,d}, \varphi_{\mathrm{Dol}}) \cong \mathrm{H}^{*}(\mathcal{L}\Loc^{\omega_{d}}_{\GL_r}(\Sigma_g), \varphi_{\mathrm{B}}).
\end{equation}
\end{thm}
\begin{proof}
Assume first that $g \geq 2$. We start by using Theorem \ref{loop_dim_red} for the stacks of Higgs bundles and local systems to deduce that 
\begin{equation}
    \BPS^{\mathrm{Dol},m}_{g,r,d} \cong \BPS^{\mathrm{Dol}}_{g,r,d} \boxtimes \mathrm{IC}(\mathbb{G}_{m}) ,\quad \BPS^{\mathrm{B},m}_{g,r,d} \cong \BPS^{\mathrm{B}}_{g,r,d} \boxtimes \mathrm{IC}(\mathbb{G}_{m})
\end{equation}
Then we have by classical nonabelian Hodge that
\begin{equation}
    \BPSo^{\mathrm{B}}_{g,r,d}  \cong \BPSo^{\mathrm{Dol}}_{g,r,d}
\end{equation}
therefore also 
\begin{equation}
     \BPSo^{\mathrm{Dol},m}_{g,r,d} \cong \BPSo^{\mathrm{B},m}_{g,r,d}
\end{equation}
proving the first part. Now using equation \eqref{coh_int_slope_cohomology} we get an isomorphism 
\begin{equation}
    \mathrm{H}^{*}(\mathcal{L} \mathrm{Higgs}^{\mathrm{ss}}_{\theta},\varphi^{\theta}_{\mathrm{Dol}}) \cong \mathrm{H}^{*}(\mathcal{L} \mathrm{Loc}_{\theta},\varphi^{\theta}_{\mathrm{B}})
\end{equation}
restricting to each dimension component we get an equivalence of $\mathbb{Z}$-graded vector-spaces
\begin{equation}
    \mathrm{H}^{*}(\mathcal{L} \mathrm{Higgs}^{\mathrm{ss}}_{na,nb},\varphi^{nr,nd}_{\mathrm{Dol}}) \cong \mathrm{H}^{*}(\mathcal{L} \mathrm{Loc}^{\omega^{nb}}_{\GL_{na}},\varphi^{na,nb}_{\mathrm{B}}).
\end{equation}
for some integers $a$,$b$, as  in \ref{coprime_ab}, such that $r = na$ and $d = nb$. Now if $g=1$, then the $2d$ BPS cohomology is computed in \cite{davison2023nonabelian} and \cite[Section 5.5]{davison2022integrality}. Then we can use the same arguments to deduce the result. 
\end{proof}
\begin{remark}
    Note that this theorem is weaker than the corresponding $2$-dimesnional one \cite[Theorem 14.4]{davison2022bps}. This is because we have not given a homeomorphism between $M^{\mathrm{Dol},m}_{a,b}$ and $M^{\mathrm{B},m}_{a,b}$. However we believe this fact can be proven using the tools of the usual nonabelian Hodge correspondence.
\end{remark}

\printbibliography

@article {ben2015darboux,
    AUTHOR = {Ben-Bassat, Oren and Brav, Christopher and Bussi, Vittoria and
              Joyce, Dominic},
     TITLE = {A `{D}arboux theorem' for shifted symplectic structures on
              derived {A}rtin stacks, with applications},
   JOURNAL = {Geom. Topol.},
  FJOURNAL = {Geometry \& Topology},
    VOLUME = {19},
      YEAR = {2015},
    NUMBER = {3},
     PAGES = {1287--1359},
      ISSN = {1465-3060},
   MRCLASS = {14A20 (14D23 14F05 14N35 32S30)},
  MRNUMBER = {3352237},
MRREVIEWER = {Hsian-Hua Tseng},
       DOI = {10.2140/gt.2015.19.1287},
       URL = {https://doi.org/10.2140/gt.2015.19.1287},
}

@Article{Joyce_dcrit,
 Author = {Joyce, Dominic},
 Title = {A classical model for derived critical loci},
 FJournal = {Journal of Differential Geometry},
 Journal = {J. Differ. Geom.},
 ISSN = {0022-040X},
 Volume = {101},
 Number = {2},
 Pages = {289--367},
 Year = {2015},
 Language = {English},
 DOI = {10.4310/jdg/1442364653},
 zbMATH = {6507799},
 Zbl = {1368.14027}
}

@article {PTVV,
    AUTHOR = {Pantev, Tony and To\"{e}n, Bertrand and Vaqui\'{e}, Michel and
              Vezzosi, Gabriele},
     TITLE = {Shifted symplectic structures},
   JOURNAL = {Publ. Math. Inst. Hautes \'{E}tudes Sci.},
  FJOURNAL = {Publications Math\'{e}matiques. Institut de Hautes \'{E}tudes
              Scientifiques},
    VOLUME = {117},
      YEAR = {2013},
     PAGES = {271--328},
      ISSN = {0073-8301},
   MRCLASS = {14F05 (14A15 18F20 18G30 53D05 53D12)},
  MRNUMBER = {3090262},
MRREVIEWER = {Andrey Yu. Lazarev},
       DOI = {10.1007/s10240-013-0054-1},
       URL = {https://doi.org/10.1007/s10240-013-0054-1},
}

@article{brav2019darboux,
  title={A Darboux theorem for derived schemes with shifted symplectic structure},
  author={Brav, Christopher and Bussi, Vittoria and Joyce, Dominic},
  journal={Journal of the American Mathematical Society},
  volume={32},
  number={2},
  pages={399--443},
  year={2019}
}

@article{ben2012loop,
  title={Loop spaces and connections},
  author={Ben-Zvi, David and Nadler, David},
  journal={Journal of Topology},
  volume={5},
  number={2},
  pages={377--430},
  year={2012},
  publisher={Oxford University Press}
}

@book{achar2021perverse,
  title={Perverse sheaves and applications to representation theory},
  author={Achar, Pramod N},
  volume={258},
  year={2021},
  publisher={American Mathematical Soc.}
}

@misc{naef2023torsion,
      title={Torsion volume forms}, 
      author={Florian Naef and Pavel Safronov},
      year={2023},
      eprint={2308.08369},
      archivePrefix={arXiv},
      primaryClass={math.AG}
}

@article{calaque2022aksz,
  title={The AKSZ construction in derived algebraic geometry as an extended topological field theory},
  author={Calaque, Damien and Haugseng, Rune and Scheimbauer, Claudia},
  journal={arXiv preprint arXiv:2108.02473},
  year={2021}
}

@article{calaque_cot,
   title={Shifted cotangent stacks are shifted symplectic},
   volume={28},
   ISSN={2258-7519},
   url={http://dx.doi.org/10.5802/afst.1593},
   DOI={10.5802/afst.1593},
   number={1},
   journal={Annales de la Facult \' e des sciences de Toulouse: Mathématiques},
   publisher={Cellule MathDoc/CEDRAM},
   author={Calaque, Damien},
   year={2019},
   month=mar, pages={67–90} }

@misc{ricolfi_savvas,
      title={The d-critical structure on the Quot scheme of points of a Calabi-Yau 3-fold}, 
      author={Andrea T. Ricolfi and Michail Savvas},
      year={2023},
      eprint={2106.16133},
      archivePrefix={arXiv},
      primaryClass={math.AG},
      url={https://arxiv.org/abs/2106.16133}, 
}

@Article{davison2022integrality,
 Author = {Davison, Ben},
 Title = {The integrality conjecture and the cohomology of preprojective stacks},
 FJournal = {Journal f{\"u}r die Reine und Angewandte Mathematik},
 Journal = {J. Reine Angew. Math.},
 ISSN = {0075-4102},
 Volume = {804},
 Pages = {105--154},
 Year = {2023},
 Language = {English},
 DOI = {10.1515/crelle-2023-0065},
 zbMATH = {7761116},
 Zbl = {1535.14005}
}

@Article{hausel_mixed,
 Author = {Hausel, Tam{\'a}s and Rodriguez-Villegas, Fernando},
 Title = {Mixed {Hodge} polynomials of character varieties. {With} an appendix by {Nicholas} {M}. {Katz}.},
 FJournal = {Inventiones Mathematicae},
 Journal = {Invent. Math.},
 ISSN = {0020-9910},
 Volume = {174},
 Number = {3},
 Pages = {555--624},
 Year = {2008},
 Language = {English},
 DOI = {10.1007/s00222-008-0142-x},
 zbMATH = {5530512},
 Zbl = {1213.14020}
}

@article{gunningham2023deformation,
  title={Deformation quantization and perverse sheaves},
  author={Gunningham, Sam and Safronov, Pavel},
  journal={arXiv preprint arXiv:2312.07595},
  year={2023}
}

@Article{joyce_conje_lino_ben,
 Author = {Amorim, Lino and Ben-Bassat, Oren},
 Title = {Perversely categorified {Lagrangian} correspondences},
 FJournal = {Advances in Theoretical and Mathematical Physics},
 Journal = {Adv. Theor. Math. Phys.},
 ISSN = {1095-0761},
 Volume = {21},
 Number = {2},
 Pages = {289--381},
 Year = {2017},
 Language = {English},
 DOI = {10.4310/ATMP.2017.v21.n2.a1},
 zbMATH = {6826483},
 Zbl = {1394.53079}
}

@Article{thomas_dt_og,
 Author = {Thomas, R. P.},
 Title = {A holomorphic {Casson} invariant for {Calabi}-{Yau} 3-folds, and bundles on {{\(K3\)}} fibrations.},
 FJournal = {Journal of Differential Geometry},
 Journal = {J. Differ. Geom.},
 ISSN = {0022-040X},
 Volume = {54},
 Number = {2},
 Pages = {367--438},
 Year = {2000},
 Language = {English},
 DOI = {10.4310/jdg/1214341649},
 zbMATH = {1782643},
 Zbl = {1034.14015}
}

@Article{cptvv,
 Author = {Calaque, Damien and Pantev, Tony and To{\"e}n, Bertrand and Vaqui{\'e}, Michel and Vezzosi, Gabriele},
 Title = {Shifted {Poisson} structures and deformation quantization},
 FJournal = {Journal of Topology},
 Journal = {J. Topol.},
 ISSN = {1753-8416},
 Volume = {10},
 Number = {2},
 Pages = {483--584},
 Year = {2017},
 Language = {English},
 DOI = {10.1112/topo.12012},
 zbMATH = {6778841},
 Zbl = {1428.14006}
}

@Article{dim_red_kinjo,
 Author = {Kinjo, Tasuki},
 Title = {Dimensional reduction in cohomological {Donaldson}-{Thomas} theory},
 FJournal = {Compositio Mathematica},
 Journal = {Compos. Math.},
 ISSN = {0010-437X},
 Volume = {158},
 Number = {1},
 Pages = {123--167},
 Year = {2022},
 Language = {English},
 DOI = {10.1112/S0010437X21007740},
 zbMATH = {7481800},
 Zbl = {1493.14096}
}

@Article{sun_analytic,
 Author = {Sun, Shenghao},
 Title = {Generic base change, {Artin}'s comparison theorem, and the decomposition theorem for complex {Artin} stacks},
 FJournal = {Journal of Algebraic Geometry},
 Journal = {J. Algebr. Geom.},
 ISSN = {1056-3911},
 Volume = {26},
 Number = {3},
 Pages = {513--555},
 Year = {2017},
 Language = {English},
 DOI = {10.1090/jag/683},
 zbMATH = {6707825},
 Zbl = {1453.14005}
}

@article{AHR_luna,
   title={A Luna étale slice theorem for algebraic stacks},
   volume={191},
   ISSN={0003-486X},
   url={http://dx.doi.org/10.4007/annals.2020.191.3.1},
   DOI={10.4007/annals.2020.191.3.1},
   number={3},
   journal={Annals of Mathematics},
   publisher={Annals of Mathematics},
   author={Alper, Jarod and Hall, Jack and Rydh, David},
   year={2020},
   month=may }

@article{dg_indschemes,
  title={DG indschemes},
  author={Gaitsgory, Dennis and Rozenblyum, Nick},
  journal={Perspectives in representation theory},
  volume={610},
  pages={139--251},
  year={2014},
eprint={1108.1738},
      archivePrefix={arXiv},
}

@article{kinjo2024cohomological,
  title={Cohomological Hall algebras for 3-Calabi-Yau categories},
  author={Kinjo, Tasuki and Park, Hyeonjun and Safronov, Pavel},
  journal={arXiv preprint arXiv:2406.12838},
  year={2024}
}

@article{cal_saf,
  title={Shifted cotangent bundles, symplectic groupoids and deformation to the normal cone},
  author={Calaque, Damien and Safronov, Pavel},
  journal={arXiv preprint arXiv:2407.08622},
  year={2024}
}

@article{monier2021notelinearstacks,
  title={A note on linear stacks},
  author={Monier, Ludovic},
  journal={arXiv preprint arXiv:2103.06555},
  year={2021}
}

@article{gaitsgory2017study,
  title={A study in derived algebraic geometry: Volume ii: Deformations, lie theory and formal geometry},
  author={Gaitsgory, Dennis and Rozenblyum, Nick},
  journal={American Mathematical Society},
  year={2017}
}

@article{de2009decomposition,
  title={The decomposition theorem, perverse sheaves and the topology of algebraic maps},
  author={de Cataldo, Mark and Migliorini, Luca},
  journal={Bulletin of the American Mathematical Society},
  volume={46},
  number={4},
  pages={535--633},
  year={2009}
}

@book{knapp1988lie,
  title={Lie groups, Lie algebras, and cohomology},
  author={Knapp, Anthony W},
  volume={34},
  year={1988},
  publisher={Princeton University Press}
}

@Inbook{Park2025,
author="Park, Hyeonjun
and You, Jemin",
editor="Kiem, Young-Hoon",
title="An Introduction to Shifted Symplectic Structures",
bookTitle="Moduli Spaces, Virtual Invariants and Shifted Symplectic Structures",
year="2025",
publisher="Springer Nature Singapore",
address="Singapore",
pages="37--64",
isbn="978-981-97-8249-9",
doi="10.1007/978-981-97-8249-9_2",
url="https://doi.org/10.1007/978-981-97-8249-9_2"
}

@article{hennion2024gluing,
  title={Gluing invariants of Donaldson--Thomas type--Part I: the Darboux stack},
  author={Hennion, Benjamin and Holstein, Julian and Robalo, Marco},
  journal={arXiv preprint arXiv:2407.08471},
  year={2024}
}

@article{holstein2024koszul,
  title={Koszul duality and Calabi-Yau structures},
  author={Holstein, Julian and Rivera, Manuel},
  journal={arXiv preprint arXiv:2410.03604},
  year={2024}
}

@article{holstein2018analytification,
  title={Analytification of mapping stacks},
  author={Holstein, Julian and Porta, Mauro},
  journal={arXiv preprint arXiv:1812.09300},
  year={2018}
}

@phdthesis{chen2018localization,
  title={A Localization Theorem for Derived Loop Spaces and Periodic Cyclic Homology},
  author={Chen, Harrison I and others},
  year={2018},
  school={UC Berkeley}
}

@article{kaubrys2024cohomological,
  title={Cohomological Donaldson-Thomas theory for local systems on the $3 $-torus},
  author={Kaubrys, Sarunas},
  journal={arXiv preprint arXiv:2409.16013},
  year={2024}
}

@misc{halpernleistner2022,
      title={On the structure of instability in moduli theory}, 
      author={Daniel Halpern-Leistner},
      year={2022},
      eprint={1411.0627},
      archivePrefix={arXiv},
      primaryClass={math.AG},
      url={https://arxiv.org/abs/1411.0627}, 
}

@article{davison2022bps,
  title={BPS Lie algebras for totally negative 2-Calabi-Yau categories and nonabelian Hodge theory for stacks},
  author={Davison, Ben and Hennecart, Lucien and Mejia, Sebastian Schlegel},
  journal={arXiv preprint arXiv:2212.07668},
  year={2022}
}

@article{alper2012local,
  title={Local properties of good moduli spaces},
  author={Alper, Jarod},
  journal={Tohoku Mathematical Journal, Second Series},
  volume={64},
  number={1},
  pages={105--123},
  year={2012},
  publisher={Mathematical Institute, Tohoku University}
}

@article{toda2018moduli,
  title={Moduli stacks of semistable sheaves and representations of Ext--quivers},
  author={Toda, Yukinobu},
  journal={Geometry \& Topology},
  volume={22},
  number={5},
  pages={3083--3144},
  year={2018},
  publisher={Mathematical Sciences Publishers}
}

@article{simpson1994moduli,
  title={Moduli of representations of the fundamental group of a smooth projective variety I},
  author={Simpson, Carlos T},
  journal={Publications Math{\'e}matiques de l'IH{\'E}S},
  volume={79},
  pages={47--129},
  year={1994}
}

@book{meinrenken2013clifford,
  title={Clifford algebras and Lie theory},
  author={Meinrenken, Eckhard},
  volume={58},
  year={2013},
  publisher={Springer}
}

@article{bu2025cohomology,
  title={Cohomology of symmetric stacks},
  author={Bu, Chenjing and Davison, Ben and N{\'u}nez, Andr{\'e}s Ib{\'a}nez and Kinjo, Tasuki and P{\u{a}}durariu, Tudor},
  journal={arXiv preprint arXiv:2502.04253},
  year={2025}
}

@article{descombes2025hyperbolic,
  title={Hyperbolic localization in Donaldson-Thomas theory},
  author={Descombes, Pierre},
  journal={arXiv preprint arXiv:2506.22400},
  year={2025}
}

@article{scherotzke2025fourier,
  title={Fourier--Mukai equivalences for formal groups and elliptic Hochschild homology},
  author={Scherotzke, Sarah and Sibilla, Nicol{\`o} and Tomasini, Paolo},
  journal={arXiv preprint arXiv:2505.00172},
  year={2025}
}

@inproceedings{kuhn2024atiyah,
  title={The Atiyah class on algebraic stacks},
  author={Kuhn, Nikolas},
  booktitle={Forum of Mathematics, Sigma},
  volume={12},
  pages={e100},
  year={2024},
  organization={Cambridge University Press}
}

@article{kinjokoseki2024cohomological,
  title={Cohomological chi-independence for Higgs bundles and Gopakumar--Vafa invariants},
  author={Kinjo, Tasuki and Koseki, Naoki},
  journal={Journal of the European Mathematical Society},
  year={2024}
}

@article{davison2023nonabelian,
  title={Nonabelian Hodge theory for stacks and a stacky P= W conjecture},
  author={Davison, Ben},
  journal={Advances in Mathematics},
  volume={415},
  pages={108889},
  year={2023},
  publisher={Elsevier}
}

@article{tanaka2017vafa,
  title={Vafa-Witten invariants for projective surfaces I: stable case},
  author={Tanaka, Yuuji and Thomas, Richard P},
  journal={arXiv preprint arXiv:1702.08487},
  year={2017}
}

@article{davison2023bps,
  title={BPS algebras and generalised Kac-Moody algebras from 2-Calabi-Yau categories},
  author={Davison, Ben and Hennecart, Lucien and Mejia, Sebastian Schlegel},
  journal={arXiv preprint arXiv:2303.12592},
  year={2023}
}

@article{kinjo_mult,
  title={Multiplicative dimensional reduction},
  author={Kinjo, Tasuki},
  year={2025}
}

\end{document}